\documentclass[10pt,fleqn,a4paper]{article}

\pdfoutput=1

\usepackage{graphicx}
\usepackage{latexsym}
\usepackage[small,bf]{caption}
\usepackage{amsfonts}
\usepackage{amssymb}
\usepackage{amsmath}
\usepackage{stmaryrd}
\usepackage{bbm}
\usepackage{theorem}
\usepackage{times}
\usepackage{color}

\usepackage{hyperref}

\newcommand{\be}{\begin{equation}}
\newcommand{\ee}{\end{equation}}
\newcommand{\bea}{\begin{eqnarray}}
\newcommand{\eea}{\end{eqnarray}}
\newcommand{\bef}{\begin{figure}}
\newcommand{\enf}{\end{figure}}
\newcommand{\ball}{\begin{array}{ll}}
\newcommand{\bacl}{\begin{array}{cl}}
\newcommand{\bal}{\begin{array}{l}}
\newcommand{\bac}{\begin{array}{c}}
\newcommand{\ea}{\end{array}}
\newcommand{\eps}{{\varepsilon}}
\newcommand{\1}{{\mathbbm{1}}}
\newcommand{\Q}{{\mathbb{Q}}}
\newcommand{\N}{{\mathbb{N}}}
\newcommand{\R}{{\mathbb{R}}}
\newcommand{\Z}{{\mathbb{Z}}}
\newcommand{\E}{{\mathbb{E}}}
\newcommand{\T}{{\mathbb{T}}}
\newcommand{\var}{{\mathrm{Var}}}
\renewcommand{\P}{{\mathbb{P}}}
\newcommand{\ca}{{\rm cap}}
\newcommand{\Lcal}{{\mathcal{L}}}
\newcommand{\Rcal}{{\mathcal{R}}}
\newcommand{\Ecal}{{\mathcal{E}}}
\newcommand{\Dcal}{{\mathcal{D}}}

\newcommand{\Ccal}{{\mathcal{C}}}
\newcommand{\Ical}{{\mathcal{I}}}
\newcommand{\Scal}{{\mathcal{S}}}
\newcommand{\Tcal}{{\mathcal{T}}}

\newcommand{\tx}{{\tilde x}}
\newcommand{\kbo}{{\mathbf{k}}}

\newcommand{\tr}{t_{\mathrm{rel}}}
\newcommand{\tm}{t_{\mathrm{mix}}}
\newcommand{\epsl}{{\underline{\epsilon}_L}}
\newcommand{\epsu}{{\overline{\epsilon}_L}}

\def\CQFD{\hfill\hbox{\vrule\vbox to 8pt{\hrule width 8pt\vfill\hrule}\vrule}}

\def\tN{\tilde{N}}
\def\teta{\tilde \eta}

\def\a{\alpha}
\def\tita{\tilde \eta}

\numberwithin{equation}{section}
\theoremstyle{plain}
\newtheorem{theorem}{Theorem}[section]
\newtheorem{lemma}[theorem]{Lemma}
\newtheorem{proposition}[theorem]{Proposition}
\newtheorem{corollary}[theorem]{Corollary}

\title{Metastability in a condensing zero-range process in the thermodynamic limit}
\author{In\'es Armend\'ariz\thanks{Departamento de Matem\'atica, Universidad de Buenos Aires, C1428EGA Buenos Aires, Argentina.\newline Email:\,{\tt iarmend@dm.uba.ar}}, Stefan Grosskinsky\thanks{Mathematics Institute, Zeeman Building, University of Warwick, Coventry CV4 7AL, UK.\newline Email:\,{\tt s.w.grosskinsky@warwick.ac.uk}}, Michail Loulakis\thanks{School of Applied Mathematical and Physical Sciences, National Technical University of Athens, 15780 Athens, Greece.\newline Email:\,{\tt loulakis@math.ntua.gr}} }
\date{\vspace{-5ex}}
\begin{document}
\sloppy
\maketitle

\begin{abstract}
\noindent
Zero-range processes with decreasing jump rates are known to exhibit condensation, where a finite fraction of all particles concentrates on a single lattice site when the total density exceeds a critical value. 
We study such a process on a one-dimensional lattice with periodic boundary conditions in the thermodynamic limit with fixed, super-critical particle density. We show that the process exhibits metastability with respect to the condensate location, i.e. the suitably accelerated process of the rescaled location converges to a limiting Markov process on the unit torus. This process has stationary, independent increments and the rates are characterized by the scaling limit of capacities of a single random walker on the lattice. 
Our result extends previous work for fixed lattices and diverging density in [J.~Beltran, C.~Landim, Probab.\ Theory Related Fields, \textbf{152}(3-4):781-807, 2012], and we follow the martingale approach developed there and in subsequent publications. 
Besides additional technical difficulties in estimating error bounds for transition rates, the thermodynamic limit requires new estimates for equilibration towards a suitably defined distribution in metastable wells, corresponding to a typical set of configurations with a particular condensate location. The total exit rates from individual wells turn out to diverge in the limit, which requires an intermediate regularization step using the symmetries of the process and the regularity of the limit generator. 
Another important novel contribution is a coupling construction to provide a uniform bound on the exit rates from metastable wells, which is of a general nature and can be adapted to other models.\\
\\
{\em AMS 2010 Mathematics Subject Classification}: 60K35, 82C22; 60J45, 82C26, 60J75\\
\\
{\em Keywords}: Metastability, Zero Range Process, Condensation.
\end{abstract}

\clearpage
\linespread{1.5}\selectfont

\tableofcontents

\clearpage
\linespread{1}\selectfont

\section{Introduction\label{sec:intro}}

A rigorous understanding of metastability phenomena in stochastic particle or spin systems has been a subject of major recent research. Intuitively, in such systems, the configuration space contains two or more disjoint metastable sets (called wells in the following) with an associated separation of time scales phenomenon in the scaling limit of a system parameter. In each well, the process spends a very long time which allows it to equilibrate to a metastable state. Exits from wells are then triggered by rare fluctuations, which lead to exponentially distributed waiting times in a well. Once activated, transitions to other wells occur on a much faster time scale, and do not depend on the detailed past of the sample path. So the limiting metastable dynamics corresponds to an effective Markov process on a highly reduced state space of metastable states associated to the wells. On a mathematically rigorous level, important conceptual questions of current research interest include a proper definition of metastable states in terms of probability measures, as well as a general practical framework to establish the separation of time scales. An important aspect in this context is the precise definition of metastable wells and in particular an optimal choice of their size, which is to some extent arbitrary. Intuitively, maximizing the depth and minimizing the complexity, or size, of the wells leads to a stronger separation of time scales. 
In combination they effectively characterize the free energy landscape for the chosen definition of wells in analogy to the classical framework of competition between energy (depth) and entropy (complexity).\\

Different approaches to metastability are summarized in \cite{olivierivares}, Chapter 4, including a pathwise treatment \cite{cassandroetal84} which is based on an analysis of empirical averages along typical trajectories, and has more recently been studied also in \cite{fernandezetal14}. 
For reversible systems a powerful potential theoretic approach has been developed \cite{bovieretal01,bovieretal02}, using systematically the concept of capacities to establish sharp estimates on expected transition times between metastable states. This technique was applied in various models 
and is summarized in the new monograph \cite{bovierbook} and the review papers \cite{bovier09,denhollander04}. 
Potential theoretic methods have been applied to a particular family of condensing zero-range processes in \cite{beltranetal09,beltranetal09b}, leading to the development of a martingale approach summarized in \cite{beltranetal14} which we follow in this paper. Instead of deriving exponential limit distributions of individual exit rates from metastable wells directly, the limit process is identified as a Markov process via the solution to a martingale problem.\\

The zero-range process was introduced in \cite{spitzer70} as a stochastic particle system without restriction on the local occupation numbers, where the jump rates of particles depend only on the occupation number of the departure site. The dynamics locally conserves the number of particles and the process is known to have a family of stationary product measures \cite{andjel82}, which can be indexed by the particle density $\rho$. Under certain conditions on the jump rates this family has a maximal element at a finite (critical) density $\rho_c$, and the process exhibits a condensation transition. Conditioned on densities $\rho >\rho_c$ the system phase separates into a homogeneous phase distributed at the critical product measure, and a condensate, where the remaining mass concentrates in a single lattice site for typical stationary configurations of large, finite systems. This phenomenon was first reported in the theoretical physics literature \cite{drouffe98,evans00}, and rigorous results including the equivalence of ensembles, a law of large numbers and a central limit theorem for the condensate size have been established in \cite{jeonetal00,grosskinskyetal03,armendarizetal08}. 
In spatially inhomogeneous zero-range processes a condensation phenomenon can occur that has different dynamic and static features as discussed above, and can be studied with the help of coupling techniques (see e.g. \cite{bahadoranetal14} and references therein). In this paper we focus on spatially homogeneous systems, which are necessarily non-monotone \cite{raffertyetal15} and where basic coupling techniques cannot be applied. \\

On a translation invariant lattice of size $L$ the condensed state exhibits an $L$-fold degeneracy with respect to the location of the condensate, which is uniformly distributed under the stationary measure. 
The metastable limit dynamics of the condensate location for reversible processes has been established in \cite{beltranetal09,beltranetal09b} in a slightly more general setting allowing for spatial inhomogeneities, but these results are restricted to a fixed lattice in the limit of diverging particle number $N$. In this setting, the depth of metastable wells dominates their complexity and they are relatively small sets. Repeated visits to particular configurations occur before the process exits the well (often called attractor states), which provides a renewal structure that can be used to establish the separation of time scales. The limit process of the condensate location is then a continuous-time random walk on the fixed lattice, where transition rates are proportional to the capacities of a single random walker. For the same model, matching upper and lower bounds for exit rates from metastable wells have been derived in \cite{rebecca} in the limit $L,N\to\infty$ with diverging particle density $N/L\to\infty$ (see also \cite{bovierbook}, Chapter 21). \\

In this paper we establish the full metastable limit dynamics for the condensate location for the above family of reversible zero-range processes in the thermodynamic limit, i.e. we take $L,N\to\infty$ with finite, supercritical density $N/L\to\rho >\rho_c$. We consider a one-dimensional system with periodic boundary conditions, which leads to a L\'evy-type limit process for the rescaled condensate location on the unit torus with stationary, independent increments. A major new difficulty is related to thermalization within the metastable wells, which are exponentially large in the system size and previously applied renewal techniques do not apply. To separate the slow and fast variables effectively, we describe the thermalization dynamics on the wells using the restricted process introduced in \cite{bianchietal15} and \cite{beltranetal14}. 
We can then prove relaxation time bounds for this process by a comparison argument with independent birth-death chains. 
Another key ingredient is a uniform bound on exit rates from a well which is established by a novel coupling construction and is used in several occasions in our proof, including thermalization. 
The characterization of the limiting generator of the process is more difficult in the thermodynamic limit due to increasing complexity of the free energy landscape, caused by the diverging number of wells and transition paths between them. To tackle this we have to introduce an intermediate regularization step on a coarse-grained lattice using the full symmetries of the model and regularity of the limit process. 
We note that, similarly to the results in \cite{rebecca} for diverging particle density, we are not able to establish matching bounds for transition rates between individual wells. But in order to identify the limit process it is sufficient to get matching bounds only between regularized sets of wells. 
As in previous results on this model, the leading order of transition rates between wells is polynomial in the system size $L$, which requires sharp estimates on transition rates from potential theory. Due to the increasing complexity of transition paths we need a better control than results in \cite{beltranetal09,beltranetal09b} on the leading order error terms, which constitutes an additional technical difficulty. We note that a rescaling argument in a different context has also been used in \cite{misturini} to establish metastability for the ABC model. \\

So far, there are only few metastability results that deal with infinite volume limits as summarized in \cite{bovierbook}, Part VII. Examples include kinetic Ising models at low temperatures \cite{bovieretal10} or low magnetic fields \cite{schonshlos}, results on the dynamics of critical droplets \cite{gois}, or 
dilute gases at vanishing density and temperature \cite{gaudilliereetal09}. 
All these results feature the additional scaling of a system parameter such as temperature or density, which increases the depth of metastable wells compared to their complexity. Often this is a necessary requirement for metastability to occur at all, and in \cite{rebecca} this has been used as a simplifying assumption to obtain results for the zero-range process in an infinite density limit. To our knowledge, the thermodynamic limit result derived here without any scaling of system parameters constitutes therefore one of the first metastability results of this kind.\\

The paper is organized as follows. In Section \ref{sec:not+main} we introduce notation and state our main result. The proof and a more technical discussion of the result is given in Section \ref{sec:proof_main}. The remaining sections are devoted to establishing the main ingredients of the proof, including the equilibration dynamics in wells in Section \ref{sec:equi}, uniform bounds on exit rates in Section \ref{sec:coup}, and the derivation of the generator of the limiting process in Section \ref{sec:regul}. Upper and lower bounds on expected transition rates from capacity estimates are given in Section \ref{sec:cap}, and Sections \ref{tight} and \ref{invmeas} contain auxiliary results on tightness, martingale convergence and properties of the stationary measure.

\section{Notation and main result \label{sec:not+main}}

\subsection{Notation\label{notation}}

Consider the zero-range process $\big(\eta (t):t\geq 0\big)$ on the one-dimensional discrete torus $\Lambda=\Z /L\Z$ with $N$ particles and state space 
\be
X_{L,N} =\big\{ \eta \in\N_0^{\Lambda} :S_L (\eta )=N\big\} \quad\mbox{where we denote}\quad S_L (\eta )=\sum_{x\in\Lambda} \eta_x \ .
\ee
The dynamics is defined by the generator
\be\label{zrgen}
\Lcal f(\eta )=\sum_{x,y\in\Lambda} p(x,y)\, g(\eta_x )\big( f(\eta^{x,y})-f(\eta )\big)
\ee
for all continuous functions $f:X_{L,N} \to \R$, with the usual notation $\eta^{x,y}_z =\eta_z -\delta_{z,x} +\delta_{z,y}$, $z\in\Lambda$ for the configuration $\eta^{x,y}$ where one particle has moved from site $x$ to site $y$. We focus on symmetric, nearest neighbour jump probabilities $p(x,y)=\frac12\delta_{y,x-1} +\frac12\delta_{y,x+1}$ with periodic boundary conditions, and jump rates are of the form
\begin{equation}\label{rates}
g(0)=0\ ,\quad g(1)=1\quad\mbox{and}\quad g(n)=\Big(\frac{n}{n-1}\Big)^b \mbox{ for }n\geq 2\ , \quad\mbox{with parameter }b>0\ .
\end{equation}
Without the canonical constraint $S_L =N$, the process is known to exhibit a family of stationary product measures with a maximal element $\nu$ \cite{andjel82,evans00,jeonetal00,grosskinskyetal03}, which has marginals
\be\label{tail}
\nu [\eta_x =n]=\frac{1}{z_c} \frac{1}{g!(n)}\quad\mbox{with}\quad g!(n)=\prod_{k=1}^n g(k)=n^{b}\quad\mbox{for }n\geq 1\quad\mbox{and}\quad g!(0)=1\ .
\ee
As long as $b>2$ the normalization and first moment of $\nu$ are both finite, i.e.
\be\label{zcrc}
z_c :=1+\sum_{n=1}^\infty n^{-b} <\infty\quad\mbox{and}\quad \rho_c :=\nu (\eta_x )=\frac{1}{z_c}\sum_{n=1}^\infty n^{1-b} <\infty\ .
\ee
The process is irreducible on $X_{L,N}$ and the corresponding unique stationary measures $\mu_{L,N}$ are called canonical distributions, and can be written as conditional product measures
\begin{align}
\label{muLN}
\mu_{L,N} =\nu \big[\,\cdot\,\big| S_L =N\big]\ .
\end{align}
To simplify notation we will write $\mu=\mu_{L,N}$ from Section \ref{sec:proof_main} onwards.\\

We study the large scale behaviour of the process in the thermodynamic limit with particle density $\rho\geq 0$, i.e.
\be\label{thermo}
L,N=N_L \to\infty\ ,\quad\mbox{such that }N/L \to\rho\ .
\ee
For $b>2$ and $\rho >\rho_c$ (\ref{zcrc}), the process is known to exhibit a condensation transition in the following sense. Denoting by
\begin{align}\label{max}
M_L =\max_{x\in\Lambda} \eta_x
\end{align}
the maximum occupation number as a relevant order parameter, we have a weak law of large numbers \cite{evans00,jeonetal00,grosskinskyetal03}
\be
M_L /L\to\left\{\bacl 0 &,\ \rho\leq \rho_c\\ \rho -\rho_c &,\ \rho >\rho_c\ea\right.\quad\mbox{as }L,N\to\infty\ (\ref{thermo})\ ,
\ee
where convergence holds in distribution w.r.t.\ the sequence $\mu_{L,N}$. 
A corresponding CLT has been established in \cite{armendarizetal08}, and the fluctuations are Gaussian and of order $\sqrt{L}$ for $b>3$, and of order $L^{1/(b-1)}$ with an associated stable law for $2<b\leq 3$.  Moreover, the distribution of the configuration outside the maximum is known to converge to the maximal product measure $\nu$ with density $\rho_c$ in the limit (\ref{thermo}), so the largest occupation number in the bulk outside the maximum is typically of order $L^{1/(b-1)}\ll L$, which holds for all $b>2$ independently of the fluctuations of the maximum. Therefore for supercritical densities $\rho >\rho_c$ configurations exhibit a unique extensive maximum with high probability, called the condensate, and the rescaled location of the maximum
\be
\psi_L (\eta ) =\frac1L\big\{ x\in\Lambda :\eta_x =M_L (\eta )\big\}
\ee
is given by a single site in the rescaled lattice $\frac1L\Lambda$. By translation invariance, $\psi_L$ is distributed uniformly in $\frac1L\Lambda$ under $\mu_{L,N}$, and ergodicity of the process on $X_{L,N}$ implies that $\psi_L \big(\eta (t) \big)$ visits the whole lattice on a long time scale. It is expected that this dynamics is metastable, i.e. $\psi_L \big(\eta (t) \big)$ is constant for a long (random) time interval, and then changes abruptly to a new value which depends only on the current location of the condensate. \\

Our main result is that, indeed, for large enough $b>2$ and $\rho >\rho_c$ the zero-range process $\big(\eta (t) :t\geq 0\big)$ in the thermodynamic limit (\ref{thermo}) exhibits metastability with respect to the observable $\psi_L$ on the time scale
\be\label{thetal}
\theta_L :=L^{1+b}\ .
\ee
This means that the sequence of processes $\big(\frac1L\,\psi_L (\eta_{t\theta_L} ):t\geq 0\big)$ converges to a Markovian limit process on the unit torus $\T =\R /\Z$, the scaling limit of $\frac1L\Lambda$. A rigorous version of this result is provided in the next subsection in Theorem \ref{main2}, including an exact formulation of the mode of convergence and required assumptions.

\subsection{The trace process and main result \label{trace_process}}

The rigorous formulation and proof of our main result on the large scale metastable dynamics of the condensate follows the martingale approach developed in \cite{beltranetal09,beltranetal09b,beltranetal14}, which requires a partition of the state space as a first step. We define the well $\Ecal^x$ as the set of configurations where the condensate is located at $x\in\Lambda$,
\be\label{valley}
\Ecal^x :=\big\{\eta\in X_{L,N} :\eta_x =M_L \geq N-\rho_c L-\alpha_L ,\eta_y \leq \beta_L \mbox{ for all }y\neq x \big\}\ .
\ee
We choose the scales in this definition as
\be\label{alphabeta}
\alpha_L=L^{\frac{1}{2}+\frac{5}{2b}} \quad\mbox{and}\quad\beta_L=2\lfloor L^{\frac{4}{b-1}}\rfloor\ ,
\ee
which allow for typical fluctuations of the condensate size of order $\sqrt{L}$, and of bulk occupation numbers of order $L^{1/(b-1)}$ for any $b>3$. Note that under our conditions on $b$ in Theorem \ref{main2} we also have $\alpha_L ,\beta_L \ll L$. We choose $\beta_L$ as a sequence of integers for later convenience, while the exponent in $\alpha_L$ is optimal for the estimates in Section \ref{Prop_3.4}. 
We denote by
\be
\Ecal :=\cup_{x\in\Lambda} \Ecal^x\quad\mbox{and}\quad \Delta:= X_{L,N}\setminus\Ecal
\ee
the set of all wells and its complement. The first result states that the process started from a well spends only a negligible amount of time on the set $\Delta$ on the timescale $\theta_L$. As usual, we denote by $\P_\eta$ and $\E_\eta$ the path measure and expectation of the process $\eta (\cdot )$ with generator (\ref{zrgen}) and with initial condition $\eta (0)=\eta$.

\begin{proposition}\label{trace} \textbf{Replacement by the trace process}\\
For the process defined in (\ref{zrgen}) with rates (\ref{rates}) and $b>20$, we have for all $t>0$
\be\label{trasub}
\sup_{\eta \in \Ecal} \E_{\eta} \Big[ \int_0^{\theta_L t} \1_\Delta (\eta (s))\, ds\Big]\Big/\theta_L \to 0\quad\mbox{as }L,N\to\infty\ (\ref{thermo})\mbox{ with }\rho >\rho_c\ .
\ee
\end{proposition}

 The proof of this result is given in Section \ref{nullDelta}. The larger lower bound on the parameter $b$ is a purely technical restriction which we discuss in Section \ref{sec:discussion}. 
Since the process spends only negligible time outside $\Ecal$, re-parametrizing time by the local time on $\Ecal$ will not change the process in the limit. Denote by 
\be\label{ltime}
\Tcal_t :=\int_0^t \1_\Ecal (\eta_s )\, ds\quad\mbox{and}\quad \Scal_t :=\sup\big\{ s\geq 0:\Tcal_s \leq t\big\}
\ee
the local time on $\Ecal$ and its generalized inverse, respectively. The trace process is then defined as
\be\label{tracepr}
\big(\eta^\Ecal (t) :t\geq 0\big)\quad\mbox{with}\quad \eta^\Ecal (t) :=\eta_{\Scal_t}\ ,
\ee
which takes values in the set $\Ecal$ and is well defined since $\Tcal_t$ diverges $\P_\eta -a.s.$ as $t\to\infty$ due to the irreducibility on the finite state space $X_{L,N}$. As is shown in \cite{beltranetal09b}, Section 6.1, (\ref{tracepr}) is in fact an irreducible Markov process on $\Ecal$ with jump rates
\begin{align}
\label{trace_rates}
r^\Ecal (\eta ,\xi ):=r(\eta ,\xi )+\sum_{\zeta\in\Delta} r(\eta ,\zeta )\,\P_\zeta [T_\Ecal =T_\xi ]\,.
\end{align}
Here $r(\eta ,\xi ):=\sum_{x\in\Lambda} g(\eta_x )\, \frac{1}{2}\big(\1_{\xi} (\eta^{x,x+1} )+\1_{\xi} (\eta^{x,x-1} )\big)$ are the jump rates of the process $\eta (\cdot)$ with generator (\ref{zrgen}), and
\be
T_{\cal F} :=\inf\big\{ t\geq 0: \eta_t \in {\cal F}\big\}
\ee
denotes the hitting time of a set ${\cal F}\subset X_{L,N}$. For point subsets ${\cal F}=\{\xi\}$ we use the shorthand $T_{\{\xi\}} =T_\xi$. The trace process has generator
\begin{align}
\label{trace_generator}
{\cal L}^\Ecal f(\eta)=\sum_{\xi \in \Ecal} r^\Ecal(\eta, \xi)\,\big[f(\xi)-f(\eta)\big]\,,
\end{align}
and reversible invariant measure $\mu_{L,N} \big[\,\cdot\,\big|\Ecal \big]$ restricted to $\Ecal$. To simplify notation, we will call
\begin{equation}
\label{mu_ecal}
\mu^{\Ecal}[\cdot]=\mu_{L,N} \big[\,\cdot\,\big|\Ecal \big]. 
\end{equation}
On the set of wells $\Ecal$ the rescaled location of the maximum can be written as
\be\label{maxloc}
\psi_L (\eta ):=\frac{1}{L}\sum_{x\in\Lambda} x\,\1_{\Ecal^x} (\eta )\in\T =\R /\Z\ ,
\ee
and is well defined without degeneracy. We can now state our main result.

\begin{theorem}\label{main2}
Consider the zero-range process with generator (\ref{zrgen}) and rates (\ref{rates}) with $b>20$. 
Fix a sequence of initial conditions $\eta (0)\in\Ecal$ such that $\psi_L (\eta (0))\to Y_0 \in\T$. In the thermodynamic limit (\ref{thermo}) with $N/L\to\rho >\rho_c$ we have for the trace process (\ref{tracepr}) on the time scale $\theta_L =L^{1+b}$
\be
\Big( \psi_L \big(\eta^\Ecal (t\theta_L )\big) :t\geq 0\Big)\Rightarrow \big( Y_t :t\geq 0\big)\ ,
\ee
where convergence holds weakly on the path space $D\big( [0,\infty ),\T\big)$ in the usual Skorohod topology.\\
The process $(Y_t :t\geq 0)$ has stationary and independent increments on $\T$, initial condition $Y_0$ and generator
\be\label{limgen}
\Lcal^\T f(u)=\int_{\T\setminus\{ 0\}} r^\T (v) \big( f(u+v)-f(u)\big)\, dv
\ee
for all Lipschitz continuous functions $f:\T\to\R$, with rates
\be\label{rlimi}
r^\T (v)=\frac{1}{z_c\, I_b \, (\rho -\rho_c )^{b+1}} \frac{1}{d_\T (0,v)}\ ,\quad v\in\T\ .
\ee
Embedding the unit torus in $[0,1)\subset\R$, the distance on $\T$ is given by $d_\T (0,v)={|v|(1-|v|)}$. $z_c$ is the normalization of the invariant measure (\ref{zcrc}), and
\be
I_b :=\int_0^1 x^b (1-x)^b \, dx =\frac{\Gamma (b+1)^2}{\Gamma (2b+2)}
\ee
is a constant depending only on $b$.
\end{theorem}

Note that (\ref{trasub}) and our main result do not apply for initial conditions $\eta (0)\in\Delta$, which are, however, untypical if the process starts from the stationary measure $\mu_{L,N}$ (cf.\ Corollary \ref{shaft}). The mode of convergence in terms of the trace process as presented here has been introduced in \cite{beltranetal09,beltranetal09b}, and extended recently in \cite{beltranetal14} to a more general context. The (random) time change in the definition of the trace process, which is negligible in the limit $L\to\infty$, can also be absorbed in a definition of a suitable topology on the path space of the limiting process. Further discussion, including possible extensions of our result, is provided in Section \ref{sec:discussion}.\\

\section{Proof of the main result \label{sec:proof_main}}

The proof of Theorem \ref{main2} uses a standard approach to establish existence of a limit process by a tightness argument, and to identify the limit by the solution of a martingale problem. We follow the method outlined in \cite{beltranetal14}, proofs of auxiliary results used below are given in the rest of the paper. We also discuss the main novelties and possible extensions.\\

Here and in the following sections we adopt a few shorthands and conventions to avoid an overload of notation. We write $\mu =\mu_{L,N}$ for the invariant measure of the full process, and $\mu^\Ecal$ for the invariant measure of the trace process \eqref{mu_ecal}. Constants denoted by $C$ are independent of $L$ and $N$, and can vary from line to line.

\subsection{Proof of Theorem \ref{main2} \label{proof_main}}

The proof of convergence holds on arbitrarily large compact time intervals for the limit process, and throughout this section we denote the length of this interval by $T>0$. Let
\begin{equation}
\label{YLproc}
(Y_t^L :t\geq 0)\in D\big( [0,T],\T\big)\quad\mbox{with}\quad Y_t^L :=\psi_L \big(\eta^\Ecal (\theta_L t)\big)\ ,
\end{equation}
be the speeded up process of the rescaled maximum location, and let $\Q^L$ be its distribution on path space $D\big( [0,T],\T\big)$. 

\begin{proposition}\label{lprocess} \textbf{Tightness}\\
Under the conditions of Theorem \ref{main2} with $b>20$ the sequence $\Q^L$ of path space distributions is tight on $D\big( [0,T],\T\big)$.
\end{proposition}

The proof is given in Section \ref{tight} where a control of the quadratic variation excludes the accumulation of jumps and ensures that limit points have right-continuous paths. Tightness implies existence of sub-sequential weak limits of $\Q^L$ in the Skorohod topology, and we denote any such weak limit by $\Q$. In order to identify the limit we need to show that for all $t\leq T$ and all Lipschitz-continuous functions $f\in\mathrm{Lip}(\T )$
\be\label{lmarteq}
f(\omega_t )-f(\omega_0)-\int_0^t \Lcal^\T f(\omega_s )ds\quad\mbox{is a martingale}\ ,
\ee
where $\omega_t :D\big( [0,\infty ),\T\big)\to \T$ is the coordinate process on path space. Together with the uniqueness result for the martingale problem associated with ${\Lcal^\T}$ proved in Subsection \ref{mapro}, this implies convergence of $\Q^L$ and characterizes the limit $\Q$ as the law of the Markov process $\big( Y_t :t\in [0,T]\big)$ with generator $\Lcal^\T$ (\ref{limgen}), because Lipschitz functions form a core for this generator. 
Since $T>0$ is arbitrary, this implies Theorem \ref{main2}. Precisely, we need to show that
\be\label{nts}
\E^\Q \bigg[g\big( (\omega_u :0\leq u\leq s) \big)\Big( f(\omega_t )-f(\omega_s )-\int_s^t \Lcal^\T f(\omega_u )\, du\Big)\bigg] =0\ ,
\ee
for all $0\leq s<t\leq T$ and all bounded, continuous functions $g:D\big( [0,T],\T\big)\to\R$. Since $\big(\eta^\Ecal (\theta_L t):t\in [0,T]\big)$ is a Markov process with generator $\theta_L \Lcal^\Ecal$, we know that
\begin{align*}
\label{marprop}
& f(Y^L_t )-f(Y^L_0 ) -\theta_L \int_0^t \Lcal^\Ecal ( f\circ\psi_L ) (\eta^\Ecal (\theta_L s))\, ds=\nonumber\\
&=f(Y^L_t )-f(Y^L_0 ) -\int_0^t \Lcal^\T f(Y^L_s)\, ds +\int_0^t \Big( \Lcal^\T f(Y^L_s)-\theta_L \Lcal^\Ecal ( f\circ\psi_L ) (\eta^\Ecal (\theta_L s))\Big)\, ds
\end{align*}
is a martingale for all $t\in [0,T]$ and $L\in\N$. We will establish below that
\be\label{suppe}
\sup_{\eta\in\Ecal} \E_\eta \bigg| \int_0^t \Big( \Lcal^\T f(Y^L_s)-\theta_L\Lcal^\Ecal ( f\circ\psi_L ) (\eta^\Ecal (\theta_L s))\Big)\, ds\bigg|\to 0
\ee
as $L\to\infty$, which implies that
\be
\E^{\Q^L} \bigg[ g\big( (\omega_u :0\leq u\leq s) \big)\Big( f(\omega_t )-f(\omega_s )-\int_s^t \Lcal^\T f(\omega_u )\, du\Big)\bigg] \to 0\ ,
\ee
for all $0\leq s<t\leq T$ and bounded, continuous $g:D\big( [0,T],\T\big)\to\R$. To identify the limit of the left-hand side with (\ref{nts}) we use the following result, which is immediate from Lemma \ref{chara} in Section \ref{sec:martingale}.
\begin{proposition}\label{martcon}
Let $M_t^f(\omega)=f(\omega_t)-f(\omega_s)-\int_s^t {\cal L}^{\T} f(\omega_u)\ du$. Then
\be
\E^{\Q^L} \Big[ g\big( (\omega_u :0\leq u\leq s) \big)\, M_t^f(\omega )\Big]\to \E^{\Q} \Big[ g\big( (\omega_u :0\leq u\leq s) \big)\, M_t^f(\omega )\Big]
\ee
as $L\to\infty$ for any $t\in [0,T]$.
\end{proposition}

To conclude the proof we have to show (\ref{suppe}), i.e. that we can replace the generator of the trace process with that of the limit process (cf. Section 3 in \cite{beltranetal14}). This is the main part of the paper and is divided into several steps. Since we cannot compare the generators $\Lcal^\Ecal$ and $\Lcal^\T$ directly, we will introduce an auxiliary processes on $\T$ with generator $\Lcal^\Lambda$ which is explained in detail below, and we rewrite the time integral in (\ref{suppe}) as

\bea
\lefteqn{\int_0^t \Big( \Lcal^\T f(Y^L_s)-\theta_L \Lcal^\Ecal ( f\circ\psi_L ) (\eta^\Ecal (\theta_L s))\Big)\, ds=}\nonumber\\
& &\qquad =\int_0^t \Big(\Lcal^\T f(Y^L_s ) -\theta_L \Lcal^{\Lambda} f (Y^L_s )\Big)\, ds\label{interv}\\
& &\qquad\quad +\theta_L \int_0^t \Big( \Lcal^\Lambda f(Y^L_s ) -\Lcal^\Ecal (f\circ\psi_L )(\eta^\Ecal (\theta_L s))\Big)\, ds \label{equili}\ .
\eea
The goal is now to show that the terms in (\ref{interv}) and (\ref{equili}) vanish individually as $L\to\infty$ in $L^1$-norm. 
In (\ref{equili}) we compare the trace process to the auxiliary process on the rescaled lattice $\frac1L\Lambda\subset\T$, which describes the effective jumps of the condensate location. Its generator is defined as
\be\label{auxgen}
\Lcal^\Lambda f(x/L):=\sum_{z\in\Lambda ,z\neq 0} r^\Lambda (z)\big( f\big( (x+z)/L\big) -f(x/L)\big)\ ,
\ee
with jump rates (using translation invariance)
\be\label{auxrate}
r^\Lambda (z)=\sum_{\xi\in\Ecal^z }\mu |_{\Ecal^0} \big( r^\Ecal (.,\xi )\big) =\frac{1}{\mu [\Ecal^0 ]} \sum_{\eta\in\Ecal^0\atop \xi\in\Ecal^z} \mu [\eta]\, r^\Ecal (\eta ,\xi )\ ,
\ee
given by the expected rate between wells for the trace process. Therefore, we can also write the generator (\ref{auxgen}) as the expectation $\Lcal^\Lambda f(x/L) =\mu |_{\Ecal^x} \big(\Lcal^\Ecal (f\circ\psi_L )\big)$\ . Before the location of the condensate changes, the process remains in the same well $\Ecal^x$ for a long enough time to equilibrate, and the transition between wells is effectively described by stationary averages of jump rates as in (\ref{auxrate}), which is established in the next result.

\begin{proposition}\label{equilibration} \textbf{Equilibration in the wells}\\
Under the conditions of Theorem \ref{main2} with $b>20$,
\be\label{equi}
\sup_{\eta\in\Ecal} \E_\eta \bigg|\theta_L \int_0^{t} \Big(\Lcal^\Lambda f(Y^L_s )-\Lcal^\Ecal \big( f\circ\psi_L \big)\big(\eta^\Ecal (\theta_L s) \big) \Big) ds\bigg|\to 0
\ee
as $L\to\infty$, for all $t\in [0,T]$.
\end{proposition}

In addition to using Lemma \ref{unibound} given below, the proof requires an estimate on the relaxation and mixing times within a well, which have to be strictly smaller than $\theta_L$, and is provided in Section \ref{sec:equi}. \\

In (\ref{interv}) we replace the generator of the auxiliary process with that of the limit process using the following result.

\begin{proposition}\label{inter} \textbf{Dynamics between wells}\\
Under the conditions of Theorem \ref{main2} with $b>5$,
\be\label{cap}
\sup_{\eta\in\Ecal} \E_\eta \bigg|\int_0^{t}\Big( \Lcal^\T f\big(Y_s^L \big) -\theta_L \Lcal^\Lambda f\big(Y_s^L\big)\Big) ds \bigg|\to 0
\ee
as $L\to\infty$, for all $t\in [0,T]$.
\end{proposition}

The proof is given in Section \ref{sec:regul} and requires sharp bounds on the transition rates of the auxiliary process, which are provided by capacity estimates in Section \ref{sec:cap}. In order to get matching upper and lower bounds, an important new step in this part of the proof is to regularize the rates $r^\Lambda (z)$ on an intermediate scale and use the regularity of the test function $f$, which is explained in detail in Section \ref{sec:regul}. This finishes the proof of our main result, Theorem \ref{main2}.\CQFD\\

An important estimate that is used in the proof of tightness, equilibration and replacement of the trace process is the following uniform bound on the exit rate from a well. The proof of this Lemma is given by a coupling argument in Section \ref{sec:coup}, which is one of the crucial new results of this paper. 

\begin{lemma}\label{unibound} \textbf{Uniform bound on the exit rate}\\
Under the conditions of Theorem \ref{main2} with $b>20$, there exists a constant $C>0$ such that the exit rate from a well is uniformly bounded by
\begin{equation}
\label{unif_rates}
\sup_{\eta\in\Ecal^0} \sum_{\xi\not\in\Ecal^0} r^\Ecal (\eta ,\xi )\leq 
 \frac{C}{L^5\log^2 L}\ .
\end{equation}
\end{lemma}

\subsection{Discussion\label{sec:discussion}}

\paragraph{Main new ideas.}
Our proof follows the martingale approach outlined in \cite{beltranetal14}, which was previously applied to zero-range processes on a fixed lattice in the limit $N\to\infty$ \cite{beltranetal09}. In contrast to this case, the thermodynamic limit (\ref{thermo}) considered here involves a significant change in the complexity of metastable wells, since the sizes of the wells and the number of transition paths between them increase with $L$. Following the discussion in \cite{fernandezetal14b}, this presents a technically more challenging metastability scenario, in particular since the free energy barriers of the metastable wells in the zero-range process are only logarithmic. 
We quickly summarize the resulting conceptual and technical difficulties and the three main novel contributions of this paper to overcome them.

\begin{itemize}
	\item In previous work with limits of diverging particle density \cite{beltranetal09,rebecca} the depth of the metastable wells was dominating their size, and they could effectively be replaced by individual configurations, so-called attractor states. The repeated visits to those configurations before hitting another well lead to a relatively simple renewal-type proof for equilibration in the wells. In our case of finite particle density the size of metastable wells increases exponentially with $L$ and these types of arguments do not apply. Instead, we have to describe the metastable states as probability distributions on wells. As outlined in \cite{beltranetal14,fernandezetal14b} a suitable candidate is simply the restriction of the stationary measure, and to prove Proposition \ref{equilibration} we use suitable dynamics restricted to a well following \cite{beltranetal14}, see Section \ref{sec:equi}. 
Using a standard jacknife estimate \cite{es1981} we establish a Poincar\'e inequality comparing these dynamics to independent birth-death chains in Section \ref{trelsub}, and obtain a bound on the relaxation time on the metastable well of order $L^4$ in Lemma \ref{trelest}. This is clearly not optimal, and a sharp bound is expected to be of order $L^2$, which would moderately improve our conditions on the parameter $b$ in Proposition \ref{equilibration} and Lemma \ref{unibound}, currently $b>20$, to $b>13$.
	
	\item Resulting from capacity bounds presented in Section \ref{sec:cap}, the transition rates between wells $r^\Lambda (z)$ (\ref{auxrate}) scale like $1/(|z| (L-|z|))$ (\ref{laca}). The total exit rate from a well on the scale $\theta_L$ diverges as $\log L$ (\ref{rabo}) in the thermodynamic limit, and the exit time distribution from a well degenerates in the scaling limit and does not converge to an exponential random variable. This is due to small jumps which accumulate at diverging rate, but are still negligible on the macroscopic scale for the rescaled limit dynamics, which is established rigorously by showing tightness in Section \ref{tight}. Related to this, the errors in (\ref{laca}) also do not allow to get matching upper and lower bounds on transition rates between individual wells, in analogy with results in \cite{rebecca} for diverging density. 
We address both issues in the standard thermodynamic limit for fixed density, using the symmetry of the system and the regularity of the limit dynamics on the unit torus, which are fully characterized by Lipschitz-continuous test functions. We regularize the dynamics on an intermediate scale in Section \ref{sec:regul} to get matching bounds on rates between regularized sets of wells, which are sufficient to derive the limiting generator and prove Proposition \ref{inter}. As an additional technical difficulty we have to keep track of corrections of capacity bounds to leading order in $L$, which require very precise estimates on the stationary measure summarized in Section \ref{invmeas}. 
Sections \ref{sec:regul} and \ref{sec:cap} are independent of the rest of the proof, and the only restrictions on the parameter $b$ arising there are given by equations \eqref{epsul} and \eqref{ellbarell}, resulting in a much weaker condition of $b>5$.
	
	\item While the limiting metastable dynamics are determined by stationary averages of transition rates $r^\Lambda$ given in (\ref{auxrate}), a uniform control of the exit rates from a well is important to estimate error terms as is discussed in general in \cite{beltranetal14}. Some particular examples of spin systems where this has been achieved are mentioned in \cite{bovierbook}, but note that using capacities between individual configurations in the thermodynamic limit would lead to bounds that diverge exponentially in the system size. Here we derive a uniform upper bound on exit rates scaling as $L^{-5}\log^{-2}L $ in Lemma \ref{unibound}, which is proved in Section \ref{sec:coup} by a novel coupling construction with a growing number of birth-death chains. The number of chains increases only linearly in time, and the construction ensures that in the event of changing well, at least one of the chains has grown a condensate. This is a central auxiliary result of the paper, and is used in the proof of equilibration in the wells (Prop.~\ref{equilibration}), replacement by the trace process (Prop.~\ref{trace}), as well as for tightness (Prop.~\ref{lprocess}). 
The proof of Lemma \ref{unibound} requires $b>20$, and sharper estimates on equilibration times in wells would slightly improve this condition as mentioned above.

\end{itemize}
	
The bottleneck of the method leading to restrictive conditions on the parameter $b$ results from the inherent competition between depth and complexity of metastable wells. Both quantities increase with the size of the wells, and the aim to maximize depth and minimize complexity leads to an optimal choice of their size that enters most prominently in the uniform bounds on exit rates. 
The coupling argument in Section \ref{coupling} gives a lower bound on the expected time to change wells for the full zero-range process. To turn this into an estimate for the trace process, we have to bound the time spent outside the wells on the set $\Delta$, which is controlled by the invariant measure in Corollary \ref{shaft}. Larger wells of increased depth improve this bound, but at the same time lead to an increase in the number of transition paths to other wells. Both effects compete and affect the uniform bound on the exit rates in Lemma \ref{unibound}. It turns out that the optimal size is controlled by the parameter $\beta_L$ (\ref{valley}), and the crucial estimate in this context is (\ref{bootstrapacal}) for the probability that the trace process has changed well on the time scale $\theta_L$, where terms of the form $\beta_L^{b-1}$ and $\beta_L^{1-b}$ appear in bounds of the right-hand side. The best choice of $\beta_L$ in (\ref{alphabeta}) leads to a bound in Lemma \ref{unibound} which is small enough for the required estimates in the proofs of Propositions \ref{trace}, \ref{lprocess} and \ref{equilibration} as long as $b>20$. The mechanism leading to this constraint is of a fundamental nature, and it seems very hard to significantly improve this with the techniques used in this paper.\\

The optimal choice of the second parameter $\alpha_L$ in (\ref{valley}), (\ref{alphabeta}) controlling the fluctuations of the condensate is of a more technical nature, and arises from two conditions on the upper bound error of the capacity estimates in (\ref{epsu}). To obtain the upper bound, a test function is constructed that concentrates on wells and tubes (\ref{tubedef}), which are subsets of $\Delta$ where the transition paths between wells are expected to concentrate on. The size of tubes and wells both increase with the parameter $\alpha$, and increasing $\alpha$ therefore reduces the distances between wells and increases the corresponding capacity. At the same time, increasing tube size improves the approximation and therefore decreases the upper bound on those capacities. The optimal choice of $\alpha_L$ results from estimates in Lemmas \ref{neglig1} and \ref{neglig3} and only leads to the weak requirement $b>5$.\\

\paragraph{Possible extensions.}
We focus on symmetric, nearest-neighbour probabilities $p(x,y)=\frac{1}{2}\delta_{y,x-1}+\frac{1}{2}\delta_{y, x+1}$ in one dimension with periodic boundary conditions, two properties of which are essential for our proof.

\begin{itemize}
	\item Symmetry; this leads to reversible dynamics, which is a necessary condition for the potential theoretic estimates on transition rates we use in Section \ref{sec:cap}. There is significant recent research interest on extended Dirichlet principles for non-reversible systems which involve double variational formulas \cite{beltranetal12,benoisetal13,gaudilliereetal14}, see also references in \cite{bovierbook}, Chapter 7. This has been applied to the totally asymmetric zero-range process on a fixed lattice in \cite{landim13}, but since a result in the thermodynamic limit requires much better control on error bounds an extenstion to this case would be a significant technical complication. While a (non-optimal) relaxation time estimate for non-reversible systems can probably be obtained, we also make critical use of symmetry of the jump rates in the regularization step in Section \ref{sec:regul}, which is not obviously adapted to asymmetric situations.
	\item Translation invariance; the results in \cite{beltranetal09,beltranetal09b} apply to zero-range processes without this property, leading to spatially inhomogeneous limit dynamics on a fixed lattice which are directly related to the choice of $p(x,y)$. We use translation invariance in our proof for equilibration and also in the regularization step in Sections \ref{sec:equi} and \ref{sec:regul}. Specific simple examples of non-translation invariance, such as alternating $p(x,y)$ or isolated inhomogeneities, can be treated as a direct extension of our result on a case-by-case basis. However, it is not clear how to formulate a result in the generality covered in \cite{beltranetal09,beltranetal09b}, where the first problem already may arise when defining the limiting process if  the probabilities $p(x,y)$ do not have good scaling properties.
\end{itemize}

Whenever the $p(x,y)$ are symmetric and translation invariant and admit a well defined scaling limit of the capacities of a single random walk on $\Lambda$ analogously to (\ref{laca}), our results can be extended without much effort to lead to limit processes with stationary, independent increments. General finite range symmetric $p(x,y)$ in one or higher dimensions which scale to Brownian motion should all give the same result related to the corresponding harmonic functions of a single walker, appropriately modified on the torus. Note that in three and more dimensions these functions have a constant scaling limit leading to uniform displacement of the condensate, with expected logarithmic corrections in two dimensions. Also if $p(x,y)$ has range diverging with $L$ with well defined scaling limits for capacities our result can be directly adapted, including for example uniform $p(x,y)$ which leads to uniform condensate dynamics on the limiting torus. Due to the special properties of one-dimensional diffusion the case covered here is already one of the most interesting. Note also that even for finite-range $p(x,y)$ the condensate dynamics will be non-local in all dimensions, in contrast to an analogous result for inclusion processes \cite{grosskinskyetal13}.\\

In addition to the jump rates $g(n)$ (\ref{rates}) considered here, there are various other choices that lead to condensing zero-range processes (see e.g.\ \cite{raffertyetal15} and references therein). A well studied example is for rates with asymptotic behaviour $g(n)\sim 1+b/n^\lambda$ with $\lambda\in (0,1)$ leading to stretched exponential tails for the stationary measure \cite{chlebounetal10,agl}. The lighter tail of the measure increases the depth of the metastable wells and leads to free energy barriers that grow sublinearly in the system size $L$, which is much faster than the logarithmic growth for the present model corresponding to $\lambda =1$. Therefore, even though an actual generalization would require considerable work, we expect that all our techniques can be applied and some estimates, in particular the ones in Section \ref{sec:cap}, should even get easier.\\

The depth of metastable wells increases with particle density, and with the parameter $b$ which determines the tail of the stationary measure. In contrast to the restrictive conditions on $b$, our proof is robust in this system parameter and does not require any additional constraints on the particle density except $\rho >\rho_c$. As long as the excess mass $N-\rho_c L\gg\sqrt{L\log L}$ it is known that the zero-range process still exhibits the condensation phenomenon on the level of stationary measures \cite{agl}. It would therefore be interesting to investigate, but is beyond the scope of this paper, in how far our results can be applied also for subextensive excess mass. 

\paragraph{Related work.}
The results presented here and other metastability results for zero-range processes only concern the  stationary dynamics of a single condensate, corresponding to the slowest time scale in the system. In particular, these results only apply when the process is started on one of the metastable wells. Starting the process from a uniform initial condition with supercritical or diverging density leads to a dynamic formation of the single condensate on a different, faster time scale. This has been discussed heuristically in \cite{grosskinskyetal03,godreche}, where it was found that after a rapid formation of several large clusters, they exchange particles in a coarsening process on the time scale $L^3 \ll L^{1+b}$ which leads to the formation of a single condensate. Note that the location of the clusters on this time scale is fixed. Recently the first rigorous result in this context has been obtained on a fixed lattice in the limit $N\to\infty$ in \cite{jaraetal}. The coarsening dynamics takes place outside the set $\Ecal$ and is therefore approached with entirely different methods than the ones in this paper. For the condensation phenomenon in the inclusion process the dynamics of the condensate and the coarsening process take place on the same time scale, and both have been rigorously understood in \cite{grosskinskyetal13}.\\

\section{Equilibration dynamics in the wells\label{sec:equi}}

To establish an upper bound for the thermalization time scale of the trace process in a metastable well we follow the procedure outlined in \cite{beltranetal14}. We introduce a process restricted to the well (called reflected process in \cite{beltranetal14}), which is reversible w.r.t.\ the restricted measure $\mu$, and estimate the relaxation time of this process. A general result from \cite{beltranetal14} can then be applied to yield a simple proof of Proposition \ref{equilibration}.

\subsection{Proof of Proposition \ref{equilibration}\label{sec:equi2}}

For each well $\Ecal^x$, $x\in\Lambda$, the restricted process is defined by the generator
\be\label{engen}
\Lcal^x f(\eta )=\sum_{\zeta\in\Ecal^x } r(\eta ,\zeta )\big( f(\zeta )-f(\eta )\big)\quad\mbox{for all }\eta\in\Ecal^x\ .
\ee
As before, $r(\eta ,\zeta )$ denote the jump rates of the full zero-range process, and jumps outside the well $\Ecal^x$ are suppressed. Note that this is not equal to the trace process on $\Ecal^x$, which has additional rates at the boundary of $\Ecal^x$. It is easy to see that this process is irreducible on $\Ecal^x$ and that it is reversible w.r.t.\ the restricted measure
\begin{equation}
\label{rest_meas}
\mu^x :=\mu [.|\Ecal^x ]\ .
\end{equation}
The following estimate on the relaxation and mixing times of the restricted process will be used to prove Proposition \ref{equilibration}.

\begin{lemma}\label{trelest}
The relaxation time $\tr$ and the $\epsilon$-mixing time $\tm (\epsilon )$ of the restricted process $\Lcal^x$ (\ref{engen}) on $\Ecal^x$ are independent of $x$ and bounded by
\begin{equation}\label{trele}
\tr \leq C L^4 \quad\mbox{and}\quad \tm (\epsilon )\leq C L^5(1+L^{-1} \log\tfrac{1}{\epsilon})\ ,
\end{equation}
for some constant $C>0$, depending only on the fixed parameters $\rho$ and $b$.
\end{lemma}

The proof of this Lemma uses path counting techniques to establish a Poincar\'e inequality for the restricted process, which is independent of the rest of this section and is therefore postponed to Subsection \ref{trelsub}. 

We will use the following $L^2$ estimate for the ergodic average of a function that has mean zero on all wells, starting from the stationary measure $\mu^\Ecal =\mu [.|\Ecal ]$ restricted to the wells. The proof of can be found in \cite{beltranetal14}, Section 3.1.

\begin{lemma}\label{ergbound}
For every function $f:\Ecal \to\R$ which has vanishing mean $\mu^x (f)=0$ for all $x\in\Lambda$, and for all $t>0$
\begin{equation}
\label{times}
\E_{\mu^\Ecal} \bigg|\int_0^t f(\eta^\Ecal (s))\, ds\bigg|^2 \leq 24 t\sum_{x\in\Lambda} \mu [\Ecal^x ]\, \tr^x \,\mu^x (f^2)\ ,
\end{equation}
where $\tr^x =\tr$ is the relaxation time on $\Ecal^x$ of the restricted process, which does not depend on $x$ in our case.
\end{lemma}

By translation invariance we can simply focus on initial conditions in a chosen well $\Ecal^0$ in the following. To prove Proposition \ref{equilibration} we have to show (\ref{equi}), i.e. prove that
\bea
\lefteqn{\E_\eta \bigg|\int_0^{\theta_L t} \Big(\Lcal^\Lambda f\big(\psi_L (\eta^\Ecal (s))\big) -\Lcal^\Ecal \big( f\circ\psi_L \big)\big(\eta^\Ecal (s) \big) \Big) ds\bigg| }\nonumber\\
& &=\E_\eta \bigg|\int_0^{\theta_L t} \sum_{z\in\Lambda\atop z\neq 0} \big( r^\Lambda (z)-r^\Ecal (\eta (s),\Ecal^{z} )\big) \Big( f\big(\psi_L (\eta (s))+z/L\big) -f\big(\psi_L (\eta (s))\big)\Big)\bigg| \to 0
\eea
as $L\to\infty$ for all Lipschitz functions $f:\T\to\R$ and $\eta\in\Ecal^0$. For the total jump rate of the trace process into another well $z\neq 0$, we have used the obvious notation
\[
r^\Ecal (\eta ,\Ecal^{z} ) =\sum_{\zeta \in\Ecal^z} r^\Ecal (\eta ,\zeta )\quad\mbox{for all }\eta\in\Ecal^0 \ .
\]
With the definition of $r^\Lambda (z)$ in (\ref{auxrate}) we have
\[
r^\Lambda (z)=\mu^x \Big( r^\Ecal (.,\Ecal^{x+z} )\Big)\quad\mbox{for all }x,z\in\Lambda\ ,
\]
and since the function $f\circ\psi_L$ is constant on all wells, the function
\be
h_f (\eta ):=\sum_{z\neq 0} \big( r^\Lambda (z)-r^\Ecal (\eta ,\Ecal^{x+z} )\big) \Big( f\big(\psi_L (\eta )+z/L\big) -f\big(\psi_L (\eta )\big)\Big)
\ee
has mean zero under $\mu^x$ for all $x\in\Lambda$, independently of $f$. On every well $\Ecal^x$, we can estimate its second moment as
\bea
\mu^x \big( h_f^2 \big) &\leq &\mu^x \bigg( \bigg( \sum_{z\neq 0} r^\Ecal (\cdot ,\Ecal^{x+z} ) \Big( f\big(\psi_L (\cdot )+\frac{z}{L}\big) -f\big(\psi_L (\cdot )\big)\Big)\bigg)^2 \bigg)
\leq C_f^2 \mu^x \bigg( \Big( \sum_{z\neq 0} r^\Ecal (\cdot ,\Ecal^{x+z} ) \Big)^2 \bigg) \nonumber\\
&\leq &C_f \Big\| \sum_{z\neq 0} r^\Ecal (\cdot ,\Ecal^z )\Big\|_\infty \sum_{z\neq 0} r^\Lambda (z)\leq C C_f \frac{1}{L^5 \log^2L}\frac{\log L}{\theta_L} \ , 
\eea
Here, $C_f$ is the Lipschitz constant of $f$. The last inequality follows from \eqref{unif_rates} in Lemma \ref{unibound} and equation \eqref{rabo}, derived independently from capacity estimates in Section \ref{sec:cap}. The latter implies that the espected total exit rate
\be\label{heqbound}
\theta_L \sum_{z\neq 0}r^\Lambda (z)
\leq C(1+\epsu ) \log L, \quad\mbox{with}\quad \epsu\to 0\quad \mbox{ as }\quad L\to\infty\ .
\ee
With Lemmas  Lemma \ref{ergbound} and \ref{trelest}, this implies that
\be
\E_{\mu^\Ecal} \bigg[ \Big(\int_0^{t\theta_L} h_f (\eta(s))\, ds\Big)^2 \bigg]\leq C C_f\, t\theta_L\,\frac{1}{L^5\log^2L}\, L^4\,\frac{\log L}{\theta_L}=CC_f t\, \frac{1}{L\log L}\ ,
\ee
with initial condition under the stationary distribution $\mu^\Ecal$. Now
\begin{equation*}
\sup_{x\in\Lambda} \E_{\mu^x} \bigg[ \Big(\int_0^{t\theta_L} \hspace{-.3cm} h_f (\eta(s))\, ds\Big)^2 \bigg] \leq \sum_{x\in\Lambda}\E_{\mu^x} \bigg[ \Big(\int_0^{t\theta_L}
\hspace{-.3cm} h_f (\eta(s))\, ds\Big)^2 \bigg] =L\E_{\mu^\Ecal} \bigg[ \Big(\int_0^{t\theta_L}\hspace{-.3cm} h_f (\eta(s))\, ds\Big)^2 \bigg]\ ,
\end{equation*}
and therefore
\be
\sup_{x\in\Lambda} \E_{\mu^x} \bigg[ \Big(\int_0^{t\theta_L} h_f (\eta(s))\, ds\Big)^2 \bigg]\leq C C_f t\, \frac{1}{\log L}\ .
\ee
Finally, we use Lemma \ref{unibound} and our estimate on the mixing time in Lemma \ref{trelest} to get for $\epsilon$ of order $1/\theta_L$
\bea
\lefteqn{\E_\eta \bigg|\int_0^{t\theta_L} h_f (\eta(s))\, ds\bigg|\leq \E_\eta \bigg|\int_0^{\tm (\epsilon )} h_f (\eta(s))\, ds\bigg| +\E_\eta \bigg|\int_{\tm (\epsilon )}^{t\theta_L} h_f (\eta(s))\, ds\bigg| }\nonumber\\
& &\leq  C_f \Big\| \sum_{z\neq 0} r^\Ecal (\cdot ,\Ecal^z )\Big\|_\infty \big( \tm (\epsilon )+\epsilon t\theta_L \big) +\sup_{x\in\Lambda} \bigg(\E_{\mu^x} \Big| \int_0^{t\theta_L} h_f (\eta(s))\, ds\Big|^2 \bigg)^{1/2}\nonumber\\
& &\leq CC_f \Big( \big( L^5 +t\big) \frac{1}{L^5\log^2L} +\sqrt{t/\log L}\Big) \longrightarrow 0 \quad\mbox{as}\quad L \to \infty.
\eea
This finishes the proof of Proposition \ref{equilibration}.\CQFD

\subsection{Proof of Lemma \ref{trelest} \label{trelsub}}

We will derive an upper bound on the relaxation time of the restricted process in a well by proving a Poincar\'e inequality, due to translation invariance it is enough to focus on the well $\Ecal^0$. We will use that the stationary measure $\mu^0$ outside the condensate location is essentially given by a product measure at the critical density in the limit $L\to\infty$. With $\mu [\Ecal^0 ]=\mu [\Ecal ]/L$ we can write
\be\label{envmeas}
\mu^0 [\eta ]=\mu^0 \big[ ( \eta_0 ,\eta_{\Lambda\setminus 0} )\big] =\nu^{L-1} [\eta_{\Lambda\setminus 0} ]\frac{\nu^1 [\eta_0 ]}{\nu^L [S_L =N]} \frac{L}{\mu [\Ecal ]}\ ,
\ee
where we use the notation $\eta =(\eta_0 ,\eta_{\Lambda\setminus 0} )$ for $\eta\in\Ecal^0$ to indicate the condensate size $\eta_0$ and the bulk configuration $\eta_{\Lambda\setminus 0}$ outside the condensate. For fixed particle number $N$ (i.e. under the measure $\mu$), $\eta_{\Lambda\setminus 0}$ obviously uniquely determines $\eta_0$. To simplify notation we identify measures with their mass function, and to avoid possible confusion in this section we will indicate the dimension of the product measure $\nu$ in each term with a superscript.\\

We interpret the product measure $\nu^{L-1}$ in (\ref{envmeas}) as the stationary measure of an auxiliary system of $L-1$ independent birth death chains with birth rate $1$ and death rate $g(n)$. Corresponding to our definition of the metastable wells in (\ref{valley}), we restrict the state space of the chains to $X =\{ 0,\ldots ,\beta_L\}$. Each chain has therefore the modified stationary measure
\be\label{barnu}
\bar\nu^1 =\nu^1 [\, .\, |X]\quad\mbox{and}\quad \bar\nu^1 [n]=\frac{\nu^1 [n]}{\nu^1 [X]}\quad \mbox{for all }n\in X\ .
\ee
Note that the state space of the auxiliary chains $X^{L-1}$ contains in particular all bulk configurations $\eta_{\Lambda\setminus 0}$ for $\eta\in\Ecal^0$, and recall that the product measure $\nu$ is stationary for the zero-range process.

\begin{lemma}\label{variance}
There exists a $C>0$ such that for all $f:X^{L-1}\to\R$ and all $L$ large enough
\be\label{varone}
\var_{\mu^0} (f)\leq C\nu^{L-1} [X^{L-1}] \var_{\bar\nu^{L-1}} (f\,\1_{B^0})\ .
\ee
Here
\be\label{b0}
B^0 =\Big\{ \zeta\in X^{L-1} :\Big( N-\sum_{x=1}^{L-1} \zeta_x ,\zeta\Big)\in\Ecal^0 \Big\} =\Big\{ \zeta\in X^{L-1} :\sum_{x=1}^{L-1} \zeta_x \leq \rho_c L+\alpha_L \Big\}
\ee
is the subset of configurations in $X^{L-1}$ which are compatible bulk configurations in the well $\Ecal^0$. For $\eta\in\Ecal^0$ we use the obvious extension $f(\eta ):=f(\eta_{\Lambda\setminus 0})$ with a slight abuse of notation, and the claim holds in particular for all functions defined on $\Ecal^0$.
\end{lemma}

\noindent\textbf{Proof.} Since $\mu^0 (f^2 ) -\mu^0 (f)^2 =\inf_c \mu^0 \big( (f-c)^2 \big)$, using the notation $\bar f=f-\bar\nu^{L-1} (f)$ we have with (\ref{envmeas})
\be
\var_{\mu^0} (f)\leq \mu^0 (\bar f^2 ) =\frac{L}{\mu [\Ecal ]}\sum_{\eta \in\Ecal^0} \frac{\nu^{L-1} [\eta_{\Lambda\setminus 0}]\,\nu^1 [\eta_0 ]}{\nu^L [S_L =N]}\, {\bar f}^2(\eta_{\Lambda\setminus 0} ) \ .
\ee
Following a result by Doney \cite{doney01}, $\nu^L [S_L =N]=L\,\nu^1 \big[ N-[\rho_c L]\big] \big( 1+o(1)\big)$ as $L,N\to\infty$ and $N/L\to\rho >\rho_c$, and therefore for $L$ large enough there exists $C>0$ such that
\bea\label{varbound}
\var_{\mu^0} (f)&\leq &C\sum_{\eta \in\Ecal^0} \frac{\nu^{L-1} [\eta_{\Lambda\setminus 0}]\,\nu^1 [\eta_0 ]}{\nu^1 \big[ N-[\rho_c L]\big]}\, {\bar f}^2(\eta_{\Lambda\setminus 0}) \nonumber\\
&\leq &C\,\frac{\nu^1 \big[N-[\rho_c L]-[\alpha_L ]\big]}{\nu^1 \big[N-[\rho_c L]\big]} \sum_{\eta \in\Ecal^0}\nu^{L-1} [\eta_{\Lambda\setminus 0}]{\bar f}^2 (\eta_{\Lambda\setminus 0} )\ ,
\eea
where we have also used $\mu [\Ecal ]\to 1$ and monotonicity of $n\mapsto\nu^1 [n]$. Since $\nu$ has power-law tails (\ref{tail}) the first ratio converges to $1$, and with the notation (\ref{b0}) we get for large enough $L$
\be
\var_{\mu^0} (f)\leq C\sum_{\zeta\in X^{L-1}} \nu^{L-1} [\zeta ]\bar f^2 (\zeta )\1_{B^0} (\zeta )\ ,
\ee
which finishes the proof with the definition of $\bar\nu$ (\ref{barnu}).\CQFD\\

The Dirichlet form of a single birth-death chain is given by
\bea\label{dsingle}
\Dcal (f)&=&\frac12 \sum_{n=0}^{\beta_L } \bar\nu^1 [n] \big[ g(n)\big( f(n-1)-f(n)\big)^2 +\1_X (n+1)\big( f(n+1)-f(n)\big)^2\big] \nonumber\\
&=&\sum_{n=0}^{\beta_L -1} \bar\nu^1 [n] \big( f(n+1)-f(n)\big)^2 \ ,
\eea
where we have used $\bar\nu [n]\, g(n) =\bar\nu [n-1]$ for $n\geq 1$ (\ref{tail}) and $g(0)=0$. The Dirichlet form for $L-1$ independent chains is therefore given by
\be\label{dind}
\Dcal_{L-1} (f)=\sum_{x=1}^{L-1}\sum_{\zeta\in X^{L-1}} \bar\nu^{L-1} [\zeta ] \big( f(\zeta^{x} )-f(\zeta )\big)^2 \,\1_X (\zeta_x +1)\ ,
\ee
where we write $\zeta^{x}$ for the configuration where a particle is added to $\zeta$ at site $x$, $\zeta^x_z=\zeta_z+\delta_{z,x}$.

\begin{lemma}\label{equichains}
For all $L>1$ we have a Poincar\'e inequality for $L-1$ independent chains, i.e. for all $f:X^{L-1}\to\R$
\be\label{equione}
\var_{\bar\nu^{L-1}} (f)\leq \frac{\beta_L^2}{4}\Dcal_{L-1} (f)\ .
\ee
For functions $f=f\1_{B^0}$ concentrating on bulk configurations of $\Ecal^0$ we can make a stronger statement,
\bea\label{equitwo}
\var_{\bar\nu^{L-1}} (f\1_{B^0} )&\leq &\frac{\beta_L^2}{4} \sum_{x=1}^{L-1}\sum_{\zeta_x =0}^{\beta_L -1} \bar\nu^{L-1} [\zeta ] \big( f(\zeta^{x} )-f(\zeta )\big)^2 \1_{B^0} (\zeta^{x})
\nonumber\\
&:=&\frac{\beta_L^2}{4}  \Dcal_{L-1}^{B^0} (f)\ .
\eea
\end{lemma}

This implies that the relaxation time for the independent birth-death chains is bounded above by $\beta_L^2$. This might seem like a crude bound, but it can in fact be shown that the relaxation time scales like that of a symmetric random walk \cite{inprep}, even though the chains are driven to the origin and have stationary measure $\bar\nu^1$. So our upper bound is sharp in the scaling $\beta_L^2$, and prefactors are not important for us here.\\
Note that the stronger statement for bulk configurations ensures that in all terms the function $f$ is only evaluated on $B^0$, which is important to avoid contributions from the boundary of $\Ecal^0$ in the final estimate by the Dirichlet form of the restricted process in Lemma \ref{dirichlets} below.\\

\noindent\textbf{Proof.} We use a standard Efron-Stein estimate \cite{es1981,bt2012} to bound the variance for $L-1$ iid random variables. Writing $\Lambda' =\Lambda\setminus 0$ for the bulk sites $1$ to $L-1$ this is given by
\be\label{varbegi}
\var_{\bar\nu^{L-1}} (f) \leq \sum_{x=1}^{L-1} \sum_{\zeta\in X^{L-1}} \bar\nu^{L-1} [\zeta ] \big( f(\zeta )-f_x (\zeta_{\Lambda'\setminus x} )\big)^2 \ ,
\ee
where we may choose any measurable function $f_x :X^{L-2}\to\R$. To show the first statement (\ref{equione}) we can simply choose
\be
f_x (\zeta_{\Lambda'\setminus x} ):=f\big( (\zeta_{\Lambda'\setminus x} ,0)\big)\ ,
\ee
where as before $(\zeta_{\Lambda'\setminus x} ,0)$ denotes the configuration $\zeta$ where $\zeta_x$ is replaced by $0$. Using the Cauchy-Schwarz inequality
\be\label{csu}
\big( f(\zeta )-f(\zeta_{\Lambda'\setminus x} ,0)\big)^2 \leq \zeta_x \sum_{l=0}^{\zeta_x -1} \big( f(\zeta_{\Lambda'\setminus x} ,l+1)-f(\zeta_{\Lambda'\setminus x} ,l)\big)^2\ ,
\ee
this leads to
\be\label{csu2}
\var_{\bar\nu^{L-1}} (f) \leq\sum_{x=1}^{L-1} \sum_{\zeta_{\Lambda'\setminus x}\atop\in X^{L-2}} \bar\nu^{L-2} [\zeta_{\Lambda'\setminus x} ]\sum_{\zeta_x =0}^{\beta_L } \!\!\zeta_x \bar\nu^1 [\zeta_x ] \sum_{l=0}^{\zeta_x -1} \!\!\big( f(\zeta_{\Lambda'\setminus x} ,l{+}1){-}f(\zeta_{\Lambda'\setminus x} ,l)\big)^2\ .
\ee
Reordering the sum and using that $n\mapsto n\bar\nu^1 [n]$ is monotone decreasing, we get
\bea\label{csu3}
\var_{\bar\nu^{L-1}} (f)&\leq &\sum_{x=1}^{L-1} \sum_{\zeta_{\Lambda'\setminus x}\atop\in X^{L-2}} \bar\nu^{L-2} [\zeta_{\Lambda'\setminus x} ] \sum_{l=0}^{\beta_L-1} \!\!\big( f(\zeta_{\Lambda'\setminus x} ,l{+}1)-f(\zeta_{\Lambda'\setminus x} ,l)\big)^2 \sum_{\zeta_x =l+1}^{\beta_L} \!\!\zeta_x \bar\nu^1 [\zeta_x ]\nonumber\\
&\leq &\sum_{x=1}^{L-1} \sum_{\zeta_{\Lambda'\setminus x}\atop\in X^{L-2}} \bar\nu^{L-2} [\zeta_{\Lambda'\setminus x} ] \sum_{l=0}^{\beta_L-1} \!\!\big( f(\zeta_{\Lambda'\setminus x} ,l{+}1)-f(\zeta_{\Lambda'\setminus x} ,l)\big)^2 \underbrace{(\beta_L {-}l)\, l}_{\leq\beta_L^2 /4}\,\bar\nu^1 [l]\nonumber\\
&\leq &\frac{\beta_L^2}{4} \Dcal_{L-1} (f)\ .
\eea
The same argument works when we restrict to functions $f\1_{B^0}$. Note that $\zeta\in B^0$ implies $(\zeta_{\Lambda'\setminus x} ,0)\in B^0$, and also all configurations appearing in the Cauchy-Schwarz decomposition in (\ref{csu}) are in $B^0$. Therefore, restricting the sum in (\ref{varbegi}) to $\zeta\in B^0$ leads to (\ref{equitwo}) with a completely analogous computation, which finishes the proof.\CQFD\\

To finish the proof of Lemma \ref{trelest} we will bound the Dirichlet form of $L-1$ independent birth-death chains restricted to $B^0$ by the Dirichlet form of the restricted process using a standard path counting argument. The bounds we get here are certainly not optimal, and one of the reasons for our conditions on the parameter $b$.

\begin{lemma}\label{dirichlets}
There exists $C>0$ such that for all $f:B^0\to\R$ (cf.\ (\ref{b0})) and all $L$ large enough
\be
\Dcal_{L-1}^{B^0} (f)\leq C\frac{1}{\nu [X^{L-1} ]}\, \frac{L^4}{\beta_L^2} \Dcal^0 (f)\ ,
\ee
where $\Dcal^0 (f)$ is the Dirichlet form of the restricted process
\be
\Dcal^0 (f)=\frac12\sum_{\eta ,\xi\in\Ecal^0} \mu^0 [\eta ] r^0 (\eta ,\xi )\big( f(\xi )-f(\eta )\big)^2 \ .
\ee
\end{lemma}

\noindent\textbf{Proof.} For each $f:B^0 \to\R$ we write $f(\eta ):=f(\eta_{\Lambda\setminus 0} )$ for its unique extension to $\eta\in\Ecal^0$ as before. The Dirichlet form of the restricted process is simply given by
\bea
\Dcal^0 (f)=\frac12\sum_{\eta\in\Ecal^0} \mu^0 [\eta ] \sum_{x\in\Lambda} \frac{g(\eta_x )}{2}&\Big( &\big( f(\eta^{x,x+1})-f(\eta )\big)^2 \1_{\Ecal^0} (\eta^{x,x+1})\nonumber\\
& &+\big( f(\eta^{x,x-1})-f(\eta )\big)^2 \1_{\Ecal^0} (\eta^{x,x-1})\Big)\ ,
\eea
since all jumps leading outside $\Ecal^0$ are suppressed. We change the summation to the set
\be
\Ecal^0_{N-1} =\big\{\eta\in X^L :\eta^z \in\Ecal^0 \mbox{ for some }z\in\Lambda\big\}\ ,
\ee
and use $\mu_{L,N}^0 [\eta^x ]g(\eta_x +1)=\mu_{L,N-1}^0 [\eta ]$ with the canonical measure for $N-1$ particles, which follows from (\ref{tail}). The Dirichlet form can then be written as
\be
\Dcal^0 (f)=\frac12\sum_{\eta\in\Ecal^0_{N-1}} \mu_{L,N-1}^0 [\eta ] \sum_{x\in\Lambda} \big( f(\eta^{x+1})-f(\eta^x )\big)^2 \1_{\Ecal^0} (\eta^{x})\1_{\Ecal^0} (\eta^{x+1})\ .
\ee
With the above notation the restricted Dirichlet form (\ref{equitwo}) of the independent chains can be written as
\bea
\Dcal_{L-1}^{B^0} (f)&=&\sum_{\eta \in \Ecal^0}\sum_{x\in\Lambda} \bar\nu^{L-1} [\eta_{\Lambda\setminus 0} ] \big( f(\eta^{0,x})-f(\eta )\big)^2 \1_{\Ecal^0} (\eta^{0,x} ) \nonumber\\
&=&\sum_{\eta \in \Ecal^0_{N-1}}\sum_{x\in\Lambda} \bar\nu^{L-1} [\eta_{\Lambda\setminus 0} ] \big( f(\eta^{x})-f(\eta^0 )\big)^2 \1_{\Ecal^0} (\eta^{x} )\1_{\Ecal^0} (\eta^{0} )\ ,
\label{chaindir}
\eea
where we used the same change of summation variable as above in the second line, and the fact that $\bar\nu^{L-1} [(\eta^0 )_{\Lambda\setminus 0} ]=\bar\nu^{L-1} [\eta_{\Lambda\setminus 0} ]$. We decompose the transport of a particle from the condensate site $0$ to $x$ into nearest neighbour jumps and use 
the Cauchy-Schwarz inequality on the telescoping sum to get
\be\label{cssum}
\big( f(\eta^{x})-f(\eta^0 )\big)^2 \leq L\,\sum_{y=0}^{x-1} \big( f(\eta^{y+1})-f(\eta^{y} )\big)^2
\ee
for all $x$, since the longest path of a particle is clearly bounded by $L$. We can bound every such term in (\ref{chaindir}) that way, and as long as $\eta_y +1\leq\beta_L$ for all $y=1,\ldots x-1$, all $\eta^{y} \in \Ecal^0$ and the terms in the sum (\ref{cssum}) correspond to 'allowed' transitions that appear also in $\Dcal^0 (f)$. The 'flow' of an allowed transition is the number of times it  appears in (\ref{chaindir}) using (\ref{cssum}), and summing over all target positions $x$ this is bounded by $L$.\\

On the other hand, if there exist sites $0<y_1 <\ldots <y_m <x$, $m>0$, with $\eta_{y_i} =\beta_L$, the generic path in (\ref{chaindir}) contains non-allowed transitions and has to be re-routed, increasing the flow of certain allowed transitions. 
To bound this increase, we introduce the notation
\be
\sigma^y_z \eta :=\eta +\delta_y -\delta_z \quad\mbox{for }y,z\in\Lambda\ ,
\ee
where $\delta_y$ is the configuration with a single particle at $y$. The path corresponding to the sum (\ref{cssum}) can then be represented by the equation
\[
\eta^x =\sigma^x_{x-1} \cdots \sigma^2_1\sigma^1_0\eta^0 \ .
\]
If there is an isolated site $y$ with $\eta_y =\beta_L$ and $\eta_{y\pm 1} <\beta_L$, we re-route the path of a particle from $y-1$ to $y+1$ from
\[
\eta^{y+1} =\sigma^{y+1}_y \sigma^y_{y-1} \eta^{y-1}\quad\mbox{to}\quad \eta^{y+1} = \sigma^y_{y-1} \sigma^{y+1}_y \eta^{y-1}\ .
\]
Instead of moving the particle to site $y$, it remains in site $y-1$ and a particle moves from $y$ to $y+1$ first. In the next step the particle follows from site $y-1$ to $y$ reaching $\eta^{y+1}$ only via allowed transitions. If there is a block of $k$ consecutive sites with $\eta_{y} =\ldots =\eta_{y+k-1} =\beta_L$, the re-routed path of a particle from $y-1$ to $y+k$ along valid transitions is
\be\label{ppath}
\eta^{y+k} =\sigma^y_{y-1} \cdots\sigma^{y+k-1}_{y+k-2} \sigma^{y+k}_{y+k-1} \eta^{y-1} \ .
\ee
Possibly combining re-routing over several blocks of sites with occupation number $\beta_L$, we associate a unique particle path from $0$ to $x$ to each base configuration $\eta\in X_{L,N-1}$, using only allowed transitions. The flow of a transition is then multiplied by the number of associated base configurations that use it for some $x\in\Lambda$, and every transition with multiplicity higher than $1$ involves at least one site with occupation number $\beta_L$. Denote by $\zeta\to\zeta'$ one of the transitions along the path in (\ref{ppath}), then one associated base configuration is obviously $\eta$, and another one is given by the minimal configuration
\[
\zeta\wedge\zeta' :=(\zeta_z \wedge\zeta'_z :z\in\Lambda )\ .
\]
Note that for all transitions $\zeta\to\zeta'$ along the path in (\ref{ppath}) we have for the maximal configuration
\[
\zeta\vee\zeta' :=(\zeta_z \vee\zeta'_z :z\in\Lambda )=\eta^{y-1,y+k} =\eta +\delta_{y-1}+\delta_{y+k}\ .
\]
It is easy to see that any possible base configuration associated to a transition $\zeta\to\zeta'$ in (\ref{ppath}) has to be of the form
\[
\eta^{y-1,y+k} -\delta_{z}-\delta_{z'} \quad\mbox{for some }z,z'\in\{ y-1,\ldots ,k+1\}\ .
\]
In many cases not all of those base configurations contribute (or are even in $\Ecal$), but this provides an upper bound of $(k+1)^2$ for the flow multiplicity of transitions $\zeta\to\zeta'$ along the path (\ref{ppath}).\\

For any base configuration $\eta$ summed over in (\ref{chaindir}) there are at most of order $L/\beta_L$ sites with occupation number $\beta_L$, and therefore, the multiplicity of the flow along any allowed transition is bounded by $C(L/\beta_L )^2$. Together with the generic flow bound of order $L$ and (\ref{cssum}) this implies that 
\[
\Dcal_{L-1}^{B^0} (f)\leq C\,\frac{L^4}{\beta_L^2} \sum_{\eta\in\Ecal^0_{N-1}} \bar\nu^{L-1} [\eta_{\Lambda\setminus 0} ] \sum_{y\in\Lambda} \big( f(\eta^{y+1})-f(\eta^y )\big)^2 \1_{\Ecal^0} (\eta^{y})\1_{\Ecal^0} (\eta^{y+1})\ .
\]
Using again (\ref{envmeas}) and the same approach as in the proof of Lemma \ref{variance} we can bound
\bea
\bar\nu^{L-1} [\eta_{\Lambda\setminus 0} ]&=&\frac{\nu^{L-1} [\eta_{\Lambda\setminus 0}] }{\nu^{L-1} [X^{L-1}]}= \frac{1}{\nu^{L-1} [X^{L-1}]}\,\mu_{L,N-1}^0 [\eta ]\, \frac{\nu^L [S_L =N{-}1]\,\mu_{L,N-1} [\Ecal ]}{L\,\nu^1 [\eta_0 ]}\nonumber\\
&\leq &C\frac{1}{\nu^{L-1} [X^{L-1}]}\,\mu_{L,N-1}^0 [\eta ]\,\frac{\nu^1 \big[ N-[\rho_c L]\big]}{\nu^1 [N]}
\eea
for a suitable $C>0$. We have used again that $\nu^L [S_L =N]=L\,\nu^1 \big[ N-[\rho_c L]\big]\big( 1+o(1)\big)$ \cite{doney01}, monotonicity of $n\mapsto\nu^1 [n]$ and $\mu_{L,N-1} [\Ecal ]\to 1$. Since
\[
\frac{\nu^1 [N]}{\nu^1 \big[ N-[\rho_c L]\big]}=\Big(\frac{\rho -\rho_c}{\rho}\Big)^b \big( 1+o(1)\big)\ ,
\]
we end up with
\[
\Dcal_{L-1}^{B^0} (f)\leq C \frac{L^4}{\beta_L^2\nu^{L-1} [X^{L-1}]}\sum_{\eta\in\Ecal^0_{N-1}} \mu_{L,N-1}^0 [\eta_{\Lambda\setminus 0} ] \sum_{x\in\Lambda} \big( f(\eta^{x+1})-f(\eta^x )\big)^2 \1_{\Ecal^0} (\eta^{x})\1_{\Ecal^0} (\eta^{x+1})\ ,
\]
which finishes the proof of Lemma \ref{dirichlets}.\CQFD\\

Together with Lemmas \ref{variance} and \ref{equichains} we get a Poincar\'e inequality for the restricted process, i.e. for all $f:\Ecal^0 \to\R$ there exists $C>0$ such that
\be
\var_{\mu^0} (f)\leq C L^4 \Dcal^0 (f)\ .
\ee
Therefore, the relaxation time for the restricted process is bounded by $\tr \leq C L^4$ on each well $\Ecal^x$, and by a standard result \cite{peres} this implies for the $\epsilon$-mixing time that
\be
\tm (\epsilon )\leq -\tr \log \big(\epsilon \min_{\eta\in\Ecal^0} \mu^0 [\eta ]\big) \leq C  L^5 \big(1+ L^{-1}\log (1/\epsilon )\big)\ ,
\ee
since $\mu^0 [\eta ]$ is at most exponentially small in $L$. This finishes the proof of Lemma \ref{trelest}. \CQFD\\

\section{Uniform bounds on exit rates via coupling\label{sec:coup}}

To derive the uniform bound of Lemma \ref{unibound} on the exit rate out of a well for the trace process, we will construct a coupling of the zero-range process with a growing number of birth-death chains. The number of chains increases only linearly in time, and the coupling ensures that in the event of changing well, at least one of the chains has grown a condensate, the probability of which can be controlled directly from metastability results on a fixed size lattice in \cite{beltranetal09}. We will need additional control on how much time the full process spends outside the well $\Ecal^0$, which we achieve by using mixing estimates on larger wells containing the original ones. This is derived first in the last subsection, together with a proof of Proposition \ref{trace} on substitution by the trace process.

\subsection{Construction of the coupling}
\label{coupling}

We construct the coupling that will be used in the next subsection to prove uniform bounds on exit rates from wells. Let $(\zeta (t):t\geq 0)$ be a birth-death chain on the state space $X=\{ 0,1,\ldots\}$ with birth/arrival rate $1$ and death/departure rate $g(\zeta )$ as given in (\ref{rates}), characterized by the generator
\be\label{bdgen}
\Lcal f(\zeta )=\big( f(\zeta +1)-f(\zeta )\big) +g(\zeta )\big( f(\zeta -1)-f(\zeta )\big)\ .
\ee
Note that the boundary condition at $\zeta =0$ is included with $g(0)=0$, and this chain has stationary measure $\nu$ as given in (\ref{tail}). 
For some fixed $\varepsilon \in (0,\rho -\rho_c )$ denote by
\be\label{hitlev}
y_L :=(\rho -\rho_c -\varepsilon )L\quad\mbox{and}\quad T_{y_L} :=\inf\{ t\geq 0:\zeta (t)\geq y_L \}
\ee
a size-dependent level that has to be crossed to grow a condensate, and the associated hitting time. 
Lower bounds for the hitting time $T_{y_L}$ are typically of order $\theta_L =L^{1+b}$ and can be derived by direct computation.

\begin{lemma}\label{twosite}
There exist constants $C_1 ,C_2 >0$, such that for all initial conditions $\zeta (0)=\zeta_0 \in \{ 0,1,\ldots ,B_L \}$ with $1\ll B_L\ll L$, we have that
\be
\E_{\zeta_0} [T_{y_L} ]\geq C_1 \theta_L \quad\mbox{and}\quad \P_{\zeta_0} [T_{y_L} \leq t]\leq C_2 \,\frac{tB_L^{b-1}}{\theta_L}\ .
\ee
\end{lemma}

\noindent\textbf{Proof.} It is easy to show that the expected hitting time $\tau_y^x =\E_x [T_y ]$ with $x<y\in\N$ of a birth-death chain with birth rates $h(\zeta )$, death rates $g(\zeta )$, and stationary measure $\nu$ is given by
\be\label{bdform}
\tau_y^x =\E_x [T_y ]=\sum_{\zeta =x}^{y-1} \frac{1}{h(\zeta )\nu[\zeta ]}\sum_{n=0}^\zeta \nu [n]\ .
\ee
For a reference see e.g. \cite{inprep}. We have $h(\zeta )=1$ and due to monotonicity of $\nu$ in this case we can use simple integral bounds for sums. We get
\be\label{firstsum}
\sum_{n=0}^\zeta \nu [n]\geq \frac{1}{z_c} \int_1^\zeta u^{-b}\, du =\frac{1-\zeta^{1-b}}{z_c (b-1)}\ ,
\ee
which then analogously leads to
\be\label{secondsum}
\tau_y^x \geq \sum_{\zeta =x}^{y-1} \frac{\zeta^b -\zeta}{b-1} \geq \frac{(y-1)^{b+1}-x^{b+1}}{(b-1)(b+1)}-\frac{y^2 -x^2}{2(b-1)}\ .
\ee
With $x\leq B_L \ll y=y_L \sim L$ this directly implies the first statement.\\

To derive the second statement, for a given $\zeta_0$ we couple the chain with a modified chain $\zeta'(\cdot )$ that cannot jump below $\zeta_0$, i.e. it has death rates
\[
g'(\zeta) =g(\zeta )\mbox{ for }\zeta >\zeta_0\quad\mbox{and}\quad g'(\zeta_0 )=0\ .
\]
It is clear that the chain $\zeta'(\cdot )$ will reach $y_L$ before the original one, so its hitting time $T'_{y_L}$ will provide a lower bound for $T_{y_L}$. Furthermore, the lowest hitting time is clearly achieved for $\zeta_0 =B_L$ and we can focus on this case. Since the point $\zeta_0 =B_L$ is the left end of the state space for the $\zeta'(\cdot)$ chain, the Markov property implies the following sub-multiplicity,
\[
\P'_{B_L} [ T'_{y_L} > s] \leq \P'_{B_L} [ T'_{y_L} > t ] ^{ [s/t] }\ .
\]
Integrating over $s$ and re-arranging yields
\[
\P'_{B_L} [ T'_{y_L} \leq t ] \leq t / \E'_{B_L} [T'_{y_L} ]\ ,
\]
and it remains to estimate the expectation of $T'_{y_L}$ from below. Note that the chain $\zeta'(\cdot)$ has the same stationary measure $\nu$ restricted to $\zeta\geq B_L$, changing the normalization to $z'_c$. The latter cancels in (\ref{bdform}), and we simply have to adapt (\ref{firstsum}) as
\[
\sum_{n=B_L}^\zeta \nu [n]\geq \frac{1}{z'_c} \int_{B_L}^\zeta u^{-b}\, du =\frac{B_L^{1-b}-\zeta^{1-b}}{z'_c (b-1)}\ .
\]
Analogously to (\ref{secondsum}) this implies
\[
\E'_{B_L} [T'_{y_L} ]\geq \sum_{\zeta =B_L}^{y_L-1} \frac{B_L^{1-b}\zeta^b -\zeta}{b-1}\geq C_2 B_L^{1-b} y_L^{b+1} \ ,
\]
for a suitable constant $C_2 >0$, finishing the proof.\CQFD\\

Now let us fix a configuration $\eta\in\Ecal^0$. For the original process $(\eta(t):t\geq 0)$ the occupation numbers outside the condensate $\eta_x (t)$, $x=\{ 1,\ldots ,L-1\}$ are birth-death processes with Markovian departure processes at rate $g(\eta_x )$, but non-Markovian arrival processes that depend on the neighbouring occupation numbers for each site. Conditioned on the configuration $\eta (t)$ at time $t$, the arrival rate at site $x$ is given by
\be
\label{arrival_rate}
a_x (\eta (t) ):=\frac{1}{2}\, g\big(\eta_{x-1}(t)\big)+\frac12 \, g\big(\eta_{x+1}(t)\big),
\ee
and if both neighbouring sites are occupied this rate can be as large as $2^b$. In order to dominate $\eta_x (t)$ by a Markovian birth-death chain with arrival rate $1$ to apply Lemma \ref{twosite}, 
we couple it with an increasing number of chains. At any given time, at least one of those chains will dominate $\eta_x (t)$ and lead to an estimate for the probability of leaving the well $\Ecal^0$.\\

\begin{figure}
\begin{center}
\includegraphics[width=0.8\textwidth]{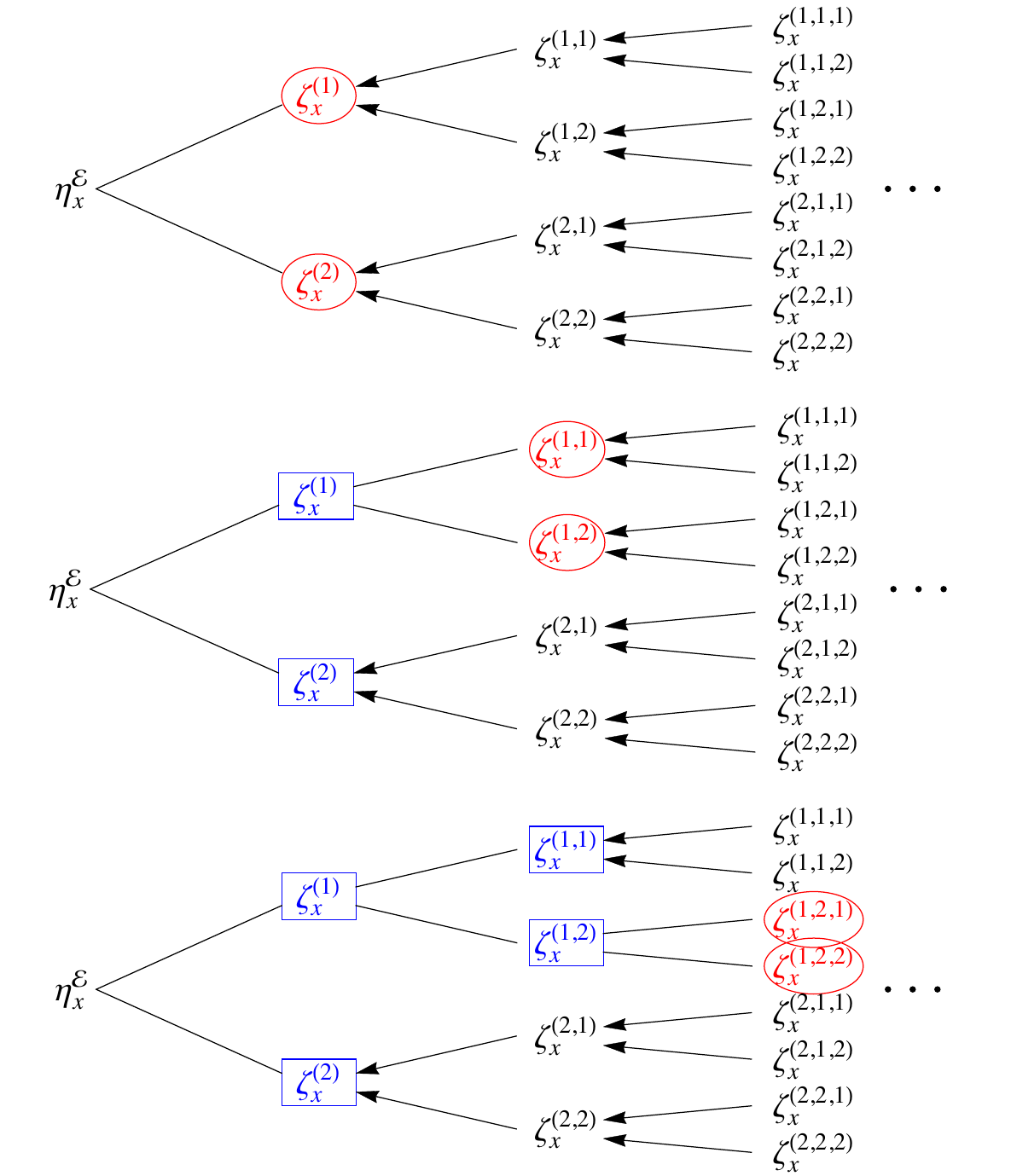}\\
\end{center}
\caption{\label{figcoup}
Coupling procedure illustrated for $m=2$, associated chains are shown in red and are encircled, independent chains are shown in blue in a rectangular box. Arrows indicate identical copies. Initially (top) only the chains in generation $n=1$ are associated. After the first particle arrives at site $x$ (middle) both chains in the first generation become independent, and the $m$ descendants of the the chain $\zeta_x^{(1)}$ get associated. This process is repeated after the second particle arrives (bottom).
}
\end{figure}

Let $m\geq 2^b$ be the smallest integer greater than or equal to  the maximal jump rate $g(2)=2^b$ (\ref{rates}) of the zero range process. The coupling as described below is applied for all times $t\geq 0$ and sites $x=1,\ldots ,L-1$, and is illustrated in Figure \ref{figcoup} for the simplest case $m=2$. 
To each site $x$ we associate an infinite number of birth-death chains $(\zeta^\kbo_x (t):t\geq 0)$ with generator (\ref{bdgen}), where we index the chains by vectors of variable length of the form $\kbo =(k_1 ,\ldots ,k_n )$ with $k_i \in \{ 1,\ldots ,m\}$ and $n=1,2,\ldots$. This corresponds to indexing the chains by the nodes of an $m$-ary regular tree $\Rcal_m$ without root, with generations indexed by $n$, and we write $\kbo\in\Rcal_m$. At any given time, each chain $\zeta^\kbo_x (t)$ in the tree  {\it 1)} is an identical copy of its unique parent chain, {\it 2)} evolves independently of $\eta_x (t)$ and all other chains, or {\it 3)} is associated to $\eta_x (t)$ as described below. The assignment of each chain to one of these three groups changes in time, and we denote by
\bea
\Ccal_x (t)&:= &\big\{ \kbo \in\Rcal_m :\zeta^\kbo_x (t)\mbox{ is associated to }\eta_x (t), \mbox{ and no ancestor of } \zeta_x^\kbo(t) \mbox{ is associated to } \eta_x (t) \big\}\nonumber\\
\Ical_x (t)&:= &\big\{ \kbo \in\Rcal_m :\zeta^\kbo_x (t)\mbox{ evolves independently}\big\}
\eea
the index sets of chains which are not identical copies of their parent. At any time $t\geq 0$, the number of  chains in $\Ccal_x(t)$ is $|\Ccal_x (t)|=m$, for all sites $x$.\\

Initially, we set the $m$ chains in generation $n=1$ equal to $\eta_x$, i.e. $\Ccal_x (0)=\{ (1),\ldots ,(m)\}$, $\Ical_x =\emptyset$ and all other chains are identical copies of their parent. We use identical initial conditions, that is, for each site $x$
\be
\zeta^{(1)}_x (0)=\ldots =\zeta^{(m)}_x (0)=\eta_x (0)\in \{ 0,1,\ldots ,\beta_L \}\ ,
\ee
and our coupling will ensure that $\eta_x (t)\leq \zeta^\kbo_x (t)$ for all $\kbo\in \Ccal_x (t)$. 
For the departure process of associated chains with $\kbo\in\Ccal_x (t)$ we simply use a basic coupling for all $x$ and $t\geq 0$, i.e. particles in $\zeta^\kbo_x (t)$ leave together with particles in $\eta_x (t)$ with probability $g\big(\zeta^\kbo_x (t) \big) /g\big(\eta_x (t)\big)\leq 1$ for $\eta_x (t)>1$, and they additionally leave, independently of particles in $\eta_x(t)$, at rate $g\big(\zeta^\kbo_x(t)\big) -g\big(\eta_x (t)\big)$ in case this quantity is positive for $\eta_x (t)\leq 1$. 
Note that the departure dynamics preserve the order
\be\label{order}
\eta_x (t)\leq\zeta^\kbo_x (t)\quad\mbox{for all }\kbo\in\Ccal_x (t) \mbox{ and } t\geq 0\ ,
\ee
and we will couple the arrival processes in such a way that this is true also for the full process. 
To achieve this, we change the structure summarized by the sets $\Ccal_x (t)$ and $\Ical_x (t)$ at every jump event on the arrival site. When a particle arrives at site $x$ in the 
$\eta$-process at time $t$ we pick one of the $m$ chains in $\Ccal_x (t)$ uniformly at random, add a particle to all of its $m$ children (which up to this point have evolved as identical copies of their parent), and disassociate the chains in $\Ccal_x(t-)$ so that from this time on they run independently of $\eta_x(s), s\ge t$. That is, sample $ \kbo^* \sim U(\Ccal_x (t))$, and let
\begin{equation*}
\Ccal_x (t)=\big\{ (\kbo^*,1) ,(\kbo^* ,2),\ldots ,(\kbo^* ,m)\big\}\ ,\quad \Ical_x(t)=\Ical_x(t-)\cup \Ccal_x(t-)\ .
\end{equation*}
So far the coupling leads to an effective arrival rate of $a_x (\eta(t) )/m\leq 1$, $a_x(\eta(t))$ as in \eqref{arrival_rate}, for each associated chain, which is typically strictly smaller than $1$. To each associated process we independently add particles at rate $1-a_x (\eta(t) )/m$, leading to a total arrival rate of $1$ as required. Therefore all chains $(\zeta^\kbo_x (t):t\geq 0)$, $\kbo\in\Rcal_m$, have the desired marginal dynamics of a birth-death chain with generator (\ref{bdgen}). Note that the total number of particles in the associated chains is not conserved and is growing in time, but the main point is that the coupling fulfills (\ref{order}).\\

This coupling construction leads to increasing sets $\Ical_x (t)$ of independently evolving chains, but at any time there is only a finite number of chains which are not identical copies of their parent. This number grows only linearly in time with high probability, and we will use this in the next subsection to prove a uniform bound on the exit rate from a well.

\subsection{Proof of Lemma \ref{unibound}\label{proofER}}

For the trace process $(\eta^\Ecal (s):s\geq 0)$, $\eta^\Ecal (0)=\eta(0)\in\Ecal^0$, to reach another well before time $t>0$, one of the occupation numbers for $x\neq 0$ has to grow larger than $y_L =(\rho-\rho_c -\varepsilon )L$, $\varepsilon \in (0,\rho -\rho_c )$ given in (\ref{hitlev}), before the full process $(\eta(s): s\geq 0)$ has spent a total of time $t$ in the well $\Ecal^0$. By construction of the above coupling, this implies that at least one of the associated chains $\zeta^\kbo_x $ with $\kbo\in\Rcal_m$ and $x\neq 0$ must also have reached this level. 
Again, it is important to note that for each $x\neq 0$ and $s\geq 0$ there is only the finite number of chains in $\Ccal_x (s)\cup\Ical_x (s)\subset\Rcal_m$ which are not identical copies of their parent. 
Since all the chains associated to site $x$ are birth-death chains with generator (\ref{bdgen}), we can use Lemma \ref{twosite} to estimate the time it takes for the chains to reach level $y_L$. \\

Note that the dynamics of the birth-death chain $(\zeta(t), t\ge 0)$ does not depend on the parameter $L$, hence Lemma \ref{twosite} in fact says that for any parameter 
$K\nearrow \infty$, letting $T_K:=\inf\{t\ge 0: \zeta(t)\ge K\}$ and for an initial condition $\zeta_0\le B_K$, $1\ll B_K \ll K$ we have
\be
\label{onechain}
\P_{\zeta_0}\big[ T_K\le t\big]\le C\,\frac{tB_K^{b-1}}{K^{1+b}}.
\ee

We apply this bound to gain control on the time it takes the process $(\eta(s): s\ge 0)$ to exit a larger well $\tilde{\Ecal}^0\supset \Ecal^0$, which we choose as
\begin{equation}
\label{phil}
\tilde{\Ecal}^0:=\big\{ \eta \in X_{L,N}: \eta_0\ge N-\rho_c L-\alpha_L,\, \eta_y\le \phi_L \mbox{ for all } y\neq 0 \big\}, \mbox{ with }\phi_L:=\frac{L}{\log L}\ .
\end{equation}
Clearly, the larger wells define an analogous partition of the state space that uniquely characterises the condensate location for sufficiently large $L$, so the definition of the restricted process \eqref{engen} in Section \ref{sec:equi} can be adapted to the extended wells. Lemma \ref{trelest} directly applies to $\tilde{\Ecal}^0$ replacing $\beta_L$ by $\phi_L$ in the derivation of all estimates, and therefore yields the following bounds on the relaxation time $\tilde{t}_{\mathrm{rel}}$ and the $\delta$-mixing time $\tilde{t}_{\mathrm{mix}}(\delta)$ of the restricted process in $\tilde{\Ecal}^0$,
\begin{equation}
\label{newmix}
\tilde{t}_{\mathrm{rel}}\le C L^4 \quad\mbox{and}\quad t_\delta :=\tilde{t}_{\mathrm{mix}}(\delta)\le C L^4\,\big(L+ \log(1/\delta)\big).
\end{equation}

Let us denote by $\tilde{T}$ the exit time of the process $(\eta(s): s\ge 0)$ from the set $\tilde{\Ecal}^0$ and for $\eta\in\tilde{\Ecal}^0$ let us write
\[
B(\eta)=\max_{x\neq 0}\eta_x.
\]
\begin{lemma}
\label{extexit}
There exist constants $C,\gamma>0$, such that for any $\eta\in\tilde{\Ecal}_0$ and $t\ge 0$
\begin{equation}
\P_\eta\big[\tilde{T}\le t\big]\le C\frac{Lt^2 B(\eta)^{b-1}}{\phi_L^{1+b}}+Le^{-\gamma t}.
\end{equation}
\end{lemma}
{\bf Proof.} On account of the coupling argument of Subsection \ref{coupling}, if we consider the system of birth-death chains $(\zeta_x^{\kbo}(s): s\ge 0, x\in \Lambda\setminus \{0\}, \kbo \in \Rcal_m)$ with initial condition associated to $\eta \in \tilde{\Ecal}^0$, so that in particular $\zeta^\kbo_x(0)\le B(\eta)$, $x\neq 0$, $\kbo \in \Ccal_x(0)$, then 
\[
\big\{\tilde{T}\le t\big\}\subseteq \big\{\exists\, \kbo \in \Rcal_m,\,x\neq 0 \mbox{ and } s\in [0,t] \mbox{ such that }\,\zeta_x^{\kbo}(s) \ge \phi_L\big\},
\]
so
\[
\P_\eta\big[ \tilde{T}\le t \big]  \le \P_{\zeta(\eta)}\big[\zeta_x^\kbo(s)\ge \phi_L \mbox{ for some } s\in [0, t],\,x\neq 0,\,\kbo \in \Ccal_x(s)\cup\Ical_x(s)  \big] .
\]
The number of chains in $\Ccal _x(t) \cup \,\Ical_x(t)$ increases with arrival processes of particles which are independent Poisson with bounded rates, so the probability that 
$|\Ccal_x (t )|+|\Ical_x (t )|> Ct$ decays exponentially in $t$ for $C$ large enough. Applying \eqref{onechain} to each of these chains independently with $K=\phi_L$ we obtain
\begin{align}
\label{schance}
\P_\eta\big[ \tilde{T}\le t \big] & \le C \frac{Lt^2 B(\eta)^{b-1}}{\phi_L^{1+b}} +\P_\eta \big[ |\Ccal_x (t)|+|\Ical_x (t )|> Ct \mbox{ for some }x\neq 0 \big] \notag\\
&\leq C \frac{Lt^2 B(\eta)^{b-1}}{\phi_L^{1+b}}+Le^{-\gamma t},
\end{align}
which completes the proof of the Lemma. \CQFD\\

We turn now to the proof \eqref{unif_rates} in the statement of Lemma \ref{unibound}. Let us first note that for any $t>0$
\begin{align}
\label{cota}
\P_{\eta}\big[\psi_L(\eta^\Ecal(s)) \neq 0 &\text{ for some } s \in [0,t] \big] 
\ge 
\Big(1-e^{-t\sum_{\xi \neq \eta} r^\Ecal(\eta,\xi)}\Big) \ \frac{\sum_{\xi \notin \Ecal^0 } r^\Ecal (\eta, \xi)}{\sum_{\xi \neq \eta} r^\Ecal(\eta, \xi)}, 
\end{align}
where the ratio on the right is the probability that the first jump out of the configuration $\eta$ is to a configuration in  the complement of $\Ecal^0$. Let now $t=\frac{1}{L}$. Since  
\[\sum_{\xi \neq \eta} r^\Ecal(\eta,\xi)\le \sum_{\xi \neq \eta} r(\eta,\xi)=\sum_{x\in \Lambda} g(\eta_x)\le mL\] we get that 
\begin{align}
\label{cota2}
\sum_{\xi \notin \Ecal^0 } r^\Ecal (\eta, \xi)\le CL\
\P_{\eta}\big[\psi_L(\eta^\Ecal(s)) \neq 0 &\text{ for some } s\le \frac{1}{L} \big],
\end{align}
where 
\[
C=\frac{m}{1-e^{-m}}=\sup_{x\in (0,m]}\frac{x}{1-e^{-x}}.
\]
We next prove that there exists a constant $C>0$ such that 
\be
\label{bunibound}
\sup_{\eta\in\Ecal^0} \P_\eta \big[\psi_L (\eta^\Ecal (s))\neq 0\mbox{ for some }s\le \frac{1}{L}\big] \leq C\,\frac{1}{L^6\log^2 L},
\ee
which together with (\ref{cota2}) immediately implies the assertion of Lemma \ref{unibound}. For this, write
\begin{align}
\big\{\psi_L(\eta^\Ecal(s))\neq 0 \mbox{ for some $s\le \frac{1}{L}$} \big\}&=\big\{\eta(s)\in \Ecal\setminus\Ecal^0\mbox{ for some } s\in [0,S_{L^{-1}}]\big\}\notag\\
&=\big\{\int_0^{T_{\Ecal\setminus\Ecal^0}}\mathbbm{1}_{\Ecal^0}\big(\eta(s)\big) ds\le \frac{1}{L}\big\},
\label{event1L}
\end{align}
where $S_t:=\sup\{s\ge 0: {\cal T}_s\le t\}$ is the inverse of the local time $\Tcal_s$ in $\Ecal ^0$ defined in Section \ref{trace_process} and $T_{\Ecal\setminus\Ecal^0}$ is the hitting time of the process $(\eta(s): s\ge 0)$ on $\Ecal\setminus\Ecal^0=\cup _{y\neq 0} \Ecal^y$. To control the probability of the event in (\ref{event1L}) we define an intermediate, deterministic time $t_L$ which is small enough for the process $\eta(\cdot )$ started at $\eta\in\Ecal^0$ to remain in $\tilde{\Ecal}^0$ by $t_L$ with high probability, but large enough for the restricted process on $\tilde{\Ecal}^0$, which we denote by $(\xi(s): s\ge 0)$, to have mixed. Precisely, we require
\be
\frac{t_L}{t_\delta}\ge CL^{\kappa} \mbox{ for some } \kappa>0\quad\text{and}\quad \sup_{\eta\in\Ecal^0}\,\P_\eta\big[\tilde{T}\le t_L\big]\le \frac{C}{L^6\log^2L} \qquad\text{ as } L\to\infty\ .
\label{tLprop}
\ee
By (\ref{newmix}) and Lemma \ref{extexit} these conditions are simultaneously satisfied if $t_L$ is chosen as
\be\label{tLdef}
t_L:=\frac{L^{(b-10)/2}}{(\log L)^\frac{b+5}{2}}.
\ee
In fact, the condition $b>20$ arises from the need that $t_L \gg t_\delta$ in (\ref{tLprop}). Define now
\begin{equation}
{\cal A}=\big\{\int_0^{\tilde{T}} \mathbbm{1}_{\Ecal^0}\big(\eta(s)\big)ds\le \frac{1}{L}\big\}\supseteq\big\{\psi_L(\eta^\Ecal(s))\neq 0 \mbox{ for some $s\le \frac{1}{L}$} \big\}.
\label{calAdef}
\end{equation}
For $\eta\in\tilde{\Ecal}^0$ we have
\begin{align}
\P_\eta\big[{\cal A}\big] &\le \P_\eta\big[\tilde{T}\le t_L\big]+\P_\eta\big[{\cal A}\cap\{\tilde{T}> t_L\}\big]\notag\\
&\le \P_\eta\big[\tilde{T}\le t_L\big]+\P_\eta\big[\int_{t_L}^{\tilde{T}} \mathbbm{1}_{\Ecal^0}\big(\eta(s)\big)ds\le \frac{1}{L};\ \tilde{T}> t_L\big]\notag\\
&\le \P_\eta\big[\tilde{T}\le t_L\big]+\E_\eta\Big[\P_{\eta(t_L)}\big[{\cal A}\big];\ \tilde{T}> t_L\Big],
\label{bootshalf}
\end{align}
where the third line follows from the Markov property. We may couple a restricted process to $\tilde{\Ecal}^0$, which we will denote by $(\xi(t): t\ge 0)$, to the original process, so that they jump together up to time $\tilde{T}$. If $\P_{(\eta,\eta)}^{\text{coup}}$ is the coupling measure with marginals $\P_\eta$ (original process) and $\tilde{\P}_\eta$ (restricted process), then
\begin{align*}
\E_\eta\Big[\P_{\eta(t_L)}\big[{\cal A}\big];\ \tilde{T}> t_L\Big]&=\E_{(\eta,\eta)}^{\text{coup}}\Big[\P_{\eta(t_L)}\big[{\cal A}\big];\ \tilde{T}> t_L\Big]=\E_{(\eta,\eta)}^{\text{coup}}\Big[\P_{\xi(t_L)}\big[{\cal A}\big];\ \tilde{T}> t_L\Big]\\&\le \E_{(\eta,\eta)}^{\text{coup}}\Big[\P_{\xi(t_L)}\big[{\cal A}\big]\Big]
=\tilde{\E}_\eta\Big[\P_{\xi(t_L)}\big[{\cal A}\big]\Big].
\end{align*}
We now substitute this estimate in (\ref{bootshalf}) to get that for any $\eta\in\tilde{\Ecal}^0$
\begin{align}
\P_\eta\big[{\cal A}\big] &\le \P_\eta\big[\tilde{T}\le t_L\big]+\tilde{\E}_\eta\Big[\P_{\xi(t_L)}\big[{\cal A}\big]\Big]\notag\\
&\le \P_\eta\big[\tilde{T}\le t_L\big]+ \P_{\tilde{\mu}^0}\big[{\cal A}\big]+d_L,
\label{bootstrapacal}
\end{align}
where $d_L$ is the total variation distance between the distribution of $\xi({t_L})$ and the invariant measure $\tilde{\mu}^0$ of the process $(\xi(s): s\ge 0)$. By \eqref{rest_meas}, $\tilde{\mu}^0$ is simply the invariant distribution $\mu$ of $\eta(\cdot)$ restricted to $\tilde{\Ecal}^0$, that is 
\[
\tilde{\mu}^0\big[\xi]=\mathbbm{1}_{\tilde{\Ecal}^0}(\xi)\frac{\mu[\xi]}{\mu\big[\tilde{\Ecal}^0\big]}
\]
and by (4.35) in \cite{peres}, choosing $\delta<1/2$ we have 
\[
d_L\le (2\delta)^{[\frac{t_L}{t_\delta}]}.
\]
When $\eta\in\Ecal^0$ (\ref{tLprop}) applies and, in view of (\ref{calAdef}), (\ref{bootstrapacal}) and the observation following (\ref{bunibound}), the assertion of Lemma \ref{unibound} reduces to showing that
\begin{equation}
\P_{\tilde{\mu}^0}\big[{\cal A}\big]\le \frac{C}{L^6\log^2 L}.
\label{Pmu0calA}
\end{equation}

To estimate $\P_{\tilde{\mu}^0}\big[{\cal A}\big]$ we partition $\tilde{\Ecal}^0$ as
$
\tilde{\Ecal}^0=\Ecal^{0,1}\cup\Ecal^{0,2}\cup\Ecal^{0,3},
$
where
\[
\Ecal^{0,1}=\{\eta\in\tilde{\Ecal}^0: \eta_0\ge N-\rho_cL-{\alpha_L}+\frac{\beta_L}{2}, \ B(\eta)\le \frac{\beta_L}{2}\}\subset\Ecal^0,
\]
\[
\Ecal^{0,2}=\{\eta\in\tilde{\Ecal}^0: \ B(\eta)\le {\chi_L}\}\setminus \Ecal^{0,1} \qquad\text{and}\qquad \Ecal^{0,3}=\{\eta\in\tilde{\Ecal}^0: \ B(\eta)> {\chi_L}\},
\]
for some $\chi_L$ to be determined later. We may now write 
\begin{align}
\P_{\tilde{\mu}^0}\big[{\cal A}\big]\le \sup_{\eta\in\Ecal^{0,1}}\P_\eta\big[{\cal A}\big]+\sup_{\eta\in\Ecal^{0,2}}\P_\eta\big[{\cal A}\big]\tilde{\mu}^0\big[\Ecal^{0,2}\big]+\tilde{\mu}^0\big[\Ecal^{0,3}\big]
\label{bootstrapacal2}
\end{align}
and use (\ref{bootstrapacal}) to estimate the middle term of the sum on the right hand side. This gives
\begin{equation}
\big(1-\tilde{\mu}^0\big[\Ecal^{0,2}\big]\big)\P_{\tilde{\mu}^0}\big[{\cal A}\big]\le \sup_{\eta\in\Ecal^{0,1}}\P_\eta\big[{\cal A}\big]+ \big(\sup_{\eta\in\Ecal^{0,2}}\P_\eta\big[\tilde{T}\le t_L\big]+d_L\big)\tilde{\mu}^0\big[\Ecal^{0,2}\big]+\tilde{\mu}^0\big[\Ecal^{0,3}\big].
\label{tmu0pre}
\end{equation}
By Proposition \ref{secmax} we have
\[
\tilde{\mu}^0\big[\Ecal^{0,2}\big]\le CL\beta_L^{1-b}\quad\text{and}\quad\tilde{\mu}^0\big[\Ecal^{0,3}\big]\le CL\chi_L^{1-b}
\]
and by Lemma \ref{extexit} 
\[
\sup_{\eta\in\Ecal^{0,2}}\P_\eta\big[\tilde{T}\le t_L\big] \le C\frac{Lt_L^2\chi_L^{b-1}}{\phi_L^{1+b}}.
\]
When $\eta\in\Ecal^{0,1}$ for the process to leave $\Ecal^0$ at least $\beta_L/2$ particles need to be moved, either away from the condensate or onto the site that already contains $\beta_L/2$ particles. If fewer jumps occur by time $\frac{1}{L}$, the process will necessarily remain in $\Ecal^0$ up to time $\frac{1}{L}$ and the event ${\cal A}$ will not be realised. Since the total rate at which jumps occur is bounded by $mL$, the probability that at least $\beta_L/2$ jumps occur by time $\frac{1}{L}$ is dominated by $\P[X>\beta_L/2]$, where $X$ has Poisson distribution with mean $m$. Hence,
\[
\sup_{\eta\in\Ecal^{0,1}}\P_\eta\big[{\cal A}\big]\le \sum_{k\ge\beta_L/2}e^{-m}\frac{m^k}{k!}\le \frac{m^{\beta_L/2}}{(\beta_L/2)!},
\]
which decays faster than any power of $L^{-1}$. If we put together the preceding estimates, (\ref{tmu0pre}) gives
\[
\P_{\tilde{\mu}^0}\big[{\cal A}\big]\le CL\left(\frac{Lt_L^2\beta_L^{1-b}}{\phi_L^{1+b}}\chi_L^{b-1}+\chi_L^{1-b}\right).
\]
The estimate (\ref{Pmu0calA}) now follows, if we optimise the preceding bound by choosing
\[
\chi_L=\left(\frac{\phi_L^{1+b}\beta_L^{b-1}}{Lt_L^2}\right)^{\frac{1}{2(b-1)}}=\left(L^7\log^2 L\right)^{\frac{1}{b-1}}
\]
and this concludes the proof of Lemma \ref{unibound}.\CQFD\\

\subsection{Proof of Proposition \ref{trace} -- Replacement by the trace process\label{nullDelta}}

We will make use of the extended well $\tilde{\Ecal}^0$ and the exit time $\tilde{T}$ of the process $\eta(\cdot)$ from $\tilde{\Ecal}^0$, as well as the process restricted to $\tilde{\Ecal}^0$, denoted by $\xi(\cdot)$, and its $\delta$-mixing time $t_\delta$. All are defined in the preceding subsection. We begin the proof of Proposition \ref{trace} with the following estimate. 
\begin{lemma} 
\label{notimeout}
For any $u>0$ we have
\be
\sup_{\eta\in\Ecal^0}\E_\eta\Big[\int_{t_\delta}^{t_\delta+u}\mathbbm{1}_\Delta\big(\eta(s)\big)\ ds;\ \tilde{T}>t_\delta\Big]\le u\Big(\delta+\frac{\mu\big[\Delta\big]}{1-\mu\big[\Delta\big]}\Big).
\ee
\end{lemma}
{\bf Proof.} For $\eta\in X_{L,N}$ let us define
$\displaystyle
W(\eta)=\E_\eta\Big[\int_0^u\mathbbm{1}_\Delta\big(\eta(s)\big)\ ds\Big].
$
For any $\eta\in\Ecal^0$ using the Markov property we have
\begin{align*}
\E_\eta\Big[\int_{{t}_\delta}^{{t}_\delta+u}\mathbbm{1}_\Delta\big(\eta(s)\big)\ ds;\ \tilde{T}>{t}_\delta\Big]= \E_\eta\big[W\big(\eta({t}_\delta)\big);\ \tilde{T}>{t}_\delta\big]\le\tilde{\E}_\eta\big[W\big(\xi({t}_\delta)\big)\big]=\int W d\pi_\eta\notag,
\end{align*}
where $\pi_\eta$ is the distribution of the random variable $\xi({t}_\delta)$ under $\P_\eta$. The inequality above follows by coupling $\xi(\cdot)$ to $\eta(\cdot)$ up to time $\tilde{T}$, just as in the argument following (\ref{bootshalf}). Then,
\begin{align*}
\int W d\pi_\eta&\le \delta\  \sup_{\xi\in\tilde{\Ecal}^0}W(\xi)+\int W\ d\tilde{\mu}^0=\delta\  \sup_{\xi\in\tilde{\Ecal}^0}W(\xi)+\frac{1}{\mu\big[\tilde{\Ecal}^0\big]}\int_{\tilde{\Ecal}^0} W\ d\mu\notag\\
&=\delta\  \sup_{\xi\in\tilde{\Ecal}^0}W(\xi)+\frac{1}{\mu\big[\bigcup_{x\in\Lambda}\tilde{\Ecal}^x\big]}\int_{\bigcup_{x\in\Lambda}\tilde{\Ecal}^x}W\ d\mu\notag\\
&\le \delta u+\frac{1}{1-\mu\big[\Delta\big]}\int_{X_{L,N}} W\ d\mu=u\Big(\delta+\frac{\mu\big[\Delta\big]}{1-\mu\big[\Delta\big]}\Big).
\end{align*}
The first line follows from the definition of $t_\delta$ and (\ref{rest_meas}), the second line follows from translation invariance of the dynamics and the set $\Delta$ and the final step uses the invariance of $\mu$ under the dynamics of $\eta(\cdot)$. 
\CQFD\\

It is now straightforward to prove the replacement by the trace process.\\

\noindent
{\bf Proof of Proposition \ref{trace}.}
In view of Lemma \ref{notimeout} we have
\begin{align*}
\E_\eta\Big[\int_0^{t\theta_L}\hspace{-2mm}&\mathbbm{1}_\Delta\big(\eta(s)\big)ds\Big]
=\E_\eta\Big[\int_0^{t\theta_L}\hspace{-2mm}\mathbbm{1}_\Delta\big(\eta(s)\big)\ ds; \ \tilde{T}\le t_L\Big]+\E_\eta\Big[\int_0^{t\theta_L}\hspace{-2mm}\mathbbm{1}_\Delta\big(\eta(s)\big)\ ds;\ \tilde{T}>t_L\Big]\notag\\
&\le  t\theta_L\P_\eta\big[\tilde{T}\le t_L\big]+{t}_\delta+\E_\eta\Big[\int_{{t}_\delta}^{t\theta_L}\hspace{-1mm}\mathbbm{1}_\Delta\big(\eta(s)\big)\ ds;\ \tilde{T}>{t}_\delta\Big]\\
&\le t_\delta+t\theta_L\Big(\P_\eta\big[\tilde{T}\le t_L\big]+\delta+\frac{\mu\big[\Delta\big]}{1-\mu\big[\Delta\big]}\Big).
\end{align*}
Dividing by $\theta_L$, the assertion now follows from (\ref{tLprop}) and Corollary \ref{shaft}, which provides a vanishing bound on $\mu [\Delta ]$.
\CQFD

\section{Regularization and inter wells dynamics\label{sec:regul}}

In order to get matching upper and lower bounds on the transition rates of the auxiliary process (\ref{auxgen}), we have to replace it with a regularized version on a renormalized lattice, making use of Lipschitz continuity of the test functions $f$. In fact, we will show Proposition \ref{inter} for all $f\in C^2 (\T ,\R )$, which can be used to uniformly approximate Lipschitz functions.\\

Sections \ref{sec:regul} and \ref{sec:cap} are independent of the rest of the article, the proofs here rely on some invariant measure estimates provided in Section 
\ref{invmeas}. The only restrictions arising from the results of these sections on the parameter $b$ are given by equations \eqref{epsul} and \eqref{ellbarell}, resulting on a lower bound $b>5$. Clearly this is much better than the condition $b>20$ following from the uniform bounds on the exit rates (Lemma \ref{unibound}), the estimates required to prove equilibration in the wells 
(Proposition \ref{equilibration}) and tightness of the condensate dynamics (Proposition \ref{lprocess}).\\

 For future use where a different scaling may be needed, we choose to keep $\alpha_L$ and $\beta_L$ and other quantities derived from them  explicit, instead of replacing them by the values obtained with the choices in \eqref{alphabeta}.

\subsection{Rate estimates from capacity bounds \label{rate_estim}}

Generalizing the rates $r^\Lambda (z)$ given in (\ref{auxrate}) to non-empty subsets $A_1 ,A_2 \subset\Lambda$ with $A_1 \cap A_2 =\emptyset$, we write
\be
r^\Lambda (A_1 ,A_2 )=\frac{1}{|A_1 |}\sum_{x\in A_1\atop y\in A_2} r^\Lambda (y-x)\ ,
\ee
using that under $\mu$, conditioned on $\eta\in \Ecal^{A_1} =\bigcup_{x\in A_1} \Ecal^x$, the location of the condensate is uniformly distributed in $A_1$. In this notation we can identify $r^\Lambda (z)=r^\Lambda \big(\{ 0\} ,\{ z\}\big)$. Following \cite{beltranetal09b}, Lemma 6.8, we have
\begin{eqnarray} \label{rate-capacity}
\lefteqn{\mu^\Ecal \big[\Ecal^{ A_1}\big]\ r^{\Lambda}( A_1 , A_2 )=}\nonumber\\
& &\qquad\frac12 \bigg(\ca^\Ecal \big(\Ecal^{ A_1}, \Ecal \setminus \Ecal^{ A_1} \big)+ 
\ca^\Ecal \big(\Ecal^{ A_2}, \Ecal \setminus \Ecal^{ A_2} \big)-\ca^\Ecal \big(\Ecal^{ A_1 \cup  A_2},  \Ecal \setminus \Ecal^{ A_1 \cup  A_2} \big)\bigg)\ .
\end{eqnarray}
Here $\ca^\Ecal \big(\Ecal^{A}, \Ecal^{B} \big)$ denotes the capacity of the trace process between the sets of wells $\Ecal^{A}$ and $\Ecal^{B}$ for any $A,B\subset\Lambda$ with $A\cap B=\emptyset$. By Lemma 6.9 in \cite{beltranetal09b}, the trace process capacities satisfy  
\[
\mu\big[\Ecal\big]\,\ca^{\Ecal}\big(\Ecal^x, \Ecal^y\big)=\ca\big(\Ecal^x, \Ecal^y\big),
\]
where the latter are the full zero-range process capacities.
By Corollary \ref{shaft} in Section 
\ref{invmeas} $\mu[\Delta]\le CL\big( \alpha_L^{1-b}+\beta_L^{1-b}\big) $, 
hence we may replace $\mu^\Ecal[\cdot]$ and $\ca^\Ecal(\cdot, \cdot)$ throughout by the full zero-range process invariant measure and capacities $\mu[\cdot]$ and $\ca(\cdot, \cdot)$, at the cost of an $o(1/L)$ multiplicative error.\\

Propositions \ref{lower} and \ref{upper} provide the following lower and upper bounds on capacities between complementary sets. For each $A\subset\Lambda$ we have
\bea\label{comple}
\theta_L\ca (\Ecal^A ,\Ecal\setminus\Ecal^A )&\geq& H(1-\epsl ) \sum_{x\in A \atop y\not\in A} \ca_\Lambda (x,y)\ ,\nonumber\\
\theta_L\ca (\Ecal^A ,\Ecal\setminus\Ecal^A )&\leq& H(1+\epsu ) \sum_{x\in A \atop y\not\in A} \ca_\Lambda (x,y)\ ,
\eea
with a constant $H=H(b,\rho )=\frac{1}{(\rho -\rho_c)^{1+b} z_c\ I_b}$ and error terms
\be\label{epsul}
\epsl =C\big({L}{\beta_L^{1-b}}+{L}^{-1}+Le^{-\frac{\alpha_L}{2\beta_L\vee \sqrt{L}}}\big)\quad\mbox{and}\quad  \epsu\gg\frac{\alpha_L}{L},\,\,
\epsu^{3}\ge \frac{|A| L^{b+3}}{\alpha_L^{2b-2}(L-|A|)}.
\ee
With our choice of $\alpha_L$ and $\beta_L$ in (\ref{alphabeta}) $\epsu\gg\epsl$, and using (\ref{rate-capacity}) this leads to
\bea\label{boundss}
\frac{\theta_L \mu\big[\Ecal^{ A_1}\big]}{H} r^\Lambda ( A_1 , A_2 )&\leq&\sum_{x\in A_1 \atop y\in A_2} \ca_\Lambda (x,y)+\epsu \bigg(\sum_{x\in A_1 \atop y\not\in  A_1} \ca_\Lambda (x,y) +\sum_{x\in A_2 \atop y\not\in A_2} \ca_\Lambda (x,y)\bigg)\nonumber\\
\frac{\theta_L \mu\big[\Ecal^{ A_1}\big]}{H} r^\Lambda ( A_1 , A_2 )&\geq &\sum_{x\in A_1 \atop y\in A_2} \ca_\Lambda (x,y)-\epsu \bigg(\sum_{x\in A_1 \atop y\not\in  A_1} \ca_\Lambda (x,y) +\sum_{x\in A_2 \atop y\not\in A_2} \ca_\Lambda (x,y)\bigg)\ .
\eea
The bounds are given in terms of capacities of a simple symmetric random walk on the lattice $\Lambda$, which are
\begin{align}\label{laca}
\ca_\Lambda (x,y)=&\ca_\Lambda (0,y-x) =1/d_\Lambda (0,y-x) \nonumber\\
=&\frac{1}{|y-x|(L-|y-x|)} =\frac{1}{L}\bigg(\frac{1}{|y-x|}+\frac{1}{L-|y-x|}\bigg)\ .
\end{align}
using the standard embedding $\{ 0,\ldots ,L-1\}$ of the torus $\Lambda$ in $\Z$. 
In the following we will use a double embedding $\{ 0,\ldots ,L-1\}$ and $\{ -L, \ldots ,-1\}$ shifted by $-L$, where we identify sites $0$ and $-L$. The combined use leads to a more intuitive notation of the regularization procedure which is formulated symmetrically around $0$, and should not cause any confusion. Note also that we write $|x-y|\in \{ 0,\ldots ,L-1\}$ when $x,y$ are chosen consistently in the same embedding, and that this is not the distance $d_\Lambda$ on the discrete torus $\Lambda$. 
Some particular bounds we will make use of later are
\begin{align}\label{rabo}
\theta_L r^\Lambda (z)\leq \ & HL \Big(\frac{1}{|z|(L-|z|)}+2\epsu\sum_{x=1}^{L-1} \frac{1}{xL}\Big) \leq C\Big(\frac{1}{|z|} +\frac{1}{L-|z|} +\epsu\log L\Big)\nonumber\\
\mbox{and}\quad &\theta_L \sum_{z\neq 0} r^\Lambda (z)\leq C\big(1+\epsu \big)\log L
\end{align}
using that $\mu [\Ecal^z ]=(1-\mu [\Delta])/L =1/L (1+o(1))$. Note that the second line in \eqref{rabo} does not follow from the first; instead, we exploit the fact that $\sum_{z\neq 0} r^\Lambda(z)=r^\Lambda(0, \Lambda\setminus 0)$ and apply Proposition \ref{upper} directly.

\subsection{Proof of Proposition \ref{inter} \label{Prop_3.4}}

To establish (\ref{cap}) in Proposition \ref{inter} we will show that
\be\label{tosho}
\sup_{x\in\Lambda}\Big|\theta_L \Lcal^\Lambda f(x/L)-\Lcal^\T f(x/L)\Big|\to 0\quad\mbox{as }L\to\infty\ .
\ee
To achieve this we have to regularize the auxiliary generator $\Lcal^\Lambda$ (\ref{auxgen}), since using bounds of type (\ref{rabo}) directly will lead to diverging error terms. 
Fix an intermediate scale $\ell$ depending on $L$, which is a sequence of integer numbers such that
\be\label{lscale}
\ell \to\infty\quad\mbox{and}\quad\ell /L\to 0\quad\mbox{as }L\to\infty\ .
\ee
We partition the lattice $\Lambda$ into subsets (or boxes) of size $2\ell +1$, with only the box $\bar V_0$ centred in the origin $0$ being of larger size. For this we fix a second sequence $\bar\ell$ such that
\be\label{llscale}
\bar\ell /\ell\to \infty\ ,\quad \bar\ell /L\to 0\quad\mbox{as }L\to\infty\ ,
\ee
and take $\bar V_0 =\{ -\bar\ell ,\ldots ,\bar\ell\}$. We choose $\bar\ell$ such that the remaining lattice $\Lambda\setminus \bar V_0$ can be partitioned into boxes of size $2\ell +1$, i.e. there exists a sequence of  integers $\bar M$ such that
\be
L-(2\bar\ell +1)=\bar M (2\ell +1)\ .
\ee
The integer sequence $\bar M$ is bounded above by the real-valued sequence $M:=L/(2\ell +1)$, which characterizes the asymptotic number of boxes, since $\bar M /M\to 1$ as $L\to\infty$. The partition of the lattice is then given by $\bar V_0$ and the boxes
\be
V_m := \big\{ y\in\Lambda :\big| y-\bar\ell -(2\ell +1)m\big| \leq\ell\big\}\ ,\quad m=1,\ldots ,\bar M\ .
\ee
A choice of scales consistent with (\ref{lscale}) and (\ref{llscale}) which is sufficient for the proof is 
\begin{equation}
\label{ellbarell}
\ell =O(\alpha_L \log^3 L)\quad\mbox{and}\quad\bar\ell =O(\alpha_L \log^4 L )\quad\mbox{as } L\to\infty\ .
\end{equation}
This is possible with our choice of $\alpha_L$ in (\ref{alphabeta}), compatible with \eqref{epsul}, and implies in particular that
\be\label{scalecond}
\epsu M\log L\to 0\quad\mbox{and}\quad \frac{\bar\ell}{L} \log L\to 0
\quad\mbox{as }L\to\infty\ .
\ee
Here we have used that we can choose the error in (\ref{epsul}) as $\epsu =\frac{\alpha_L}{L} \log L$, since the size of the sets of wells we consider is $|A|\leq\bar\ell$. 
In fact, a quick computation shows that among the sequences $\alpha_L=L^{\frac{1}{2}+\gamma}$, 
$\gamma>0$, and $\ell$ satisfying \eqref{lscale}, \eqref{llscale}, the choices $\alpha_L=L^{\frac{1}{2}+\frac{5}{2b}}$ and $\ell=\alpha_L \log^3 L$ provide the smallest possible 
$\epsu$ (up to powers of $\log L$) which is consistent with \eqref{epsu}, and require a lower bound on the parameter $b$ of only $b>5$.\\

To regularize the rates $r^\Lambda$ (\ref{auxrate}), we rewrite the generator of the auxiliary process (\ref{auxgen}) as
\bea\label{zorro}
\Lcal^\Lambda f(x/L)&=&\sum_{z\in \Lambda} r^\Lambda (z) \Big( f\Big(\frac{x+z}{L}\Big) -f\Big(\frac{x}{L}\Big)\Big) \nonumber\\
&=&\frac{1}{2\ell +1} \sum_{y\in V_0} \sum_{z\in \Lambda} r^\Lambda (z-y) \Big( f\Big(\frac{x+z-y}{L}\Big) -f\Big(\frac{x}{L}\Big)\Big)\nonumber\\
&=&\frac{1}{2\ell +1} \sum_{y\in V_0} \sum_{z\in \bar V_0} r^\Lambda (z-y) \Big( f\Big(\frac{x+z-y}{L}\Big) -f\Big(\frac{x}{L}\Big)\Big) \nonumber\\
& &+\sum_{m=1}^{\bar M} \frac{1}{2\ell +1} \sum_{y\in V_0} \sum_{z\in V_m}  r^\Lambda (z-y) \Big( f\Big(\frac{x+z-y}{L}\Big) -f\Big(\frac{x}{L}\Big)\Big)\nonumber\\
&=:& \Lcal_1^\Lambda f(x/L) +\Lcal_2^\Lambda f(x/L)\ .
\eea
In the first step we used translation invariance of the rates $r^\Lambda$ and averaged them over the box $V_0 =\{ -\ell ,\ldots ,\ell\}$ of size $2\ell +1$ around the origin. 
For the first term $\Lcal_1^\Lambda$, we use the Lipschitz property of $f$ with constant $C_f$ to get
\bea\label{l1}
\theta_L \Lcal_1^\Lambda f(x/L) &\leq &C_f\frac{1}{2\ell +1}\sum_{z\in \bar V_0} \sum_{y\in V_0} \frac{|z-y|}{L} \theta_L r^\Lambda (z-y)\leq C_f\frac{2\bar\ell}{L} \theta_L \sum_{ 0\neq z \in \bar V_0} r^\Lambda (z)\nonumber\\
&\leq &C \frac{\bar\ell}{L}(1+\epsu)\log L\ ,
\eea
where we used (\ref{rabo}) in the last line. 
So with (\ref{scalecond}), $|\Lcal_1^\Lambda f(x/L)|\to 0$ as $L\to\infty$ uniformly for $x\in\Lambda$.\\

The second term of (\ref{zorro}) can be split as $\Lcal_2^\Lambda f(x/L)=\Lcal_{2A}^\Lambda f(x/L)+\Lcal_{2B}^\Lambda f(x/L)$ with
\bea\label{batman}
\Lcal_{2A}^\Lambda f(x/L) &:=&\sum_{m=1}^{\bar M} \frac{1}{2\ell +1} \sum_{y\in V_0} \sum_{z\in V_m}  r^\Lambda (z-y) \Big( f\Big(\frac{x+z-y}{L}\Big) -f\Big(\frac{x+\bar\ell +(2\ell +1)m}{L} \Big)\Big)\ ,\nonumber\\
\Lcal_{2B}^\Lambda f(x/L) &:=&\sum_{m=1}^{\bar M} r^\Lambda (V_0 ,V_m) \Big( f\Big(\frac{x+\bar\ell +(2\ell +1)m}{L} \Big) -f\Big(\frac{x}{L}\Big)\Big)\ .
\eea
For the first term we can use the Lipschitz continuity of $f$ with constant $C_f$ to get
\be
\big|\Lcal_{2A}^\Lambda f(x/L)\big| \leq C_f \,\frac{2\ell +1}{L}\sum_{m=1}^{\bar M} r^\Lambda (V_0 ,V_m)\ .
\ee
Using (\ref{boundss}) for complementary sets $ A_1 =V_0$ and $ A_2 =\Lambda\setminus V_0$ we get 
\bea
\theta_L \sum_{m=1}^{\bar M} r^\Lambda (V_0 ,V_m) &\leq &\theta_L r^\Lambda (V_0 ,\Lambda\setminus V_0)\leq K(1+\epsu )\,\frac{M}{L} \sum_{x\in V_0\atop y\not\in V_0} \Big( \frac{1}{|x{-}y|} +\frac{1}{L-|x{-}y|}\Big)  \nonumber\\
&\leq & CK(1+\epsu ) \frac{M(2\ell {+}1)\log L}{L}\leq C\log L\ ,
\eea
which leads to
\be
\theta_L \big|\Lcal_{2A}^\Lambda f(x/L)\big| \leq C_f (2\ell +1)\frac{\log L}{L}\to 0\quad\mbox{as }L\to\infty\mbox{ uniformly in }x\in\Lambda\ .
\ee

In the main contribution $\Lcal_{2B}^\Lambda$ in (\ref{batman}) we can use (\ref{boundss}) to obtain for the renormalized rates
\be
\theta_L r^\Lambda (V_0 ,V_m) = \frac{M}{z_c I_b}\bigg(\sum_{y\in V_0} \sum_{z\in V_m} \ca_\Lambda (0,z-y)+\epsu \mathcal{R}_L (V_0 ,V_m)\bigg)\ .
\ee
Analogously to the above, we can bound the remainder term as
\bea
\mathcal{R}_L (V_0 ,V_m)&:=&\sum_{x\in V_0 \atop y\not\in V_0 } \ca_\Lambda (0,y-x) +\sum_{x\in V_m\atop y\not\in V_m} \ca_\Lambda (0,y-x)\nonumber\\
&\leq & 
\frac{C}{L}\sum_{x\in V_0\atop y\not\in V_0}\Big( \frac{1}{|x-y|} +\frac{1}{L-|x-y|}\Big) \leq C(2\ell +1)\frac{\log L}{L}\ .
\eea
This leads to
\bea\label{sclim}
\lefteqn{\bigg|\theta_L\Lcal_{2B}^\Lambda f(x/L) -\frac{M}{z_c I_b}\sum_{m=1}^{\bar M} \sum_{y\in V_0\atop z\in V_m} \ca_\Lambda (0,z-y) \Big( f\Big(\frac{x+\bar\ell +(2\ell +1)m}{L} \Big) -f\Big(\frac{x}{L}\Big)\Big) \bigg|}\nonumber\\
& &\qquad\leq C\epsu \bar M \frac{M}{z_c I_b} \frac{(2\ell +1)\log L}{L}=C\epsu \bar M \log L\to 0
\eea
with (\ref{scalecond}), since $\bar M\leq M$. To conclude the proof it remains to identify the second term in the first line with $\Lcal^\T f(u)$ in the limit $L\to\infty$ when $x/L\to u\in\T$. To this end, we use the representation $\ca_\Lambda (0,z-y)=\frac{1}{L}\big(\frac{1}{|z-y|} +\frac{1}{L-|z-y|}\big)$ and note that
\be\label{toni1}
\sum_{y\in V_0\atop z\in V_m} \frac{1}{|z-y|} =\sum_{y,z\in V_0} \frac{1}{m(2\ell {+}1)+\bar\ell +z-y} =\frac{2\ell +1}{m+\bar\ell /(2\ell {+}1)}\Big( 1+O\Big(\frac{1}{m+\bar\ell /(2\ell {+}1)}\Big)\Big)\ .
\ee
Analogously, we have
\be\label{toni2}
\sum_{y\in V_0\atop z\in V_m} \frac{1}{L-|z-y|} =\frac{2\ell +1}{M-m-\bar\ell /(2\ell {+}1)}\Big( 1+O\Big(\frac{1}{M-m-\bar\ell /(2\ell {+}1)}\Big)\Big)\ .
\ee
So the contributions to (\ref{sclim}) of the correction terms in (\ref{toni1}) and (\ref{toni2}) vanish, since we have
\be
\frac{M}{L}\sum_{m=1}^{\bar M} \frac{2\ell +1}{(m+\bar\ell /(2\ell {+}1))^2}\leq C\frac{2\ell +1}{\bar\ell}\to 0\ ,
\ee
and an analogous bound holds for the second correction term. Therefore we can replace the second term in the first line of (\ref{sclim}) by
\be \label{arghh}
\frac{M}{z_c I_b}\frac{2\ell {+}1}{L}\sum_{m=1}^{\bar M} \Big(\frac{1}{m{+}\bar\ell /(2\ell {+}1)} +\frac{1}{M-m-\bar\ell /(2\ell {+}1)} \Big) \Big( f\Big(\frac{x{+}\bar\ell {+}(2\ell {+}1)m}{L} \Big) -f\Big(\frac{x}{L}\Big)\Big)
\ee
Rewriting the rates in this expression as
\[
\frac{1}{m{+}\bar\ell /(2\ell {+}1)} +\frac{1}{M-m-\bar\ell /(2\ell {+}1)}  = \frac{1}{M}\,\frac{1}{(m/M +\bar\ell /L)\,\,(1-m/M -\bar\ell /L)}
\]
and using that $\bar\ell /L\to 0$ and $\bar M/M\to 1$, (\ref{arghh}) converges to
\be
\frac{1}{z_c I_b}\int_0^1 \frac{1}{v(1-v)} \big( f(u+v)-f(u)\big)\, dv =\Lcal^\T f(u)
\ee
as $L\to\infty$ and $x/L\to u\in\T$. By regularity of $f$, $u\mapsto\Lcal^\T f(u)$ is a uniformly continuous function on $\T$, and with (\ref{sclim}) this implies (\ref{tosho}) 
and finishes the proof of Proposition \ref{inter}. \CQFD\\

\section{Capacity estimates\label{sec:cap}}

Upper and lower bounds for the capacities appearing in the proof of Proposition \ref{inter} in the previous section are obtained following closely Sections 4 and 5 in \cite{beltranetal09}. The extension of these methods in \cite{rebecca} lead to matching upper and lower bounds for complementary sets of wells also in the limit $L\to\infty$ with diverging particle density $N/L\to\infty$. Using precise estimates on the stationary measure summarized in Section \ref{invmeas}, we are able to improve estimates for the upper bound error that allow us to further extend these results to the the case $N/L\to\rho >\rho_c$.\\

To simplify notation, let $\tN$ be the typical number of particles at the condensate,
\begin{equation}
\label{tilde_N}
\tN:=N-\rho_cL\, .
\end{equation}
Recall that $\Lambda$ denotes the lattice $\Lambda=\Z/L\Z$, and we will use the subindex $k$ in $\Lambda_k=\Z/k\Z$ whenever we refer to a different lattice.

\subsection{Lower bound}

\begin{proposition}\label{lower}
Let $ A$ be a non-empty subset of the lattice $\Lambda$ and denote by $\Ecal^A$ the corresponding set of wells. Then, there exists a positive constant $C$ such that
\begin{equation}
\label{lower cap}
\tilde{N}^{b+1} \ca (\Ecal^A ,\Ecal\setminus \Ecal^A )\geq  \frac{1}{z_c\ I_b} \sum_{x\in A \atop \,y \notin A} \ca_\Lambda (x,y) \Big(1-C\big({L}{\beta_L^{1-b}}+{L}^{-1}+Le^{-\frac{\alpha_L}{2\beta_L\vee\sqrt{L}}}\big)\Big)
\end{equation}
provided that $L$ is large enough. Here $I_b=\int_0^1 u^b (1-u)^b\,du$, and for any $x,\,y \in \Lambda$, $\ca_\Lambda(x,y)$ denotes the capacity (\ref{laca}) of the simple symmetric random
walk on $\Lambda$. 
\end{proposition}

\noindent\textbf{Proof.} 
The variational formulation of the capacity establishes that $\ca(\Ecal^A, \Ecal \setminus \Ecal^A )=\inf_{F\in {\cal H}(\Ecal^A)} \cal{D}(F)$, where 
\[
{\cal H}(\Ecal^A)=\big\{F: F(\eta)=1, \eta \in \Ecal^A,\, F(\eta)=0, \eta \in \Ecal \setminus \Ecal^A \big\}, 
\]
and ${\cal D}$ is the Dirichlet form
\begin{align}
\label{Dirichlet}
{\cal D}(F)=\frac 1 2 \sum_{\eta\in X_{L,N}} \sum_{x\in \Lambda}\ \sum_{r=-1,\, 1}\mu(\eta)\,g(\eta_x)\,\frac 1 2 \,\big[F(\eta^{x,x+r})-F(\eta)\big]^2 .
\end{align}

Define the tube
\be\label{tubedef}
T^{x,y}_{L,N} =\Big\{\eta\in X_{L,N} :\eta_x +\eta_y \geq \tN-\frac{\alpha_L}{2} ;\,\eta_z \leq \beta_L ,\, z\neq x,y\Big\}
\ee
connecting the two wells at $x$ and $y$. For any function $F: X_{L,N}\to \R$ in ${\cal H}(\Ecal^A)$, we can bound the Dirichlet form by
\be
\label{1estimate}
{\cal D} (F)
\geq \sum_{x\in A ,\,y\notin A}\left(\frac12  \sum_{\substack{ z\in \Lambda \\ r=-1,\,1}}  \sum_{\eta\in T^{x,y}_{L,N}} \mu (\eta ) \, \frac{1}{2} \, g(\eta_z )\big( F(\eta^{z,z+r})-F(\eta )\big)^2 \right) \ .
\ee
This holds because with $\beta_L \ll\alpha_L$ the tubes only overlap within wells on which the function $F$ is constant, equal to $0$ or $1$: 
$T_{L,N}^{x,y}\cap T_{L,N}^{x',y'} \subseteq \Ecal^x \cup \Ecal^y$ and necessarily $x=x'$ or $y=y'$ for $N,L$ large enough.\\

Fix $x\in A,\,y\notin A$. Let $\delta_z$ denote the configuration with one and only one particle at $z$, and for $\eta\in X_{L,N}$ let $\xi =\eta -\delta_z$, where summation of configurations is performed componentwise.  Given $M, K \in \N$ and  $\zeta \in X_{M,K}$ we will denote by 
$g!(\zeta)=\prod_{u \in \Lambda_M} g!(\zeta_u)=\prod_{u \in \Lambda_M} \zeta_u^b$ with the choice of rates  in \eqref{rates}. Define  the set of configurations $J^{x,y}_{L,N}=\{\xi\in X_{L,N-1}:\  |{\xi}_x +{\xi}_y-\tilde{N}| \leq \alpha_L/2,\  \xi_u \leq \beta_L -1, u\neq x,\,y\}$ and note that $\xi\in J^{x,y}_{L,N} \Rightarrow \xi+\delta_z\in T_{L,N}^{x,y}$ for all $z\in\Lambda$.
Then we can rewrite the parenthesis in the right-hand side of \eqref{1estimate} as
\begin{align}
{\cal D}_{x,y}(F)&= \frac12  \sum_{\substack{z\in\Lambda\\r=-1,\,1}} \sum_{\xi \in X_{L,N-1}\atop \xi +\delta_z \in T^{x,y}_{L,N}} 
\frac{1}{Z_{L,N}} \frac{1}{g!(\xi)}\,\frac{1}{2}\,\big( F(\xi+\delta_{z+r})-F(\xi+\delta_z)\big)^2 \nonumber\\
&\ge\frac{1}{2Z_{L,N}} \sum_{\xi \in J^{x,y}_{L,N}} \frac{1}{g!(\xi)}\, \sum_{\substack{z\in\Lambda\\r=-1,1}}\,\frac{1}{2}\big( F(\xi+\delta_{z+r})-F(\xi+\delta_z)\big)^2,\label {middle}
\end{align} 
with
\begin{equation}
\label{ZLN}
Z_{L,N}=\sum_{\eta \in X_{L,N}} \frac{1}{g!(\eta)}=z_c^L \,\nu[S_L=N].
\end{equation}
If  $F(\xi +\delta_x) \neq F(\xi +\delta_y)$, then in full analogy with \cite{beltranetal09} we may define
\[
f=f_\xi: \Lambda\to \R,\quad f(z)=\frac{F(\xi+\delta_z)-F(\xi+\delta_y)}{F(\xi+\delta_x)-F(\xi+\delta_y)}.
\]
Note that 
\[ 
f(x)=1,\, f(y)=0, \text{ and}\quad f(z+r)-f(z)= \frac{F(\xi+\delta_{z+r})-F(\xi+\delta_z)}{F(\xi+\delta_x)-F(\xi+\delta_y)},
\]
so we can estimate the last sum in (\ref{middle}) from below by $2L\ca_\Lambda(x,y)$ to get
\begin{equation}
{\cal D}_{x,y} (F)\ge\,\frac{L}{Z_{L,N}}\,\ca_\Lambda (x,y) \sum_{\xi \in J^{x,y}_{L,N}} \frac{1}{g!(\xi)}\,\big( F(\xi +\delta_x) -F(\xi +\delta_y)\big)^2 \ .\label{last}
\end{equation}

Clearly, this bound holds trivially also for terms where $F(\xi +\delta_x )=F(\xi +\delta_y )$. Let $\zeta$ be a configuration in $\Lambda \setminus\{ x,y\}$. 
For such a $\zeta\in X_{L-2,k}$ we define the function $G_\zeta:\{0,1,\ldots,N-k\}\to\R$ by $G_\zeta(i)=F(\bar{\zeta})$ where $\bar{\zeta}\in X_{L,N}$ coincides with $\zeta$ outside $\{x,y\}$, and $\bar{\zeta}_x=i,\ \bar{\zeta}_y=N-k-i$. Let us also define the set ${\mathcal O}_L=\big\{k\in\N:\ |k-\rho_cL|\le \frac{\alpha_L}{2}\big\}.$ With this notation \eqref{last} can be written as 
\begin{align}
\label{on-xy-tube}
{\cal D}_{x,y} (F)\ge \frac{L\, \ca_\Lambda (x,y)}{Z_{L,N}}  \sum_{k\in {\mathcal O}_L} \,\,\sum_{\zeta\in X_{L-2,k} \atop \zeta_z\le\beta_L-1} \frac{1}{g!(\zeta)} \sum_{i=1}^{N-k-2} \frac{\big(  G_\zeta (i+1)-G_\zeta (i)\big)^2}{i^b\,(N-k-i-1)^b} \,.
\end{align}
Note that if $|k-\rho_cL|\le\frac{\alpha_L}{2}$ and $\zeta \in X_{L-2,k}$ with $\max_z\zeta_z\le \beta_L$ then $G_\zeta(i) \equiv 0$ when $i\le \beta_L$ for $L$ sufficiently large, because $\bar{\zeta}\in{\cal E}^y.$ Likewise, $G_\zeta(i) \equiv 1$ when $i\ge N-k-\beta_L$ for sufficiently large $L$, because $\bar{\zeta}\in{\cal E}^x$. Hence, the rightmost sum in (\ref{on-xy-tube}) reduces to 
\begin{align*}
\sum_{i=\beta_L}^{N-k-\beta_L-1} \frac{\big(  G_\zeta (i+1)-G_\zeta (i)\big)^2}{i^b\,(N-k-i-1)^b}\ .
\end{align*}
Minimizing the preceding expression under the constraint 
\[
\displaystyle  \sum_{i=\beta_L}^{N-k-\beta_L-1} {G_\zeta (i+1)-G_\zeta (i)}=1\,, 
\] 
the estimate in (\ref{on-xy-tube}) becomes
\begin{align}
\label{sum_on_i} 
{\cal D}_{x,y} (F)\ge \frac{L\, \ca_\Lambda (x,y)}{Z_{L,N}}\hspace{-1mm}  \sum_{k\in {\mathcal O}_L} \hspace{-1mm}\sum_{\zeta\in X_{L-2,k} \atop \zeta_z\le\beta_L-1} \frac{1}{g!(\zeta)}\left( \sum_{i=\beta_L}^{N-k-1-\beta_L} \hspace{-1mm}{i^b\,(N\!-\!k-\!1-\!i)^b} \right)^{-1}\hspace{-2mm}.
\end{align}
Using the fact that for any $C^2$ function $f:(0,1)\to\R$, $m\in\N$ and sequence $k_m \ll m$, the following bound holds
\[
\left|\frac{1}{m}\sum_{i=k_m}^{m-k_m} \hspace{-1mm}f\Big(\frac{i}{m}\Big)-\int_{\frac{k_m-1/2}{m}}^{1-\frac{k_m-1/2}{m}} f(x)\ dx\right|\le\frac{1}{24m^2}\sup_{x\in(0,1)}|f''(x)| ,
\]
we get that
\[
\sum_{i=\beta_L}^{N-k-1-\beta_L} \hspace{-1mm}{i^b\,(N\!-\!k-\!1-\!i)^b}\le (N-k-1)^{2b+1}I_b\ \Big(1+\frac{C}{(N-k-1)^2}\Big)
\]
where $I_b=\int_0^1x^b(1-x)^b dx$. 
Hence, the estimate in (\ref{sum_on_i}) further becomes (provided $N,L$ are sufficiently large)
\begin{align}
{\cal D}_{x,y} (F)&\ge \frac{L\, \ca_\Lambda (x,y)}{Z_{L,N}\ I_b}\hspace{-1mm}  \sum_{k\in {\mathcal O}_L} \hspace{-1mm}\sum_{\zeta\in X_{L-2,k} \atop \zeta_z\le\beta_L-1} \frac{1}{g!(\zeta)\,(N-k-1)^{2b+1}} \times\Big(1-\frac{C}{(N-k-1)^2}\Big)\nonumber\\
&\ge \frac{z_c^{L-2}L\, \ca_\Lambda (x,y)}{Z_{L,N}\ I_b\ }\hspace{-1mm}  \sum_{k\in {\mathcal O}_L}\frac{1}{(N-k)^{2b+1}} \nu\big[S_{L-2}=k,\ M_{L-2}\le\beta_L-1\big] \nonumber\\
&\ge \frac{z_c^{L-2}L\, \ca_\Lambda (x,y)}{Z_{L,N}\ I_b\ }\, \, \nu\Big({(N-S_{L-2})^{-(2b+1)}}; K_L\Big),\label{midDxy}
\end{align}
where $K_L$ is the event $K_L=\{ \ |S_{L-2}-\rho_cL|\le \frac{\alpha_L}{2}, M_{L-2}\le\beta_L-1\}$, $M_{L-2}=\max_{\eta \in \Lambda_{L-2}} \eta_x$. Just as in (\ref{taylor}), if we
define $k_*=\nu\big(S_{L-2};\ K_L\big)$ then
\[
\nu\Big({(N-S_{L-2})^{-(2b+1)}}; K_L\Big)\ge\frac{1}{\tilde{N}^{2b+1}}\nu\big[K_L\big]+\frac{2b+1}{\tilde{N}^{2b+2}}(k_*-\rho_cL).
\]
In a fashion similar to (\ref{tailK}) and ({\ref{kstar}) we have
\[
\nu\big[K_L^c\big]= \nu\big[M_{L-2}> \beta_L-1\big]+\nu\big[|S_{L-2}-\rho_cL|>\alpha_L, M_{L-2}\le \beta_L\big]\le C\ \big(L\beta_L^{1-b}+e^{-\frac{\alpha_L}{\beta_L\vee\sqrt{L}}}\big).
\]
and
\[
0\le \rho_c(L-2)-k_*\le C\ L^2\beta_L^{1-b}.
\]
In view of the preceding estimates $(\ref{midDxy})$ becomes
\begin{align*}
\tilde{N}^{b+1}{\cal D}_{x,y}(F)&\ge  \frac{z_c^{L-2}L\tilde{N}^{-b}\, \ca_\Lambda (x,y)}{Z_{L,N}\ I_b\ } \Big(1-C\big({L}{\beta_L^{1-b}}+{L}^{-1}+Le^{-\frac{\alpha_L}{2\beta_L\vee\sqrt{L}}}\big)\Big),
\end{align*}
and the statement of the proposition follows from Proposition \ref{NagaevRates}, 
since $F$ was an arbitrary function $F: X_{L,N}\to \R$, with $\, F|_{\Ecal^A}=1,\, F|_{\Ecal \setminus \Ecal^A}=0$.
\CQFD

\subsection{Upper bound}

\begin{proposition}\label{upper}
Let $A \subset\Lambda$ be a non-empty subset of the lattice, and denote by $\Ecal^A$ the corresponding set of wells. Then for any $\epsu$ satisfying
\begin{equation}
\label{epsu}
\epsu \gg \frac{\alpha_L}{L}\quad\mbox{and}\quad\epsu^3 \geq\frac{|A| L^4 \tN^{b-1}}{ \alpha_L^{2b-2}  (L-|A|)}
\end{equation}
we have
\begin{equation}
\label{hawkeye}
\tN^{b+1} \ca (\Ecal^A ,\Ecal \setminus \Ecal^A )\leq (1+\epsu )\frac{1}{I_b z_c} \sum_{x\in A ,\,y\notin A} \ca_\Lambda (x,y)\ .
\end{equation} 
In \eqref{hawkeye}  $I_b=\int_0^1 u^b (1-u^b)\,du$, and for any $x,\,y \in \Lambda$, $\ca_{\Lambda}(x,y)$ denotes the capacity (\ref{laca}) of the simple symmetric random walk on $\Lambda$. 
\end{proposition}

\paragraph {Construction of the test function.}

For  $1\le y <  L$, let $f_{0,y}:\Lambda \to [0,1]$ be the function that realises the capacity $\ca_\Lambda(0,y)$ between the sites $0$ and $y$ for the symmetric, nearest neighbour random walk in $\Lambda$. By elementary results of potential theory, it is known that $f_{0,y}(z)$ equals the probability that the random walk reaches $0$ before
visiting $y$, when started at $z$. The formula for $f_{0,y}$ can be easily derived,
\be\label{f0y}
 f_{0,y}(z) = \left\{ \begin{array}{lll}
         1-\frac z y\quad & \mbox{ if $\;0\le z<y$};\\
         \\
         \frac{z-y}{L- y}\,& \mbox{ if $\;y\le z <L$}.\end{array} \right.
\ee

Given $\epsilon>0$, consider the smooth function $H_{\epsilon}:[0,1]\to[0,1]$ given by
\begin{align}
\label{hepsilon}
H_{\epsilon}(t):=\frac1 {I_b} \int_0^{\phi_{\epsilon}(t)} u^{b}(1-u)^b \, du\,,
\end{align}
where $I_b=\int_0^1 u^b (1-u)^b\,du$, and $\phi_{\epsilon}:\R\to [0,1]$ is a smooth, non-decreasing function such that $\phi_{\epsilon}(t)+\phi_{\epsilon}(1-t)=1\, \forall t\in \R$, 
$\phi_{\epsilon}|_{(-\infty,3\epsilon] }\equiv 0$, $\phi_{\epsilon}|_{[1-3\epsilon,\, \infty)}\equiv 1$, $\sup_{u\in \R}|\phi_\epsilon '(u)|\le 1+C\epsilon$ and $\sup_{u\in [0,1]}|\phi_\epsilon (u)-u|\le C' \epsilon$, for some universal constants $C, C' >0$. It can be easily checked that 
\be
\label{Hfunction}
H_{\epsilon}(t)+H_{\epsilon}(1-t)=1\ \forall t\in \R,\quad
H_{\epsilon}|_{(-\infty,3\epsilon]}\equiv 0,\quad H_{\epsilon}|_{[1-3\epsilon, \infty)}\equiv 1\,.
\ee
Enumerate the sites in $\Lambda$ as $0=x_1,\, x_2,\dots,\, x_L=y$, so that $f_{0,y}(x_j)> f_{0,y}(x_{j+1})$ for all $j$.\\

Given $\eta \in X_{L,N}$, let $\tilde \eta \in \R^{\Lambda }$ be given by  
\begin{equation}
\tilde\eta_y= \eta_z -\rho_c, \ y\in \Lambda.
\end{equation}
For $y\neq 0$, define $F_{0,y}^j:X_{L,N}\to [0,1]$, $1\le j\le L-1$ as 
\bea
\label{F0yj}
&& F_{0,y}^1(\eta)=H_{\epsilon}\left(\frac{\tilde \eta_0}{\tN}\right), \notag\\
&& F_{0,y}^j(\eta)=H_{\epsilon}\left(\frac{\tilde \eta_0}{\tN} +\max\left\{ -\epsilon ;\min\Big\{ \frac 1 \tN \sum_{i=2}^{j} \tilde \eta_{x_i}; \epsilon \Big\} \right\}\right).
\eea
Finally, for $y\neq 0$, let $F_{0,y}:X_{L,N}\to \R$ be the convex combination
\begin{align}
\label{F0y}
F_{0,y}(\eta)=\sum_{j=1}^{L-1} \big[f_{0,y}(x_j)-f_{0,y}(x_{j+1})\big]\ F_{0,y}^j(\eta)\,.
\end{align}
To define $F_{0,0}$, consider a $C^1$, non-decreasing function $h:\R\to [0,1]$ such that $h|_{(-\infty, 2\epsilon]}\equiv 0$, $h|_{(1-2\epsilon, +\infty]} \equiv 1$, $\sup_x h'(x)\le 2$, and let 
\begin{equation}
\label{thor}
F_{0,0}(\eta):=h\big(\teta_0/\tN \big).
\end{equation}
In order to construct the candidate test functions for the capacities in \eqref{rate-capacity} we need to stitch the functions $\{F_{0,y}\}_{y\in \Lambda}$ together.\\
 
Let us hence define the sets 
$S=\{u\in \R^{\Lambda},\, \sum_{i\in \Lambda} u_i=1\}$, and, for $0\neq y \in \Lambda$,
\begin{align}
\label{sets2}
&{\cal T}_{\epsilon}^{0,y}=\left\{ u\in S:\ u_0+u_y\ge 1-\epsilon \right\},\notag  \\
&{\cal A}_{\epsilon}=\left\{ u\in S: \ u_0\ge 1-\epsilon\right\}\notag\\
\intertext{and the disjoint sets}
&{\cal K}^{0,y}_{\epsilon}=\left({\cal T}_{\epsilon}^{0,y} \setminus {\cal A}_{2 \epsilon}\right)\cap \left\{u\in S,\,u_z <\frac{\epsilon}{2},\, z\neq 0,\,y\right\}.
\end{align}
Let now 
\begin{equation}
\label{unity}
\big\{\Theta^{0,y}=\Theta^{0,y}(L,\epsilon)\big\}_{y \in \Lambda}
\end{equation}
be a smooth partition of unity of $S$, such that  $\Theta^{0,y}:S\to [0,1],\, y \in \Lambda$, $\sum_{y} {\Theta}^{0,y}\equiv 1$, 
\begin{align}
\label{partition}
{\Theta}^{0,y} \big|_{{\cal K}_{\epsilon}^{0,y}} \equiv 1
\quad
\mbox{and}\quad \sup_{u, u' \in S}|\Theta^{0,y}(u)-\Theta^{0,y}(u')|\le \frac{C}{\epsilon}\,\|u-u'\|_2,
\end{align}
$C$ a positive constant independent of $\epsilon$, $L$ and $y\in \Lambda$. To construct this partition of unity, consider a bump function $g:\R \to [0,1],\,g\equiv 1$ on
$(-\infty, \frac{1}{2} \epsilon],\,g\equiv 0$ on $[\epsilon, +\infty)$, $\|g'\|_{\infty} \le \frac{C}{\epsilon}$. Let $d(u, {\cal K}_\epsilon^{0,y})=\inf\{\|u-w\|_2, w \in {\cal K}^{0,y}_\epsilon\}: \R^L \to \R$ be the Euclidean distance from a point $u$ to the set ${\cal K}^{0,y}_\epsilon$,
and let $\Theta^{0,y}(u)=g(d(u, {\cal K}^{0,y}_\epsilon)),\, 1\le y\le L-1$. Note that these have disjoint supports and take values in $[0,1]$.
Define $\Theta^{0,0}(u)=1-\sum_{1\le y\le L-1} \Theta^{0,y}(u)$. \\

The candidate to solve the variational problem for the capacity 
$\ca \big(\Ecal^0, \Ecal \setminus \Ecal^0  \big)$ 
is then $F_0:X_{L,N}\to \R$,
\begin{align}
\label{F0}
F_0(\eta)=\sum_{y\in \Lambda} \Theta^{0,y}\big(\tilde \eta/\tN \big)\, F_{0,y}(\eta).
\end{align}

\begin{lemma}
\label{F0-proper}
The function $F_0$ in \eqref{F0} satisfies
\begin{align}
& F_0(\eta)=F_{0,y}(\eta),\quad\mbox{for } \tilde\eta\in {\cal K}_{\epsilon}^{0,y},\label {fitem}\\ 
& F_0(\eta)=1, \quad\mbox{for } \tilde \eta_0\ge (1-2\epsilon)\tilde N, \label{sitem}\\
& F_0(\eta)=0, \quad\mbox{for } \tilde\eta_0 \le 2\epsilon \tN, \label{titem}\\
&|F_0(\eta)-F_0(\eta')| \le C\frac{1}{\epsilon \tN}\,\|\eta-\eta'\|_2,\quad\mbox{for }\eta,\,\eta'\in X_{L,N}. \label{foitem}
 \end{align}
\end{lemma}
\noindent\textbf{Proof.}
The first assertion follows from \eqref{partition}. For the second one, note that $\phi _{\epsilon}\equiv 1$ on $[1-3\epsilon, 1]$ hence 
$F_{0,y}^j(\eta)=F_{0,y}(\eta)=1$ for all $y\in \Lambda\setminus \{0\}$, $j=1\dots L-1$, and also $F_{0,0}(\eta)=1$, if $\tilde \eta_0\ge (1-2\epsilon)\tN$. Similarly, \eqref{titem} follows from $F_{0,y}^j(\eta)=F_{0,y}(\eta)=0$ for all $y\in \Lambda\setminus \{0\}$, $j=1\dots L-1$, $F_{0,0}(\eta)=0$, if $\tilde \eta_0\le 2\epsilon \tN$. Finally, \eqref{foitem} is a consequence of \eqref{partition}, if we use that the supports of the functions $\Theta^{0,y},\,y\neq 0$ are disjoint, and hence at most four terms in the difference 
\[F_0(\eta)-F_0(\eta')=\sum_{y\in \Lambda} \big[\Theta^{0,y}\big(\tilde \eta/\tN \big)\, F_{0,y}(\eta) - \Theta^{0,y}\big(\tilde \eta'/\tN \big)\, F_{0,y}(\eta')\big]\]
do not vanish identically. \CQFD\\

By \eqref{sitem} and \eqref{titem}, $F_0\big|_{\Ecal^0}=1$, $F_0\big|_{\Ecal \setminus \Ecal^0}=0$, hence we can use it to estimate the capacity between $\Ecal^0$ and $\Ecal \setminus \Ecal^0$.

\paragraph{Single well capacities.} 

Given $f:X_{L,N} \to \R$ and $U\subseteq X_{L,N}$, define
\begin{align}
\label{CNA}
C_N(f;U):=\frac 1 2 \sum_{\eta\in U} \sum_{x\in \Lambda}\ \sum_{r=-1,\, +1}\mu(\eta)\,g(\eta_x)\,\frac 1 2 \,\big[f(\eta^{x,x+r})-f(\eta)\big]^2,
\end{align}
 and for $x,\,y \in \Lambda_{L}$ let 
 \begin{align}
 \label{IN0y}
 {\cal I}_N^{x,y}:=\{\eta\in X_{L,N},\, \tilde \eta_x +\tilde \eta_y \ge \tilde N -\alpha_L \},
 \end{align}
  where $\sqrt{L}\ll \alpha_L \ll L$ is the constant used
 in the definition of the wells \eqref{valley}.\\

 We first notice that the only relevant contributions to the Dirichlet form \eqref{Dirichlet} originate from configurations in the sets ${\cal I}^{0,y}_N,\,y\in \Lambda$.

\begin{lemma}
\label{neglig1}
\begin{align*}
C_N\big(F_0;X_{L,N}\setminus \cup_{0\neq y\in \Lambda}{\cal I}_N^{0,y}\big)\le  C\,\frac{L^2}{\epsilon^2\,\tN^2\,\alpha_L^{2b-2}\,},
\end{align*}
$C$ a constant independent of $L$ and $\epsilon$.
\end{lemma}
\noindent\textbf{Proof.}
By \eqref{sitem}, \eqref{titem}, \eqref{foitem} in Lemma \ref{F0-proper},
\begin{align*}
C_N\big(F_0;X_{L,N}\setminus \cup_{0\neq y\in \Lambda}{\cal I}_N^{0,y}\big)
&\le  C \,\frac{L}{\epsilon^2 \tN^2}\,\mu\Big[\Big\{\eta\in X_{L,N}\setminus \cup_{0\neq y\in \Lambda}{\cal I}_N^{0,y};\, 
2\epsilon\le \frac{\tilde \eta_0}{\tilde N}\le 1-2\epsilon \Big\}\Big]\\
& \le  C\, \frac{L}{\epsilon^2 \tN^2}\frac{L}{\alpha_L^{2b-2}}=C\,\frac{L^2}{\epsilon^2\,\tN^2\,\alpha_L^{2b-2}\,}\,.
\end{align*}
The last line follows from Proposition \ref{uno} in Section \ref{invmeas}.\CQFD\\

Note that 
\begin{equation}
\label{sumcap}
C_N\big(F_0; \cup_{0\neq y\in \Lambda}{\cal I}_N^{0,y}\big)\le \sum_{y\in \Lambda} C_N\big(F_0; {\cal I}_N^{0,y}\big),
\end{equation}
as all terms in the latter sum are non-negative.

\begin{lemma}
\label{neglig3}
Let $\epsilon>0$ be such that $\epsilon \tN\gg \alpha_L$.
There exist positive constants $C, C'$ which are independent of $L$, and such that for any $0\neq y \in \Lambda\,$,
\begin{align}
\label{last_step}
\big|C_N\big(F_0;{\cal I}_N^{0,y}\big)-C_N\big(F_{0,y};{\cal I}_N^{0,y}\big)\big| \le C e^{-C'\alpha_L/\sqrt{L}}.
\end{align}
\end{lemma}
\noindent\textbf{Proof.}  By Lemma \ref{F0-proper} we may just consider those configurations $\eta$ such that $2\epsilon \tilde N\le \tilde \eta_0\le (1-2\epsilon) \tilde N$. Note that the fact that 
$\eta\in {\cal I}_N^{0,y}$ implies that $\tilde\eta_y\ge \epsilon \tilde N$, and also that for 
$F_0(\eta)\neq F_{0,y}(\eta)$ it is necessary that $\Theta^{0,y}(\tilde\eta/\tilde N)\neq 1$, which combined with the previous conditions implies the existence of $z\in \Lambda,\,z\neq0,y,$ such that $\tilde \eta_z\ge \epsilon \tilde N/2$. From \eqref{CNA} and \eqref{foitem} in Lemma \ref{F0-proper}, we get
\begin{align}
\label{last_step_1}
&\big|C_N\big(F_0;{\cal I}_N^{0,y}\big)-C_N\big(F_{0,y};{\cal I}_N^{0,y}\big)\big| \notag\\
&\hspace{.5cm}\le C\,\frac{L}{\epsilon^2 \tN^2}\,\mu\Big[\Big\{\eta\in {\cal I}_N^{0,y},\,2 \epsilon \le \frac{\tilde \eta_0}{\tN} \le (1-2\epsilon) ;\,
 \frac{\teta_y}{\tN}\ge \epsilon \, ;\,\exists z\neq 0,y,\,\tilde \eta_z\ge \frac{\epsilon \tN}{2}
\Big\}\Big]\notag\\
&\hspace{.5cm}\le C\,\frac{L}{\epsilon^2 \tN^2} \mu\Big[\Big\{\exists z\neq 0,\,y, \sum_{i\ge2, x_i \neq z}^{L-2} \teta_{x_i} \le -\frac{\epsilon \tN}{3} \Big\}\Big]\,.
\end{align}
By Proposition \ref{largedev}, 
\begin{align*}
\big|C_N\big(F_0;{\cal I}_N^{0,y}\big)-C_N\big(F_{0,y};{\cal I}_N^{0,y}\big)\big| \le C\,\frac{L^2}{\epsilon^2 \tN^2} \exp\{-C' \epsilon^2\tN^2/(L-3)\}.
\end{align*}
with $\frac{\epsilon^2\tN^2}{L-3}\ge \frac{\alpha_L^2}{L}$.
\CQFD

\begin{proposition}
\label{upper-bound}
Consider $\epsilon>0$ such that $\epsilon \tN\gg \alpha_L$. 
Let $y\neq 0$ in $\Lambda$. 
There exists a constant $C>0$ independent of $L$, $\epsilon$ and $y$, such that
\begin{align}
\label{1valley}
 C_N\big(F_{0,y},{\cal I}_N^{0,y}\big)\le \left(1+C \epsilon \right) \frac{L\,z_c^{L-2}}{\tilde N^{2b+1}\,Z_{L,N}}  \,\frac{\ca_\Lambda(0,y)}{I_b}\,,
\end{align}
$\ca_\Lambda(x,y)$ the capacity between sites $x$ and $y$ for the simple symmetric random walk in $\Lambda$.
\end{proposition}
\noindent\textbf{Proof.}
We first focus on terms $r=1$ in (\ref{CNA}). 
Fix $\tx \in \Lambda$, and let $1\le i,\,j\le L$ be the indexes of $\tx,\, \tx+1$ in the enumeration of  $\Lambda$ determined by sites $0,\,y$. Let us suppose $i<j$, that is, $f_{0,y}(\tx+1) < f_{0,y}(\tx)$. By \eqref{F0y}, we get
\begin{align}
\label{bracket}
F_{0,y}(\eta^{\tx,\tx+1})-F_{0,y}(\eta)=\sum_{k=i}^{j-1} \big[f_{0,y}(x_k)-f_{0,y}(x_{k+1})\big]\big[ F_{0,y}^k(\eta^{\tx,\tx+1})- F_{0,y}^k(\eta)\big].
\end{align}
Replacing \eqref{bracket} in  definition \eqref{CNA}, 
\begin{align}
\label{inner1}
&\sum_{\eta\in {\cal I}_N^{0,y}} \mu(\eta)\, g(\eta_{\tx})\, \frac{1}{2}\, \big[F_{0,y}(\eta^{\tx,\tx+1})-F_{0,y}(\eta)\big]^2  \notag \\
& \hspace{1cm} \le \big[ f_{0,y}(\tx)-f_{0,y}(\tx+1)\big]\sum_{k=i}^{j-1} \big[ f_{0,y}(x_k)-f_{0,y}(x_{k+1})\big] \notag \\ 
&\hspace{4.5cm} \times \sum_{\eta\in {\cal I}_N^{0,y}} 
\mu(\eta) \frac{g(\eta_{\tx})}{2} \big[ F_{0,y}^k (\eta^{\tx,\tx+1})-F_{0,y}^k(\eta)\big]^2
\end{align}
by Jensen's inequality. Note that terms associated to $\eta\in {\cal I}_N^{0,y}$ and $\tilde \eta_0< 2\epsilon \tilde N$ or $\tilde \eta_0>(1-2\epsilon)\tilde N $ vanish by the definition of 
$F^k_{0,y}$, \eqref{F0yj}.\\

Let 
\begin{equation}
\label{too_small}
{\cal B}_N^{0,y}:=\Big\{\eta\in {\cal I}_N^{0,y},\, 2\epsilon \tilde N \le\tilde \eta_0\le (1-2 \epsilon) \tilde N,\, \exists \,\, 2\le l\le L-1, \big|\sum_{i=2}^l \tilde \eta_{x_i}\big|\ge \epsilon \tilde N
 \Big\}.
\end{equation}
Fix now $i\le k\le j-1$ and consider the sum in the last line of \eqref{inner1} associated to the index $k$. We have
\begin{eqnarray*}
 \sum_{ {\cal I}_N^{0,y}} 
\mu(\eta) \frac{g(\eta_\tx)}{2} \big[ F_{0,y}^k (\eta^{\tx,\tx+1})-F_{0,y}^k(\eta)\big]^2
=& \displaystyle{\sum_{{\cal I}_N^{0,y}\setminus  {\cal B}_N^{0,y}}} \mu(\eta) \frac{g(\eta_{\tx})}{2} \big[ F_{0,y}^k (\eta^{\tx,\tx+1})-F_{0,y}^k(\eta)\big]^2\\
&+\displaystyle{\sum_{{\cal B}_N^{0,y}}} \mu(\eta) \frac{g(\eta_{\tx})}{2} \big[ F_{0,y}^k (\eta^{\tx,\tx+1})-F_{0,y}^k(\eta)\big]^2,
\end{eqnarray*}
where
\begin{align}
\label{residue0}
\displaystyle{\sum_{{\cal B}_N^{0,y}}} \mu(\eta) \frac{g(\eta_{\tx})}{2} \big[ F_{0,y}^k (\eta^{\tx,\tx+1})-F_{0,y}^k(\eta)\big]^2\le \frac{C}{\epsilon^2\tilde N^2}\, \mu\big[{\cal B}_N^{0,y}\big].
\end{align}
In ${\cal B}_N^{0,y}$ we have $\sum_{i=2}^{L-1}\tita_{x_i}\le \alpha_L$. Together with the fact that there is $2\le l\le L-1$ such that $|\sum_{i=2}^l \tita_{x_i} |\ge \epsilon \tN$, this  implies that either
{\it a)} $\sum_{i=2}^l \tita_{x_i}\le -\epsilon \tN$ or {\it b)} $\sum_{i=l+1}^{L-1} \eta_{x_i}\le \alpha_L-\epsilon \tN$. Once more, by Proposition \ref{largedev}, and with the hypothesis that $\epsilon \tN\gg \alpha$, we get
\begin{equation}
\label{BN1}
\mu\big[{\cal B}_N^{0,y}]\le C_1e^{-C_2 \alpha_L/\sqrt{L}},
\end{equation}
$C_1, C_2$ positive constants that do not depend on $L$. 
From \eqref{BN1} and \eqref{residue0} it follows that  $C_N(F_{0,y}, {\cal B}_N^{0,y})$ is negligible to any polynomial order in \eqref{1valley}.\\

On the other hand,
\begin{align}
\label{main_term}
&\sum_{{\cal I}_N^{0,y}\setminus  {\cal B}_N^{0,y}} \mu(\eta) \frac{g(\eta_{\tx})}{2} \big[ F_{0,y}^k (\eta^{\tx,\tx+1})-F_{0,y}^k(\eta)\big]^2 \notag\\
&\hspace{2cm} \le
\frac{1}{2 Z_{L,N}} 
\sum_{\xi \in {\cal J}_{N-1}^{0,y}} 
\frac{1}{g!(\xi)}
\left[H_{\epsilon}\Big(\frac{1}{\tilde N} + \sum_{r=0}^k \frac{\tilde \xi _{x_r}}{\tilde N} \Big)
-H_{\epsilon}\Big( \sum_{r=0}^k \frac{\tilde \xi _{x_r}}{\tilde N} \Big) \right]^2,
\end{align}
with $g!(\xi)=\prod_{0\le i\le L-1 } g!(\xi_i)=\prod_{0\le i\le L-1} \xi_i^b$, 
where ${\cal J}_{N-1}^{0,y}$ is the set
\begin{equation*}
{\cal J}_{N-1}^{0,y}:=\Big\{\xi \in X_{L, N-1};2\epsilon \le\frac{\tilde \xi_0}{\tN}\le (1-2 \epsilon), \, \tilde \xi_0+\tilde \xi_y\ge \tilde N -\alpha_L;
\big|\sum_{i=2}^l \tilde \xi_{x_i}\big|\le \epsilon \tilde N \,\forall \,2\le l\le L-1
 \Big\}.
\end{equation*}
By the definition \eqref{hepsilon} of $H_{\epsilon}$, for $\xi \in {\cal J}_{N-1}^{0,y}$, we have
\[
\left[ H_{\epsilon}\Big(\frac{1}{\tilde N} + \sum_{r=1}^k \frac{\tilde \xi _{x_r}}{\tilde N} \Big)
-H_{\epsilon}\Big( \sum_{r=1}^k \frac{\tilde \xi _{x_r}}{\tilde N} \Big) \right]^2=
\frac{1}{I_b^2}\left(\int_{\phi_\epsilon \big(\sum_{0\le r\le k} \frac{\tilde \xi_r}{\tilde N}\big)}^{\phi_\epsilon \big(\frac{1}{\tilde N}+\sum_{0\le r\le k} \frac{\tilde \xi_r}{\tilde N}\big)} u^b\,(1-u)^b\,du\right)^2,
\]
and hence
\begin{align}
\label{main_term2}
&\sum_{\xi \in {\cal J}_{N-1}^{0,y}} 
\frac{1}{g!(\xi)}
\left[H_{\epsilon}\Big(\frac{1}{\tilde N} + \sum_{r=0}^k \frac{\tilde \xi _{x_r}}{\tilde N} \Big)
-H_{\epsilon}\Big( \sum_{r=0}^k \frac{\tilde \xi _{x_r}}{\tilde N} \Big) \right]^2 \notag\\
&\hspace{.5cm}\le \sum_{m=-\rho_c(L-2)}^{\alpha_L-1} \,\,\sum_{\substack{\xi \in X_{L-2, M}, \\ M=\rho_c(L-2)+m \\ |\sum_{i=2}^l \tilde \xi_{x_i}\big|\le \epsilon \tilde N,  \,\forall\, l}}  \frac{1}{\tilde N^{2b}}\,  \frac{1}{I_b^2}\, \frac{1}{g!(\xi)} \notag\\
&\hspace*{2.5cm}\times \sum_{\iota=2\epsilon \tilde N+\rho_c }^{(1-2\epsilon) \tilde N+\rho_c} \left(\int_{\Phi_{\iota}}^{\Phi_{\iota+1}} u^b(1-u)^b du\right)^2 \frac{1}{\big(\frac{\iota}{\tilde N}\big)^b \big(1-\frac{m+\iota-2\rho_c +1}{\tilde N}\big)^b}
\, ,
\end{align}
where $\Phi_{\iota}:=\phi_\epsilon\Big(\frac{\iota-\rho_c}{\tilde N}+\frac{1}{\tilde N}\sum_{r=1}^k \tilde \xi_{x_r}\Big)$.
In order to derive the right hand side in \eqref{main_term2}, we consider the distribution of $m+\rho_c (L-2)$ particles in the
$L-2$ sites other than $0$ and $y$, and then the distribution of the remaining particles between these two sites, such that  $\iota$ particles are 
assigned to $0$ and $\tilde N-m+2\rho_c -\iota-1$ particles to $y$.\\

Terms in \eqref{main_term2}  associated to indices $-\rho_c (L-2)\le m\le -\epsilon \tN$ correspond to moderate or large deviations events $\sum_{2\le i\le L} \tilde{\xi}_{x_i}\le m$. 
Proposition \ref{largedev} implies that their combined contribution decays faster than any power of $L$, and is hence negligible.\\

Now, if $-\epsilon \tN \le m\le \alpha_L-1$, and recalling that $\big|\sum_{i=2}^l \tilde{\xi}_{x_i}\big|\le \epsilon \tN\, \forall l$, we have
\begin{align*}
\int_{\Phi_{\iota}}^{\Phi_{\iota+1}} 
\frac{u^b(1-u)^b}{\big(\frac{\iota}{\tilde N}\big)^b \big(1-\frac{m+\iota-2\rho_c +1}{\tilde N}\big)^b}\,\, du
&\le (\Phi_{\iota+1}-\Phi_\iota) 
\left(\frac{\Phi_{\iota+1}}{\frac{\iota}{\tilde N}}
\right)^b \left( \frac{1-\Phi_\iota}{1-\frac{m+\iota-2\rho_c +1}{\tilde N}} \right)^b  \\
&\le (\Phi_{\iota+1}-\Phi_\iota) \sup_{u\in [2\epsilon,1-2\epsilon]} \left(\frac{\phi_\epsilon(u)}{u-\frac{3}{2}\epsilon} \right)^{2b} \notag\\
&\le  \frac{1}{\tilde N} \,(1+C\epsilon)^{2b}.
\end{align*}
In order to obtain the last line, we applied Cauchy's mean value theorem to $\frac{\phi_\epsilon(u)}{u-\frac{3}{2}\epsilon}=\frac{\phi_\epsilon(u)-\phi_\epsilon(\frac{3}{2}\epsilon)}{u-\frac{3}{2}\epsilon}$, using that $\phi_\epsilon(3\epsilon/2)=0$ and $\sup_u |\phi'_epsilon (u)|\leq 1+C\epsilon$ from the definition after (\ref{hepsilon}). 
Replacing this bound in \eqref{main_term2} yields
\begin{align*}
&\sum_{\xi \in {\cal J}_{N-1}^{0,y}} 
\frac{1}{g!(\xi)}
\left[H_{\epsilon}\Big(\frac{1}{\tilde N} + \sum_{r=0}^k \frac{\tilde \xi _{x_r}}{\tilde N} \Big)
-H_{\epsilon}\Big( \sum_{r=0}^k \frac{\tilde \xi _{x_r}}{\tilde N} \Big) \right]^2 \notag\\
&\hspace{.1cm}\le (1+C\epsilon) \sum_{m=-\rho_c(L-2)}^{\alpha_L-1} \,\,\sum_{\substack{\xi \in X_{L-2, M}, \\ M=\rho_c(L-2)+m \\ |\sum_{i=2}^l \tilde \xi_{x_i}\big|\le \epsilon \tilde N,  \,\forall\, l}}   \frac{1}{\tilde N^{2b+1}}\,  \frac{1}{I_b^2}\, \frac{1}{g!(\xi)}\sum_{\iota=2\epsilon \tilde N+\rho_c }^{(1-2\epsilon) \tilde N+\rho_c} \int_{\Phi_{\iota}}^{\Phi_{\iota+1}} u^b(1-u)^b \,du 
\notag\\
&\hspace{.1cm} \le \left(1+C \epsilon\right) \frac{1}{I_b\,\tilde N^{2b+1}}  \sum_{m=-\rho_c(L-2)}^{\alpha_L-1} \,\,\sum_{\substack{\xi \in X_{L-2, M}, \\ M=\rho_c(L-2)+m \\ |\sum_{i=2}^l\tilde \xi_{x_i}\big|\le \epsilon \tilde N,  \,\forall\, l}}  \frac{1}{g!(\xi)}
\,\le\,\,\left(1+C\epsilon \right)\ \frac{z_c^{L-2}}{I_b\,\tilde N^{2b+1}}.
\end{align*}
We can now put the pieces together: replace the above estimate in \eqref{main_term} to derive from \eqref{inner1} that
\begin{equation*}
\sum_{\eta\in {\cal I}_N^{0,y}} \mu(\eta)\, \frac{g(\eta_{\tx})}{2}\, \big[F_{0,y}(\eta^{\tx,\tx+1})-F_{0,y}(\eta)\big]^2  
 \le\,\left(1+C\epsilon\right)\frac{1}{2\,I_b\,} \,\frac{z_c^{L-2}}{Z_{L,N} \, \tilde N^{2b+1}} \\
 \big[ f_{0,y}(\tx)-f_{0,y}(\tx+1)\big]^2.
\end{equation*}
We conclude that
\begin{align*}
C_N(F_{0,y}, {\cal I}_N^{0,y})&\le \left(1+C\epsilon \right) \frac{z_c^{L-2}}{\tilde N^{2b+1}\,Z_{L,N}}  \, \frac{1}{I_b\,} \, \frac{1}{2}\sum_{\tilde x \in \Lambda} \,\sum_{r=-1,+1}\frac{1}{2}\big[ f_{0,y}(\tx)-f_{0,y}(\tx+r)\big]^2\\
&=\left(1+C \epsilon \right) \frac{L\,z_c^{L-2}}{\tilde N^{2b+1}\,Z_{L,N}}  \,\frac{\ca_\Lambda(0,y)}{I_b}.
\end{align*}
The extra factor $L$ in the line above is compensating for the uniform weight $\frac{1}{L}$ appearing in the computation of the capacities for the symmetric simple random walk on the torus.
\CQFD\\

Lemmas \ref{neglig1},  \ref{neglig3}, Propositions \ref{upper-bound}, \ref{NagaevRates} and observation \eqref{sumcap} combined yield an upper bound for the capacity
\begin{equation}
\label{blackwidow}
\tN^{b+1}\ca\big(\Ecal^0, \Ecal\setminus \Ecal^0 \big)\le \frac{1}{I_b z_c} \sum_{y\neq 0} \ca_{\Lambda}(0,y) \Big(1+C\epsilon+\frac{1}{L}+\frac{\tN^{b+1}}{\epsilon^2 \alpha_L^{2b-2} \log L}\Big).
\end{equation}
 By the symmetry of the model, a similar estimate holds for 
$\ca \big(\Ecal^x, \Ecal \setminus  \Ecal^x  \big)$, $x\in \Lambda$, evaluating in this case the Dirichlet form at $F_x=F_0(\tau_{-x}(\eta))$, if $\tau_{-x}: X_{L,N}\to X_{L,N}$ is the shift $[\tau_{-x}(\eta)]_z:=\eta_{z+x},\,z\in \Lambda$.

\paragraph{Multiple wells capacitites.}

As described in Section \ref{sec:regul}, to compute the rescaled rates of the process we need bounds on the capacities 
\begin{equation}
\label{multi_valley}
\ca \big( \Ecal^A,  \Ecal \setminus \Ecal^A \big),
\end{equation}
where we recall the notation $\Ecal^A=\cup_{z\in A}\Ecal^z$,
 $A \subset \Lambda$. The strategy is to compute these in terms of the single well capacities estimated in the previous paragraph.\\

At this point one would like  to take advantage of the subadditivity of the capacities to bound $\ca(\Ecal^A, \Ecal \setminus \Ecal^A)\le \sum_{x\in A, y\notin A} \ca(\Ecal^x, \Ecal^y)$. To make it work, instead of estimating the capacity between complementary sets $\ca(\Ecal^x, \Ecal \setminus \Ecal^x)$ in terms of $\sum_{y\neq x} \ca_\Lambda(x,y)$, as we did in the previous paragraph, it would be  necessary to estimate each term $\ca(\Ecal^x, \Ecal^y)$ directly from the associated capacity $\ca_{\Lambda}(x,y)$ for the simple symmetric random walk on $\Lambda$. Unfortunately we are not able to do this: we need the function $F_{x,y}$ to be the relevant term on the set ${\cal I}_N^{x,y}$, but it also contributes to the other sets ${\cal I}_N^{x,z},\,z\neq y$, in non-negligible ways. 
Indeed, it seems that the only solution to this problem is to consider a partition of unity that singles out each of the sets ${\cal K}_\epsilon^{xz}, z\neq x,$ as in \eqref{sets2}, and to define $F_x$ as in \eqref{F0}, thus reducing the Dirichlet form to the sum of the Dirichlet form restricted to these sets, while still providing a suitable test function $F_x$ such that $F_x\big|_{\Ecal^x}=1,\,F_x\big|_{\Ecal\setminus \Ecal^x}=0$; that is, the construction of the previous paragraph. Following \cite{beltranetal09}, we next show how these functions can be combined to estimate $\ca(\Ecal^A, \Ecal\setminus \Ecal^A)$.\\

\noindent{\bf Proof of Proposition \ref{upper}:}
Fix $\epsu>0$ be as in the statement of the proposition. 
Define $F^A:=\sum_{z \in A}F_z$, so that $F^A \big|_{\Ecal^A}=1$ and $F^{A}\big|_{\Ecal\setminus\Ecal^A}=0$ by \eqref{sitem} and \eqref{titem}.
As a first observation, notice that it is enough to consider  
\begin{align}
\label{cncombined0} 
C_N\Big(F^A;\,\cup_{z \in A}\big(\cup_{\substack{y\in \Lambda\\ y\neq z}} {\cal I}_N^{z,y}\,\big) \Big). 
\end{align}
Indeed, using Cauchy Schwarz,
\begin{align}
\label{neg1'}
C_N\Big( F^{A};\,X_{L,N}\setminus \cup_{z \in A}\big(\cup_{\substack{y\in \Lambda\\ y\neq z}} {\cal I}_N^{z,y}\,\big)\Big)&\le |A| \sum_{z\in A} 
C_N \Big(F^z;\,X_{L,N}\setminus \cup_{z \in A}\big(\cup_{\substack{y\in \Lambda\\ y\neq z}} {\cal I}_N^{z,y}\,\big)\Big)\notag \\
&\le C \frac{ |A|^2 L^2}{\epsu^2 \tN^2 \alpha_L^{2b-2} } \qquad \mbox{by Lemma \ref{neglig1},}
\end{align}
since $X_{L,N}\setminus \cup_{z \in A}\big(\cup_{\substack{y\in \Lambda\\ y\neq z}} {\cal I}_N^{z,y}\big)\,\, \subset \,\, X_{L,N} \setminus \cup_{y\neq u} \,{\cal I}_N^{u,y}$ for each given $ u\in A$.\\

Next, 
\begin{align}
\label{splitup}
&C_N\Big(F^A;\,\cup_{z \in A}\big(\cup_{\substack{y\in \Lambda\\ y\neq z}} {\cal I}_N^{z,y}\,\big) \Big)\le \sum_{z\in A} \sum_{y\neq z} 
C_N\Big( F^A;\, {\cal I}_N^{z,y}\Big) \notag\\
&\hspace{.3cm}=\sum_{z\in A} \sum_{\substack{y\in A\\ y\neq z}} C_N\Big(F_z+F_y +\sum_{\substack{w\neq z,\,y\\ w\in A}} F_w;\,{\cal I}_N^{z,y}\Big)
+\,\sum_{z\in A}\sum_{y\notin A} C_N\Big(F_z+\sum_{\substack{w\neq z \\ w\in A}} F_w;\,{\cal I}_N^{z,y}\Big).
\end{align}
We will consider the terms in these sums individually.\\

We have
\begin{equation*}
\Big| C_N\Big(F_z+F_y+\sum_{\substack{w\neq z,\,y\\ w\in A}} F_w; {\cal I}_N^{z,y}\Big)-C_N\Big(F_z+F_y;{\cal I}_N^{z,y}\Big)\Big| \le C \frac{L}{\epsu^2 \tN^2} 
\,\mu\Big[ {\cal I}_N^{z,y} \cap \big\{\,\tita_x\ge 2 \epsu \tN, x\neq y,z \big\} \Big],
\end{equation*}
on account of \eqref{titem} and \eqref{foitem}. Now on ${\cal I}_N^{z,y}\cap \{\tita_x\ge 2\epsu \tN\}$ we have
$
\sum_{u\neq z,y,x} \tita_u\le \tN-(\tN-\alpha_L)-\epsu \tN=\alpha_L-2\epsu\tN, 
$
and by Proposition \ref{largedev}
\begin{equation*}
\label{2dos}
\mu\Big[{\cal I}_N^{z,y} \cap \{\tita_x\ge 2 \epsu \tN,\,x\neq y, z \} \Big] \le C_1e^{-C_2 \alpha_L/\sqrt{L} },
\end{equation*}
$C_1, C_2$ positive constants that do not depend on $L$. Furthermore, the arguments in Lemma \ref{neglig3}
also yield
\begin{equation*}
\label{3tres}
\big| C_N\Big(F_z+F_y;\, {\cal I}_N^{z,y}\Big)-C_N\Big(F_{z,y}+F_{y,z};\, {\cal I}_N^{z,y}\Big)\big| \le C_1 e^{-C_2 \alpha_L/\sqrt{L}},
\end{equation*}
so that
\begin{equation}\label{final 1}
\big| C_N\Big(F_z+F_y+\sum_{w\neq z,\,y} F_w; {\cal I}_N^{z,y}\Big)-C_N\big(F_{z,y}+F_{y,z};{\cal I}_N^{z,y}\Big)\Big| \le C_1 e^{-C_2 \alpha_L/\sqrt{L}}, \ \ \  z,\,y \in A,
\end{equation}
where the values of the constants may change from line to line. By similar estimates,
\begin{align}\label{final 2}
\big| C_N\Big(F_z+\sum_{w\neq z,\,y} F_w;\, {\cal I}_N^{z,y}\Big)-C_N\Big(F_{z,y};\, {\cal I}_N^{z,y}\Big)\big| \le C_1 e^{-C_2 \alpha_L/\sqrt{L}}, \quad z \in A,y \notin A. 
\end{align}
We next show that terms $C_N\Big(F_{z,y}+F_{y,z},\,{\cal I}_N^{z,y}\Big)$,  $y, z \in A$, are negligible to first order.  Let $\eta \in {\cal I }_N^{z,y}\setminus\big({\cal B}_N^{z,y}\cup {\cal B}_N^{y,z}\big)$ (see definition \eqref{too_small}). We get
\begin{align*}
\label{term1}
F^j_{z,y}(\eta)+F_{y,z}^{L-j}(\eta)&=H_\epsu\left(\frac{\tilde \eta_z}{\tilde N}+\frac{1}{\tilde N}\sum_2^{j}\tilde \eta_{x_i}\right)+H_\epsu\left(\frac{\tilde \eta_y}{\tilde N}+\frac{1}{\tilde N}\sum_{j+1}^{L-1}\tilde \eta_{x_i}\right)\notag\\
&=H_\epsu\left(\frac{\tilde \eta_z}{\tilde N}+\frac{1}{\tilde N}\sum_2^{j}\tilde \eta_{x_i}\right)+H_\epsu\left(1-\left[ \frac{\tilde \eta_z}{\tilde N}+\frac{1}{\tilde N}\sum_2^{j}\tilde \eta_{x_i}\right]\right)=1.
\end{align*}
The first equality holds by the fact that the enumerations $\big\{x_i\big\}_{i=1}^L$ and $\big\{{\bar x}_j\big\}_{j=1}^{L}$ on $\Lambda$ determined by $f_{z,y}$ and
 $f_{y,z}$ satisfy
$x_i={\bar x}_{L-i+1}$. The second equality is due to the identity $H_\epsu(t)+H_\epsu(1-t)\equiv 1$. 
Then
\begin{align*}
&F_{z,y}(\eta)+F_{y,z}(\eta)=\\
&\qquad=\sum_{j=1}^{L-1} \big[f_{z,y}(x_j)-f_{z,y}(x_{j+1})\big]\ F_{z,y}^j(\eta)+\sum_{j=1}^{L-1} \big[f_{y,z}(\bar x_j)-f_{y,z}(\bar x_{j+1})\big]\ F_{y,z}^j(\eta) \\
&\qquad=\sum_{j=1}^{L-1} \big[f_{z,y}(x_j)-f_{z,y}(x_{j+1})\big]\ F_{z,y}^j(\eta)+\sum_{j=1}^{L-1}  \big[f_{z,y}(x_j)-f_{z,y}(x_{j+1})\big]\ F_{y,z}^{L-j}(\eta) \\
&\qquad=\sum_{j=1}^{L-1} \big[f_{z,y}(x_j)-f_{z,y}(x_{j+1})\big]\ \Big(F_{z,y}^j(\eta)+F_{y,z}^{L-j}(\eta)\Big) \\
&\qquad= \sum_{j=1}^{L-1} \big[f_{z,y}(x_j)-f_{z,y}(x_{j+1})\big]=1, \qquad \eta \in {\cal I }_N^{z,y}\setminus\big({\cal B}_N^{z,y}\cup {\cal B}_N^{y,z}\big),
\end{align*}
where  we used that $f_{y,z}(\bar x_{L-j})-f_{y,z}(\bar x_{L-j+1})=f_{z,y}(x_j)-f_{z,y}(x_{j+1})$. In particular
\begin{equation}
\label{term1bis}
C_N\Big( F_{z,y}+F_{y,z}\,,\,{\cal I }_N^{z,y}\setminus\big({\cal B}_N^{z,y}\cup {\cal B}_N^{y,z}\big)\Big)=0.
\end{equation}
Finally, the arguments leading to \eqref{BN1} yield
\begin{equation}
\label{neg2'}
C_N\Big( F_{z,y}+F_{y,z}\,,{\cal B}_N^{z,y}\cup {\cal B}_N^{y,z}\Big)\le C_1 e^{-C_2 \alpha_L/\sqrt{L}}. 
\end{equation}

Collecting all results from equations \eqref{neg1'} 
to \eqref{neg2'}, we obtain
\[
C_N(F^A;X_{L,N})\le \sum_{z\in A}\sum_{y\notin A} C_N(F_{z,y};\,{\cal I}_N^{z,y}) + C\, \frac{ |A|^2 L^2}{\epsu^2 \tN^2 \alpha_L^{2b-2}}\,,
\]
and by  Proposition \ref{upper-bound},
\begin{align*}
 C_N\Big(F^A, X_{L,N}\Big) \le
 \left(1+C \epsu \right) \frac{L\,z_c^{L-2}}{\tilde N^{2b+1}\,Z_{L,N}} \frac{1}{I_b} \sum_{z\in A} 
 \sum_{y\notin A} \ca_\Lambda(z,y) \,+
 C\, \frac{ |A|^2 L^2}{\epsu^2 \tN^2 \alpha_L^{2b-2}}\,.
\end{align*}
In particular, pulling out the first term as a common factor, and applying the rough estimates $\ca_\Lambda(z,y)\ge \frac{1}{L^2}$ and $\sum_{z\in A} 
 \sum_{y\notin A} \ca_\Lambda(z,y)\ge \frac{|A|(L-|A|)}{L^2}$, condition \eqref{epsu} yields
\begin{align*}
\tN^{b+1}\ca\big(\Ecal^A, \Ecal\setminus \Ecal^A\big)\,\le \,
 \left(1+C \epsu \right)\, \frac{L\, z_c^{L-2}}{\tilde N^{b}\,Z_{L,N}} \frac{1}{I_b} \sum_{z\in A,\,y\notin A} \ca_\Lambda(z,y)\,.
\end{align*}
The assertion of the proposition follows from Proposition \ref{NagaevRates}, as $\epsu_L \ge 1/L$. \CQFD\\

\section{Tightness and Limiting Distribution\label{tight}}

\subsection{Proof of Proposition \ref{lprocess} -- Tightness\label{tightness}}
In this section we will prove tightness of the distributions $\Q^L$ of  $\big(Y_t^L:\,t\ge0\big):=\big(\frac{1}{L}\psi\big(\eta^\Ecal({t\theta_L})\big):\,t\ge 0\big)$. By Aldous's tightness criterion (cf. Theorem 16.10 in \cite{Billingsley}), it suffices to show that for any $\eta\in\Ecal, \epsilon>0, t>0$, 
\begin{equation}
\label{aldous}
\lim_{\delta\downarrow 0}\limsup_{L\to\infty}\sup_{s\le\delta}\sup_{\tau\in{\mathfrak T}_t}\P_\eta\big[d_\T\big(Y^L_{\tau+s},Y^L_\tau\big)>\epsilon\big]=0,
\end{equation}
with $d_\T(x,y )=|x-y|\big(1-|x-y|\big)$ the distance in the torus $\T$, where, as before, all algebraic operations are performed modulo $Z$, i.e.  $x-y=(x-y)_{\,\mbox{mod}(\Z)}$.  In \eqref{aldous} ${\mathfrak T}_t$ is the set of stopping times (for the trace process) bounded by $ t$. \\

Note that by the strong Markov property, for any $\eta\in\Ecal$ we have
\bea
\P_\eta\big[d_\T\big(Y^L_{\tau+s},Y^L_\tau \big)>\epsilon\big]&=\E_\eta\big[ \P_{\eta^\Ecal(\tau)}\big[d_\T\big(Y_s^L,Y_0^L\big)>\epsilon\big]\big]\notag \\
&\le \sup_{\eta\in\Ecal} \P_{\eta}\big[d_\T\big(Y_s^L,Y_0^L\big)>\epsilon\big].
\eea
We will denote by 
\[
r^\Ecal(\eta ;z):=\sum_{x\in\Lambda}\1_{\Ecal^x} (\eta )\sum_{\xi\in \Ecal^{x+z}}r^\Ecal(\eta, \xi)
\] 
the total jump rate from a configuration $\eta$ in some well $\Ecal^x$ to a new well $\Ecal^{x+z}$. 
By It\^{o}'s formula, for any $s>0$ we have
\begin{align*}
&d_\T\big(Y_s^L,Y_0^L\big)=\\
&\quad=\int_0^{s\theta_L}
\sum_{z\neq 0}r^\Ecal\big(\eta^\Ecal (u);z\big)\,\Big[d_\T\big(\psi_L(\eta(u))+z/L, Y_0^L\big) - d_\T\big(\psi_L(\eta(u)), Y_0^L\big)\Big]\,du
+M_{s\theta_L},
\end{align*}
where $\{M_s\}_{s\ge 0}$ is a martingale with quadratic variation
\[
\langle M\rangle_s=\int_0^s\sum_{z\neq 0}
r^\Ecal\big(\eta(u);z\big)
\Big[d_\T\big(\psi_L(\eta^\Ecal (u))+z/L, Y_0^L\big) - d_\T\big(\psi_L(\eta(u)), Y_0^L\big)\Big]^2
\,du.
\]
By the triangle inequality we get
\[
\big|d_\T\big(\psi_L(\eta)+z/L, Y_0^L\big) - d_\T\big(\psi_L(\eta), Y_0^L\big)\big|
\le d_\T\big(z/L,0\big)=\frac{z}{L}\Big(1-\frac{z}{L}\Big).
\]
By Doob's inequality, for any $\eta\in\Ecal$,
\begin{align}\label{doobin}
\P_\eta\Big[\sup_{s\le \delta\theta_L}|M_s|>\frac{\epsilon}{2}\Big] & \le\frac{16}{\epsilon^2}\,\E_\eta\Big[\langle M\rangle_{\delta\theta_L}\Big] \le \frac{16}{\epsilon^2}\sum_{z\neq 0}\frac{z^2}{L^2}\Big(1-\frac{z}{L}\Big)^2\E_\eta\Big[\int_0^{\delta\theta_L}r^\Ecal\big(\eta^\Ecal (u);z\big)du\Big]\\
& \le \frac{16}{\epsilon^2}\sum_{z\neq 0}\frac{z}{L}\Big(1-\frac{z}{L}\Big)\,\E_\eta\Big[\int_0^{\delta\theta_L}r^\Ecal\big(\eta^\Ecal (u);z\big)\,du\Big].
\end{align}
On the other hand, by Markov's inequality, for any $\eta\in\Ecal$,
\begin{align}\label{markovin}
\P_\eta\Big[\sup_{s\le \delta\theta_L} &\int_0^{s}
\sum_{z\neq 0}\frac{z}{L}\Big(1-\frac{z}{L}\Big)\,r^\Ecal\big(\eta^\Ecal (u);z\big)
du>\frac{\epsilon}{2}\Big]\nonumber\\
&\le \frac{2}{\epsilon}\sum_{z\neq 0}\frac{z}{L}\Big(1-\frac{z}{L}\Big)\E_\eta\Big[\int_0^{\delta\theta_L}r^\Ecal\big(\eta^\Ecal (u);z\big)du\Big].
\end{align}
Hence, it suffices to show that
\be
\lim_{\delta\downarrow 0}\limsup_{L\to\infty}\sup_{\eta\in\Ecal}\sum_{z\neq 0}\frac{z}{L}\Big(1-\frac{z}{L}\Big)\,\E_\eta\Big[\int_0^{\delta\theta_L}r^\Ecal\big(\eta^\Ecal (u);z\big)\,du\Big]=0.
\label{suffic}
\ee
We proceed analogously to the proof of Proposition \ref{equilibration} given in Section \ref{sec:equi2}, using the bounds on the mixing time $\tm (\epsilon' )$ for the restricted process given in (\ref{trele}). We pick $\epsilon' =1/\theta_L$ and split the integral in the preceding display according to
\be
\int_0^{\delta\theta_L } r^\Ecal (\eta^\Ecal (u);z)\, du=\int_0^{\tm (\epsilon' )} r^\Ecal (\eta^\Ecal (u);z)\, du +\int_{\tm (\epsilon' )}^{\delta\theta_L} r^\Ecal (\eta^\Ecal (u);z)\, du\ .
\ee
We can estimate the first contribution to (\ref{suffic}) as
\be
\sum_{z\neq 0}\frac{z}{L}\Big(1-\frac{z}{L}\Big)\,\E_\eta\Big[\int_0^{\tm(\epsilon')}r^\Ecal\big(\eta^\Ecal (u);z\big)\,du\Big]\leq 
\frac{\tm(\epsilon')}{4}\sup_{\eta\in\Ecal}\sum_{z\neq 0}r^\Ecal\big(\eta;z\big).
\ee
By translation invariance $\sup_{\eta\in\Ecal}\sum_{z\neq 0}r^\Ecal\big(\eta ;z\big)=\sup_{\eta\in\Ecal^0}\sum_{z\neq 0}r^\Ecal\big(\eta,\Ecal^{z}\big)$, 
which vanishes as $L\to\infty$ by Lemmas \ref{unibound} and \ref{trelest}. 
For the second contribution in \eqref{suffic} we use the definition of the mixing time and the fact \eqref{rest_meas} that the invariant measure $\mu^x$ of the restricted process  
equals the invariant measure $\mu$ restricted to the well $\Ecal^x$, 
to get the upper bound
\begin{align*}
\sum_{z\neq 0}\frac{z}{L}\Big(1-\frac{z}{L}\Big) &\E_\eta\Big[\int_{\tm(\epsilon')}^{\delta\theta_L}r^\Ecal\big(\eta^\Ecal (u),z\big)\,du\Big]\\
&\leq \frac{\delta\theta_L\epsilon'}{4}\sup_{\eta\in\Ecal^0}\sum_{z\neq 0}r^\Ecal\big(\eta,\Ecal^{z}\big)+\delta\theta_L\sum_{z\neq 0}\frac{z}{L}\Big(1-\frac{z}{L}\Big)r^\Lambda (z).
\end{align*}
The first term vanishes as $L\to\infty$ with the choice of $\epsilon'$ above and Lemma \ref{unibound}, and for the second term \eqref{tosho} implies
\be
\delta\theta_L \sum_{z\neq 0} \frac{z}{L}\big(1-\frac{z}{L}\big) r^\Lambda (z)\to \frac{\delta}{z_cI_b (\rho-\rho_c)^{b+1}},
\ee
hence establishing that the family $\{\Q^L\}_{L\in\N}$ is tight.\CQFD

\subsection{Proof of Proposition \ref{martcon} -- Martingale convergence\label{sec:martingale}}

The rescaled position of the condensate $\big(Y_t^L:\,t\ge 0\big)$ 
is a random variable that takes values on the space $D\big([0,T];\ \T\big)$ of c\`adl\`ag paths on the torus $\T$, and as we proved in section \ref{tightness}, the family of the corresponding distributions $\{\Q^L\}_{L\in\N}$ is tight. In this section we prove 
Lemma \ref{chara} stated below, that we use in section \ref{proof_main} to show that in fact $\{\Q^L\}_{L\in\N}$ is convergent and to characterise its limit. 
We begin with a preliminary result.
\begin{lemma}
\label{fd}
Let $\Q$ be any subsequential weak limit of $\Q^L$, $U$ a bounded continuous function in $D\big([0,T];\ \T\big)$ and $f:\ \T\to\R$ a continuous function on the torus. Then, for any $t\ge 0$
\[
\int U(\omega) f(\omega_t)\ d\Q^L(\omega)\longrightarrow \int U(\omega) f(\omega_t)\ d\Q(\omega).
\]
\end{lemma}
{\bf Proof.} Consider the mapping $\Pi^{f,t}: D\big([0,T];\ \T\big)\to\R$ with $\Pi^{f,t}(\omega)=f(\omega_t)$. Then $\Pi^{f,t}$ is continuous at paths $\omega\in D\big([0,T];\ \T\big)$ that are continuous at $t$, i.e. $\omega_t=\omega_{t-}:=\lim_{s\uparrow t}\omega_s$. To see this, note that for any $u\in D\big([0,T];\ \T\big)$ and $\lambda_t\in[0,T]$ we have
\[
\left|\Pi^{f,t}(u)-\Pi^{f,t}(\omega)\right|\le|f(u_t)-f(\omega_{\lambda_t})|+|f(\omega_{\lambda_t})-f(\omega_t)|.
\]
Since $f$ is continuous on the torus and $\omega$ is continuous at $t$, the right hand side of the preceding display can be made arbitrarily small, provided we can control $|u_t-\omega_{\lambda_t}|$ and $|\lambda_t-t|$. If we choose $u$ sufficiently close to $\omega$ in the Skorokhod topology we may find a $\lambda_t$ that simultaneously makes these quantities suitably small. Since the set of discontinuities of the function
$\omega\mapsto U(\omega)f(\omega_t)$ is contained in the set of paths $\big\{\omega:\ \omega_t\neq \omega_{t-}\big\}$, the assertion of the Lemma will follow if we show that
\be
\Q\big[\omega_t\neq\omega_{t-}\big]=0.
\label{left}
\ee
For $t\ge 0$, and $\eps,\delta>0$ define the subset $J_t^{\delta,\eps}\subset D\big([0,T];\ \T\big)$ by
\[
J_t^{\delta,\eps}=\{\omega\in D\big([0,T];\ \T\big): \sup_{s\in(t-\delta,t)}d_\T(\omega_{t},\omega_{s})>\eps\}.
\] 
We have $\{\omega_t\neq\omega_{t-}\}=\cup_{\eps>0}\cap_{\delta>0}J_t^{\delta,\eps}$. Hence 
\[
\Q\big[\omega_t\neq\omega_{t-}\big]=0\Leftrightarrow \lim_{\delta\downarrow 0}\Q\big[J_t^{\delta,\eps}\big]=0,\ \forall \eps>0.
\]
Noting that  $J_t^{\delta,\eps}$ are open in the Skorokhod topology it would suffice to show that
\[
\lim_{\delta\downarrow 0}\limsup_{L\to\infty}\Q^L\big[J_t^{\delta,\eps}\big]=\lim_{\delta\downarrow 0}\limsup_{L\to\infty}\ 
\P\big[\sup_{s\in(t-\delta,t)}d_{\T}\big(Y_t^L,Y_s^L\big)>\eps \big]=0,\ \forall\delta>0.
\]
But this follows from tightness estimate in Section \ref{tightness}, in particular (\ref{doobin}) and (\ref{markovin}) and estimates after that. \CQFD\\

Note that for any Lipschitz function $f:\T\to\R$ we have that ${\cal L}^\T(f)$ is a continuous function on $\T$. Indeed,
using the elementary estimate
\[
|f(y+u)-f(y)-f(x+u)+f(x)|\le 2\,\text{Lip}(f) \big(d_\T(x,y)\wedge d_\T(u,0)\big)\ ,
\]
we have that
\[
|\Lcal^\T f(y)-\Lcal^\T f(x)| \le 2\,\text{Lip}(f)\int_T\left(1\wedge\frac{d_\T(x,y)}{d_\T(u,0)}\right) du\le 4\,\text{Lip}(f) d_\T(x,y)\left(2+\ln\frac{1}{d_\T(x,y)}\right).
\]
In particular, the mapping $\omega\mapsto\int_0^t \Lcal^\T f(\omega_s)ds$ is continuous in the Skorokhod topology, and combined with Lemma \ref{fd} we get the following:

\begin{lemma}\label{chara}
Let $\Q$ be any subsequential weak limit of $\Q^L$ and $f: \T\to\R$ a Lipschitz continuous function on the torus. Set
\[
M^f_t(\omega)=f(\omega_t)-f(\omega_0 )-\int_0^t\Lcal^\T f(\omega_s)\ ds.
\]
If $U$ is a bounded continuous function in $D\big([0,T];\ \T\big)$ and $t\ge 0$ we have
\[
\int U(\omega) M^f_t(\omega)\ d\Q^L(\omega)\longrightarrow \int U(\omega) M_t^f(\omega)\ d\Q(\omega).
\]
\end{lemma}

\subsection{Uniqueness for the martingale problem}\label{mapro}

In this subsection we prove uniqueness for the martingale problem associated with the operator $\Lcal^T$,
i.e. there exists a unique measure $\Q$ on the space $D(\R_+;\T)$ of c\`adl\`ag paths on $\T$, such that for the coordinate process $(\omega_s : s\ge 0)$ we have
\begin{align}
\begin{cases}
&\Q\big[\omega_0 =y_0\in\T\big]=1,\\
&f\big(\omega_t \big)-\int_0^t \Lcal^\T f\big(\omega_s \big)\, ds \quad\mbox{ is a } \Q-\mbox{martingale for all } f\in\text{Lip}(\T).
\end{cases}
\label{martproblem}
\end{align}
We may plug the test function $f_k(x)=e^{2\pi i kx}, x\in[0,1], k\in\Z$ in (\ref{martproblem}) to get that if $\Q$ solves (\ref{martproblem}), then for any $t\ge s\ge 0$ we have
\[
\E^\Q\big[e^{2\pi ik (\omega_t-\omega_s)}\,\big|\, {\cal F}_s\big]= e^{-\psi(k)(t-s)},\qquad \text{where } \psi(k)=H(b,\rho)\int_\T\frac{1-\cos(2\pi k y)}{d_\T(0,y)}\,dy.
\]
In particular, this shows that under $\Q$ the coordinate process has independent, time homogeneous increments. This determines the finite dimensional distributions of $(\omega_s: s\ge 0)$ and implies that such a $\Q$ is unique. To get a better insight on the limiting process, we may rewrite $\psi$ as
\[
\psi(k)=H(b,\rho)\int_{-1}^{1}\frac{1-\cos(2\pi k y)}{|y|}\,dy
\]
and note that if $\big(X_t : t\ge 0\big)$ is a symmetric L\'evy process on $\R$ with L\'evy measure
\[
\Pi(dy)=H(b,\rho)\frac{\mathbbm{1}\big\{|y|<1\big\}}{|y|}\, dy,
\] 
then $\Q$ is the distribution of the process $\big(X_t (\text{mod}\, 1): t\ge 0\big)$.
\vspace*{2mm}

\section{Estimates on the invariant measure\label{invmeas}}

\par  In this section we collect some auxiliary results of technical nature that we needed throughout this article. We will use $C$ to denote a constant (not always the same) that only depends on absolute constants. When $X$ is a random variable defined in a probability space $(\Omega,{\cal F},\nu)$ and $A\in{\cal F}$ we will write $\nu\big(X;A\big)$ as a shorthand for $\int_A X\ d\nu$. 
Recall also that $\tN=N-\rho_cL$, and $Z_{L,N}$ is the normalization in \eqref{ZLN}.
 \begin{proposition} 
 \label{NagaevRates}
Suppose $\delta>\rho_c$. Then, provided that $L$ is sufficiently large, we have
\begin{equation}
\sup_{N\ge \delta L}\Big|\frac{z_c\,\nu \big[S_L=N\big]}{L\tilde{N}^{-b}}-1\Big|\le\frac{C}{L}
\quad\mbox{and thus}\quad \sup_{N\ge \delta L} \Big| \frac{Z_{L,N}}{z_c ^{L-1} L\tilde N^{-b} }-1\Big|\leq\frac{C}{L}\ .
\end{equation}
\end{proposition}
{\bf Proof.} The proof essentially follows the argument in \cite{DDS} keeping track of the rate of convergence when $N$ is supercritical. For a given sequence $h_L$ we define the events $B_x =\{\eta_{x}\le h_L\}$, and $C_k= B_1^c\cap\cdots\cap B_k^{c}\cap B_{k+1}\cap\cdots\cap B_L$. Then,
\begin{align}
\nu\big[S_L=N\big]&=\sum_{k=0}^L {L \choose k} \nu\big[\{S_L=N\}\cap C_k\big].
\label{decomp}
\end{align}
Let us denote by $G$ the distribution function of the one-site marginal of the critical measure, i.e. $G(t)=\nu\big[\eta_x\le t\big]$ and let $\bar{G}=1-G$. If $h_L$ is chosen so that $L\bar{G}(h_L)\to 0\Leftrightarrow h_L\gg L^{\frac{1}{b-1}}$, then Lemmas 2.4 and 6.1 in \cite{DDS} together imply that
\[
\sum_{k=2}^L {L \choose k} \nu\big[\{S_L=N\}\cap C_k\big]\le C L\,\nu\big[\{S_L=N\}\cap C_1\big] \times \big(L\bar{G}(h_L)\big).
\]
We will choose $h_L=L^{\frac{2}{b-1}}$ so that the sum of the terms for $k\ge 2$ in (\ref{decomp}) is at most $O(L^{-1})$ times the term for $k=1$. To estimate the term for $k=0$ we may use Lemma 2.1 in \cite{DDS}. Precisely, if $\kappa=\max\{\frac{2}{b-1},\frac{1}{2}\}$, then $L^\kappa$ is a natural scale for $S_L-\rho_c L$ and
\begin{equation}
\nu\big[\ |S_{L}-\rho_c L|>x,\ M_{L}\le h_L\big]\le \nu\big[\ |S_{L}-\rho_c L|>x,\ M_{L}\le L^\kappa\big]\le Ce^{-{x}{L^{-\kappa}}}\qquad \forall x>0.
\label{DDSest}
\end{equation}
Since $\tilde{N}\ge (\delta-\rho_c)L$ and $\kappa<1$, for all sufficiently large $L$ we have $e^{-{\tilde{N}}L^{-\kappa}}\le L^{-1}$, hence
\[
L\ \nu\big[\{S_L=N\}\cap C_1\big]\le\nu\big[S_L=N\big]\le L\ \nu\big[\{S_L=N\}\cap C_1\big] \left(1+\frac{C}{L}\right).
\]
We now pick a sequence $\gamma_L=L^\gamma$, with $\kappa<\gamma<1$. Using (\ref{DDSest}) again, it is not hard to see that 
\begin{align}
\nu\big[\{S_L=N\}\cap C_1\big]&=\frac{1}{z_c}\sum_{k=0}^{N-h_L}(N-k)^{-b}\ \nu\big[S_{L-1}=k,\ M_{L-1}\le h_L\big]\nonumber\\
&=\frac{1}{z_c}\nu\big((N-S_{L-1})^{-b};\ \{S_{L-1}< N-h_L, \ M_{L-1}\le h_L\}\big)\nonumber\\
&=\frac{1}{z_c}\nu\big((N-S_{L-1})^{-b};\ K_L\big)\ \big(1+o(\frac{1}{L})\big),
\end{align}
where $K_L=\{\ |S_{L-1}-\rho_cL|\le \gamma_L,\ M_{L-1}\le h_L\}$. It also follows that
\begin{equation}
\frac{L}{z_c}\nu\big((N-S_{L-1})^{-b};\ K_L\big)\le\nu\big[S_L=N\big]\le \frac{L}{z_c}\nu\big((N-S_{L-1})^{-b};\ K_L\big) \ \big(1+\frac{C}{L}\big).
\label{Eest}
\end{equation}
If $|S_{L-1}-\rho_cL|\le\gamma_L$ we have the following pointwise inequality.
\[
0\le(N-S_{L-1})^{-b}-\tilde{N}^{-b}-b\tilde{N}^{-(b+1)}(S_{L-1}-\rho_cL)\le {C}{\tilde{N}^{-(b+2)}}(S_{L-1}-\rho_cL)^2.
\]
Integrating over $K_L$ and setting  $k_*=\nu\big(S_{L-1};\ K_L\big)=(L-1)\nu\big(\eta_1;\ K_L\big)$ we get
\begin{align}
0\le \nu\big((N-S_{L-1})^{-b}&;\ K_L\big)-\frac{\nu\big[K_L\big]}{\tilde{N}^{b}}-\frac{b\big(k_*-\rho_cL\,\nu\big[K_L\big]\big)}{\tilde{N}^{b+1}}\nonumber\\
&\le\frac{C}{\tilde{N}^{b+2}}\nu\big( (S_{L-1}-\rho_c L)^2\big)\le\frac{CL}{\tilde{N}^{b+2}}.
\label{taylor}
\end{align}
We now have to estimate how close $\nu\big[K_L\big]$ is to 1, and $k_*$ to $\rho_cL$. Using (\ref{DDSest}) again we have
\begin{align}
\nu\big[K_L^c\big]&=\nu\big[M_{L-1}>h_L\big]+\nu\big[\ |S_{L-1}-\rho_cL|> \gamma_L,\ M_{L-1}\le h_L\big]\nonumber\\
&\le (L-1)\nu [\eta_1 >h_L ]+\nu\big[\ |S_{L-1}-\rho_cL|> \gamma_L,\ M_{L-1}\le L^\kappa\big]\nonumber\\
&\le C\big(Lh_L^{1-b}+e^{-L^{\gamma-\kappa}}\big).
\label{tailK}
\end{align}
We also have
\begin{align}
\rho_c&=\nu\big(\eta_1\big)\ge\nu\big(\eta_1;\ K_L\big)\big(=k_*/(L-1)\big)\nonumber\\
&=\nu\big(\eta_1; M_{L-1}\le h_L\big)-\nu\big(\eta_1;\ |S_{L-1}-\rho_cL|>\gamma_L,\ M_{L-1}\le h_L\big)\nonumber\\
&\ge G(h_L)^{L-2}\nu\big(\eta_1; \eta_1\le h_L\big)-\nu\big(\eta_1;\ |S_{L-1}-\rho_cL|>\gamma_L,\ M_{L-1}\le L^\kappa\big)\nonumber\\
&\ge\big(1-L\bar{G}(h_L)\big)(\rho_c-Ch_L^{2-b})-\nu\big(\eta_x^2\big)^{\frac{1}{2}}e^{-L^{\gamma-\kappa}/2}\ge \rho_c\big(1-CLh_L^{1-b}\big).\nonumber
\end{align}
In the penultimate step we have used Cauchy-Schwarz, (\ref{DDSest}), and the elementary inequality $(1-x)^n\ge 1-nx$.
Hence,
\begin{equation}
0\le \rho_c(L-1)-k_*\le C\,L^2h_L^{1-b}.
\label{kstar}
\end{equation}
Plugging this estimate into (\ref{taylor}), the assertion follows from (\ref{Eest}), since $L^2h_L^{1-b}=1.$\hfill\CQFD\\

In the proof of Proposition \ref{maximum} we will need the following lemma, which holds in particular for subextensive $\tN$.
\begin{lemma}\label{parti}
Suppose $L\to\infty$ and $\tilde{N}\gg\sqrt{L\log L}$. Then,
\begin{equation}
{Z_{L,N}}=z_c^{L-1}L{\tilde{N}^{-b}}\,\big(1+o(1)\big).
\end{equation}
\end{lemma}
{\bf Proof.} This follows immediately from Theorem 2.1 in \cite{DDS} since
\[
Z_{L,N}=\sum_{\eta\in X_{L,N}}\prod_{x\in\Lambda}\frac{1}{g!(\eta_x)}=z_c^L\nu\big[S_L=N\big].
\CQFD
\]
Corollary 1 in \cite{armendarizetal08} states that for $b>3$ the fluctuations of the condensate size around $\tN$ are of order $\sqrt{L}$ and asymptotically normal. The next Proposition provides a conditional large deviation upper bound for the size of the condensate.
\begin{proposition}
\label{maximum}
Suppose $L^\gamma\le\alpha_L\le \tilde{N}$ for some $\gamma>\frac{1}{2}$. Then, for $L$ sufficiently large we have
\begin{align}
\label{max-est}
\mu\big[ M_L\le \tilde{N}-\alpha_L \big] \le  C {L}{\alpha_L^{1-b}}.
\end{align}
\end{proposition}
\noindent 
{\bf Proof.} 
We first observe that there exists a positive $\epsilon=\epsilon(b,\gamma)$ such that for $L$ sufficiently large we have
\be\label{inter2}
\mu \big[ M_L\le \epsilon \tN]\le CL\tN^{1-b}.
\ee
This follows from Lemma 5 in \cite{doney01} and Lemma \ref{parti} since
\[
\mu \big[M_L\le\epsilon \tN\big]=\frac{z_c^{L}}{Z_{L,N}}\nu\big[S_L=N,\ M_L\le \epsilon \tN\big]\le \frac{z_c^{L}}{Z_{L,N}}\big(\frac{CL}{\epsilon\tilde{N}^2}\big)^{\frac{1}{2\epsilon}}.
\]
Choosing $\epsilon<\frac{\gamma-1/2}{(2b-1)\gamma-1}$ we get that (\ref{inter2}) holds for $L$ large enough. With this observation the assertion follows easily when $\alpha_L\ge (1-\epsilon)\tN$. Indeed, in this case
\[
\mu \big[M_L\le \tN-\alpha_L\big]\le\mu \big[M_L\le\epsilon \tN\big]\le CL\tN^{1-b}\le CL\alpha_L^{1-b}.
\]
If on the other hand $\alpha_L\le (1-\epsilon)\tN$ we have
\begin{align}
\label{first_step}
\mu  \big[M_L\le \tN-\alpha_L \big] \le \mu \big[M_L\le \epsilon\tN \big]+L\mu \big[\epsilon\tN < M_L=\eta_1 \le \tN-\alpha_L\big].
\end{align}
Now,
\begin{align}
\label{long}
\mu \big[ \epsilon\tN\le M_L=\eta_1\le \tN-\alpha_L\big]&=\frac{z_c^{L-1}}{Z_{L,N}}\sum_{k=\epsilon\tilde{N}}^{\tN-\alpha_L} \frac{1}{k^b}\, \nu\big[M_{L-1}\le k,\, S_{L-1}=N-k\big] \notag\\
& \le\frac{z_c^{L-1}}{Z_{L,N}}\sum_{k=\epsilon\tN}^{\tN-\alpha_L} \frac{1}{k^b}\, \nu\big[S_{L-1}=N-k\big] \notag\\
&=\frac{(L-1)z_c^{L-2}}{Z_{L,N}}\big(1+o(1)\big)\sum_{k=\epsilon\tN}^{\tN-\alpha_L} \frac{1}{k^b}\frac{1}{(\tN-k)^b},
\end{align}
where the last line is a consequence of Theorem 2.1 in \cite{DDS} and the fact that $N-k\ge \rho_cL+\alpha_L$ if $k\le \tN-\alpha_L$. By Jensen's inequality for the convex function $x\mapsto {x^{-b}(1-x)^{-b}},\, x\in(0,1)$ we have
\[
\sum_{\epsilon\tN}^{\tN-\alpha_L} \frac{1}{k^b(\tN-k)^b}\le \tN^{1-2b}\int_{\epsilon-\frac{1}{2\tN}}^{1-\frac{\alpha_L-1/2}{\tN}}\frac{dx}{x^b(1-x)^b}\le C \tN^{-b}\alpha_L^{1-b},
\]
for $L$ sufficiently large. Applying this estimate to the last line of \eqref{long}, we conclude that
\begin{align}
\label{termtwo}
\mu \big[ \epsilon\tN\le M_L=\eta_1\le \tN-\alpha_L\big]\le C\frac{z_c^{L}L\tilde{N}^{-b}}{Z_{L,N}} {\alpha_L^{1-b}}.
\end{align}
In view of Lemma \ref{parti} the assertion follows by replacing \eqref{inter2} and \eqref{termtwo} in \eqref{first_step}.
\CQFD\\
\begin{proposition}\label{secmax}
Suppose $L^\gamma\le\alpha_L\le \frac{1}{2}\tilde{N}$ for some $\gamma>1/2$. Let $M_L^{(2)}$ stand for the second largest component of $(\eta_1,\ldots,\eta_L)$. Then for $L$ sufficiently large we have
\begin{equation}
\mu \big[M_L>\tilde{N}-\alpha_L,\, M_L^{(2)} >\beta_L\big]\le CL\beta_L^{1-b}.
\end{equation}
\end{proposition}
{\bf Proof.} We have
\begin{align*}
\mu \big[M_L>\tilde{N}&-\alpha_L,\, M_L^{(2)} >\beta_L\big]\le L \mu \big[
\eta_L>\tilde{N}-\alpha_L,\, \eta_L\ge M_{L-1}>\beta_L\big]\\
&=\frac{L\, z_c^{L-1}}{Z_{L,N}}\sum_{k>\tilde{N}-\alpha_L}\frac{1}{g!(k)}\nu\big[S_{L-1}=N-k,\ k\ge M_{L-1}>\beta_L\big]\\
&\le\frac{L\, z_c^{L-1}}{Z_{L,N} (\tilde{N}-\alpha_L)^b}\nu\big[S_{L-1}<\rho_cL+\alpha_L,\ M_{L-1}>\beta_L\big]\\
&\le\frac{L\, z_c^{L-1}}{Z_{L,N} (\tilde{N}-\alpha_L)^b}\nu\big[M_{L-1}>\beta_L\big].
\end{align*}
The assertion now follows from Lemma \ref{parti} and the elementary estimate
\[
\nu\big[M_{L-1}>\beta_L\big]\le (L-1)\bar{G}(\beta_L)\le CL\,\beta_L^{1-b}.\CQFD
\]

Propositions \ref{maximum} and \ref{secmax} show that the invariant measure $\mu $ is essentially supported in the union of the wells, as the following corollary states.
\begin{corollary}\label{shaft}
Consider the wells ${\cal E}^x$ defined in (\ref{valley}). Suppose $L^\gamma\le\alpha_L\le \frac{1}{2}\tilde{N}$ for some $\gamma>1/2$, and $\beta_L\gg L^{\frac{1}{b-1}}$. Then, for the complement $\Delta=X_{L,N}\setminus\bigcup_{x\in\Lambda}{\cal E}^x$ we have
\[
\mu \big[\Delta]\le CL\big(\alpha_L^{1-b}+\beta_L^{1-b}\big)\to 0\quad \text{as }L\to\infty.
\]
\end{corollary}
{\bf Proof.} It suffices to note that
\begin{align*}
\mu \big[\Delta\big]&=\mu \big[\{M_L\le \tilde{N}-\alpha_L\}\cup\{M_L^{(2)}>\beta_L\}\big]\\
&=\mu \big[M_L\le \tilde{N}-\alpha_L\big]+\mu \big[M_L>\tilde{N}-\alpha_L,\ M_L^{(2)}>\beta_L\big].
\end{align*}
\CQFD
\begin{proposition}
\label{uno}
Suppose $L^\gamma\le\alpha_L\le \frac{1}{2}\tilde{N}$ for some $\gamma>1/2$, and consider the sets ${\cal I}_N^{0,y}$ defined in \eqref{IN0y}. Then
\[
\mu \Big[\big\{\eta\in X_{L,N}\setminus \bigcup_{0\neq y\in \Lambda}{\cal I}_N^{0,y}:\, 
2\alpha_L\le {\eta_0}\le \tN-2\alpha_L \big\}\Big]\le \frac{C\ L}{\alpha_L^{2b-2}}.
\]
\end{proposition}
\noindent\textbf{Proof.} We have
\begin{align*}
&\mu \Big[\big\{\eta\in X_{L,N}\setminus \bigcup_{0\neq y\in \Lambda}{\cal I}_N^{0,y}:\, 
2\alpha_L\le {\eta_0}\le \tN-2\alpha_L \big\}\Big]\\
&\hspace{1.5cm}=\sum_{k=2\alpha_L}^{\tN-2\alpha_L} \frac{Z_{L-1,N-k}}{Z_{L,N}} \frac{1}{g!(k)}
\,\,\,\mu_{ L-1, N-k}\big[ M_{L-1} \le \tilde N-k-\alpha_L\big]\\
&\hspace{1.5cm}\le \frac{CL}{\a_L^{b-1}} \sum_{k=2\alpha_L}^{\tN-2\alpha_L} \frac{\tilde{N}^b}{k^b(\tN-k)^b}\\
&\hspace{1.5cm}\le\,  \frac{CL}{(\a_L\tN)^{b-1}}\int_{\frac{2\alpha_L-1}{\tN}}^{1-\frac{2\alpha_L-1}{\tN}}\frac{dx}{x^b(1-x)^b}\le\frac{C\ L}{\alpha_L^{2b-2}}.
\end{align*}
The third line follows from the second by Lemma \ref{parti} and an application of Proposition \ref{maximum}, term by term, with supercritical particle number $N-k\ge \rho_cL+2\alpha_L$, while the fourth line follows from the third by Jensen's inequality.
\CQFD\\

The following Proposition states that deviations of occupation numbers below their typical value have exponentially decaying probability.
\begin{proposition}
\label{largedev}
Suppose $\tN\ge L^\gamma$ for some $\gamma>\frac{1}{2},$ and let $A\subset\Lambda$. Let $\teta_x=\eta_x-\rho_c$ and $v=2\nu\big(\eta_x^2\big)$. For any $m>0$ and large enough $L$ we have
\begin{align*}
\mu\Big[ \sum_{x\in A} \teta_x\le -m\Big]\le Ce^{-\frac{m^2}{v|A|}}.
\end{align*}
\end{proposition}
\noindent\textbf{Proof.} 
Write
\begin{align*}
\mu \Big[ \sum_{x\in A} \teta_x \le -m\Big]&= \frac{z_{c}^{L}}{Z_{L,N}} \sum_{k\le \rho_c|A|-m} 
\nu\big[ S_{|A|}=k\big]\,\nu\big[S_{L-|A|} = N-k \big].\notag
\end{align*}
For $k$ in the range of summation the number of particles in $\Lambda\setminus A$ is supercritical. Indeed, $N-k-\rho_c(L-|A|)\ge\tN+m\ge\tN$. Hence, $\nu\big[S_{L-|A|} = N-k \big]\le {C(L-|A|)}{\tN^{-b}}$ for $L$ large enough, and
\[
\mu \Big[ \sum_{x\in A} \teta_x \le -m\Big]\le \frac{Cz_{c}^{L}}{Z_{L,N}}  \frac{(L-|A|)}{\tN^b} \sum_{k\le \rho_c|A|-m}\nu\big[ S_{|A|}=k\big]\ \le C\,\nu\Big[ \sum_{x\in A} \teta_x \le -m\Big],
\]
where the last inequality follows from Lemma \ref{parti}. The rest of the proof is now a standard Chernoff bound for the sum of the independent variables $\{\eta_x\}_{x\in A}$. For any $\lambda<0$ we have
\begin{align*}
\nu\Big[ \sum_{x\in A} \teta_x \le -m\Big]&\le e^{\lambda(m-\rho_c|A|)}\nu\big( e^{\lambda\eta_x}\big)^{|A|}
\le e^{\lambda(m-\rho_c|A|)}\bigg(\nu\Big( 1+\lambda\eta_x+\frac{\lambda^2}{2}\eta_x^2\Big)\bigg)^{|A|}
\\
&= e^{\lambda(m-\rho_c|A|)}\big(1+\lambda\rho_c+\frac{v\lambda^2}{4}\big)^{|A|}\le e^{\lambda m+\frac{v|A|\lambda^2}{4}}.
\end{align*}
The assertion now follows by choosing $\lambda=-2m/v|A|$.
\CQFD\\

\section*{Acknowledgements}
The authors are grateful to Martin Slowik and Claudio Landim for useful discussions. This work was supported by the Engineering and Physical Sciences Research Council (EPSRC) -- Grant No.\ EP/I014799/1, grant PICT 2012-2744 Stochastic Processes and Statistical Mechanics, the European Social Fund and Greek national funds through the Operational Programme "Education and Lifelong learning" -NSRF Research Funding Programmes Thales MIS377291 and Aristeia 1082.


\small

\end{document}